\documentclass[a4paper,11pt,reqno]{amsart}
\usepackage[utf8]{inputenc}
\usepackage{amsfonts}
\usepackage{amssymb}
\usepackage{amsmath}
\usepackage{amsthm}
\usepackage{pdflscape}
\usepackage{mathrsfs}
\usepackage{bigints}
\usepackage{subcaption}
\usepackage{etoolbox} 
\usepackage{graphicx}
\usepackage{caption}
\usepackage{pgfplots}
\pgfplotsset{compat=1.18}
\usepackage{xcolor}
\usepackage{algorithm}
\usepackage{algpseudocode}
\usepackage{hyperref}
\hypersetup{
colorlinks=true,
linkcolor=blue,
urlcolor=blue,
citecolor=blue
}

\usepackage{bm}
\usepackage{commath}
\usepackage{eurosym}
\usepackage{float}
\usepackage{dsfont}
\usepackage[nocompress, space]{cite}

\allowdisplaybreaks

\usepackage{ifthen}

\makeindex

\renewcommand{\div}{\operatorname{div}}

\newcommand{\diam}{\operatorname{diam}}

\def\d{{\rm d}}

\def\dt{{\rm d}t}

\numberwithin{equation}{section}

\newtheoremstyle{thmlemcorr}{10pt}{10pt}{\itshape}{}{\bfseries}{.}{10pt}{{\thmname{#1}\thmnumber{
#2}\thmnote{ (#3)}}}
\newtheoremstyle{thmlemcorr*}{10pt}{10pt}{\itshape}{}{\bfseries}{.}\newline{{\thmname{#1}\thmnumber{
#2}\thmnote{ (#3)}}}
\newtheoremstyle{defi}{10pt}{10pt}{\itshape}{}{\bfseries}{.}{10pt}{{\thmname{#1}\thmnumber{
#2}\thmnote{ (#3)}}}
\newtheoremstyle{remexample}{10pt}{10pt}{}{}{\bfseries}{.}{10pt}{{\thmname{#1}\thmnumber{
#2}\thmnote{ (#3)}}}
\newtheoremstyle{ass}{10pt}{10pt}{}{}{\bfseries}{.}{10pt}{{\thmname{#1}\thmnumber{
A#2}\thmnote{ (#3)}}}

\theoremstyle{thmlemcorr}
\newtheorem{theorem}{Theorem}
\numberwithin{theorem}{section}
\newtheorem{lemma}[theorem]{Lemma}

\newtheorem{proposition}[theorem]{Proposition}

\theoremstyle{thmlemcorr*}
\newtheorem{theorem*}{Theorem}
\newtheorem{lemma*}[theorem]{Lemma}
\newtheorem{corollary*}[theorem]{Corollary}
\newtheorem{proposition*}[theorem]{Proposition}
\newtheorem{problem*}[theorem]{Problem}
\newtheorem{conjecture*}[theorem]{Conjecture}

\theoremstyle{defi}
\newtheorem{definition}[theorem]{Definition}
\newtheorem{hyp}{Assumption}

\newtheorem{problem}{Problem}

\theoremstyle{remexample}
\newtheorem{remark}[theorem]{Remark}
\newtheorem{example}[theorem]{Example}

\theoremstyle{ass}

\theoremstyle{plain}

\newcommand{\R}{\mathbb{R}}
\newcommand{\T}{\mathbb{T}}

\newcommand{\Ptwo}{\mathcal{P}_2}

\title[MFG   Model with  Kantorovich Potential Coupling ]{Mean Field Game Problem with Kantorovich Potential Coupling: An Attrition model}

\author[A. Arakelyan]{Avetik Arakelyan}
\address{Institute of Mathematics, SNCO; and Yerevan State University;
and DiGREA lab, Center for Scientific Innovation and Education, Yerevan, Armenia}
\email{avetik.arakelyan@ysu.am}

\author[D. Alaverdyan]{Davit Alaverdyan}
\address{Center for Scientific Innovation and Education; and Yerevan State University, Yerevan, Armenia}
\email{davit.alaverdyan2@edu.ysu.am}

\author[N. Hovakimyan]{Naira Hovakimyan}
\address{Department of Mechanical Science and Engineering, University of Illinois Urbana-Champaign, Urbana, IL, USA}
\email{nhovakim@illinois.edu}

\author[I. Kaminer]{Isaac Kaminer}
\address{Department of Mechanical and Aerospace Engineering, Naval Postgraduate School, Monterey, CA, USA}
\email{kaminer@nps.edu}

\author[T. Bakaryan]{Tigran Bakaryan}
\address{Institute of Mathematics, SNCO; and Yerevan State University;
and DiGREA lab, Center for Scientific Innovation and Education, Yerevan, Armenia}
\email{tigran.bakaryan@instmath.sci.am}

\thanks{The research is supported by the Higher Education and Science Committee of the MESCS RA under Research Project No.~24IRF-1A001 and in part by AFOSR grant \#FA9550-21-1-0411 and NASA grant \#80NSSC22M0070.}

\begin{document}

\begin{abstract}
   In this work we address  a class of Mean-Field Game (MFG) systems on the flat torus $\mathbb{T}^{d}$ in which individual agent incentives arise from macroscopic, distribution-dependent objectives governed by the Wasserstein-2 ($W_{2}$) distance. Taking the first variation (Lions derivative) of the Wasserstein-based macroscopic functional yields a coupled forward-backward system of Hamilton–Jacobi–Bellman (HJB) and Fokker–Planck (FP) equations featuring a nonlocal coupling driven by the Kantorovich optimal transport potential. As a primary case study, we formulate a large-scale attacker-defender attrition model aimed at protecting High-Value Units (HVUs). By defining population-wise attrition rates via the $W_{2}$ metric and casting the defense strategy as a Pareto-optimal McKean–Vlasov control problem, we derive a concrete HJB–FP system with a Kantorovich potential coupling. We prove the local-in-time existence of classical solutions for this system  using Schauder's fixed-point theorem.  Finally, we compute the MFG dynamics by combining the Sinkhorn algorithm with ResNet-based Physics-Informed Neural Networks (PINNs), and   simulate two-dimensional attacker-defender engagements across diverse attrition parameters, control costs, and multi-HVU configurations.
\end{abstract}


\maketitle
\tableofcontents

\section{Introduction}

\subsection{Background and Organization of the Paper}

Mean-field games (MFGs) describe strategic interactions in large populations
of agents whose individual decisions depend on the evolving population
distribution \cite{lasry2007mean,carmona2018probabilistic}. For stochastic
dynamics with diffusion, an equilibrium is described by a backward
Hamilton--Jacobi--Bellman (HJB) equation for the value function and a forward
Fokker--Planck (FP) equation for the population density. The HJB equation
determines the optimal individual response, while the FP equation imposes
consistency between this response and the population evolution. In
applications such as swarm coordination and distributed defense, a central
modeling question is how to specify individual interaction costs when the
desired outcome is expressed in terms of the spatial distribution of the
whole population.

Optimal transport provides a natural framework for this purpose. The squared
Wasserstein distance $W_2^2(m,\nu)$ measures the minimum quadratic cost of
redistributing one probability measure into another. It therefore accounts
for the spatial displacement of mass and distinguishes distributions
according to both their relative position and their spatial spread. Under
appropriate assumptions, the first variation of this transport cost with
respect to $m$ is represented, up to an additive constant, by a Kantorovich
potential. Using this first variation as an individual interaction cost
connects a population-level objective to a spatially varying incentive, as
in potential mean-field games \cite{Graber2025}. The resulting HJB--FP
system has a nonlocal coupling. Evaluating the cost at a single position
requires information about the full population distributions.

In this paper, we study such couplings through an attacker--defender model
for the protection of one or more high-value units (HVUs). Existing
agent-based formulations describe attrition through individual trajectories
and pairwise interactions
\cite{walton2018optimal,tsatsanifos2021,walton2022defense,kang2025online}.
These models resolve individual engagements, but the number of interaction
terms and control variables increases with the sizes of the competing
populations. This motivates a continuum description in which attacker,
defender, and HVU configurations are represented by probability
distributions.

We introduce population-wise attrition rates that depend on the
Wasserstein distance between the relevant distributions. This formulation
keeps information about spatial spread, which can substantially affect
defensive effectiveness even when the populations have similar mean
positions; see \cite{arakelyan2026mfg}. The accumulated attrition determines
the survival probabilities of the attackers and the HVUs. A parameter
$\alpha$ balances the survival objective against a quadratic penalty on
defender control effort. Motivated by this population-level criterion, we
formulate a representative-agent MFG model for the defender population,
with prescribed attacker and HVU trajectories. Its coupling combines a
normalized Kantorovich potential with scalar factors depending on
Wasserstein distances and their time history.

Our main analytical result is the local-in-time existence of classical
solutions on the flat torus $\T^d$. The assumptions include a positive,
H\"older-continuous initial defender density, a prescribed attacker density
bounded away from zero and infinity with uniform spatial H\"older
regularity, and suitable temporal regularity and continuity conditions on
the data and terminal coupling. The periodic setting allows us to study
the interaction between optimal transport and the parabolic equations
without imposing spatial boundary conditions.

The main difficulty is to control the dependence of the Kantorovich
potential on the evolving densities. We combine spatial regularity
results \cite{caffarelli1992} with quantitative stability estimates for
optimal transport
\cite{delalande2023,gallouet2025,kitagawa2025}
to obtain the bounds and continuity properties required for the coupling.
The Hopf--Cole transformation and parabolic estimates then provide
gradient and Hessian bounds for the HJB solution, while comparison
estimates control the FP density from above and below. A key point is that
the estimates used to choose the spatial H\"older radius of the
admissible density class are independent of that radius. Together with
stability estimates for the coupled equations, these bounds yield a
continuous self-map of a compact, convex set of density trajectories.
Schauder's fixed-point theorem gives a solution for a sufficiently small
time horizon. We also discuss formation dynamics and optimal-control
flows as examples of prescribed attacker motion.

For the numerical approximation, we combine the Sinkhorn algorithm for
entropically regularized optimal transport \cite{cuturi2013sinkhorn} with
physics-informed neural networks (PINNs) \cite{pinn}. Separate residual
networks approximate the value function and the defender density.
A one-dimensional Gaussian test case provides an explicit reference
for assessing the approximation. Two-dimensional experiments then examine
the effects of relative attrition strengths, the control-cost parameter,
and the configuration of attacker groups and HVUs. They illustrate how
the defender population redistributes itself, including its division
into groups protecting different HVUs.

The paper is organized as follows. Subsection~\ref{subsec: OT-intro} recalls the optimal
transport notation and the Wasserstein distance. Section~\ref{sec: Motivation} explains the connection between distribution-dependent objectives and Kantorovich
potential couplings. Section~\ref{sec: Pop-wise-models} introduces the population-wise attrition rates and the finite-population model. Section~\ref{sec: MFG model} formulates the attacker--defender MFG system. Section~\ref{sec: existence} establishes the regularity estimates and proves local existence of classical solutions. Finally,
Section~\ref{sec: simulations} presents the explicit test case and the numerical experiments.

\subsection{Optimal Transport and the Wasserstein Distance}\label{subsec: OT-intro}

In this section, we briefly review the essential concepts of Optimal Transport (OT) theory that are utilized in our analysis. For a comprehensive treatment, we refer the reader to the monographs of Villani \cite{villani2009wasserstein, villani2021topics} and Santambrogio \cite{santambrogio2015optimal}.

Let the spatial domain be the flat torus $\mathbb{T}^d = \mathbb{R}^d / \mathbb{Z}^d$. Let $\mathcal{P}(\mathbb{T}^d)$ denote the space of probability measures on the torus. Given two probability measures $\mu, \nu \in \mathcal{P}(\mathbb{T}^d)$, the classical \textbf{Monge problem} seeks an optimal transport map $\mathcal{T}: \mathbb{T}^d \to \mathbb{T}^d$ that pushes $\mu$ forward to $\nu$, denoted as $\mathcal{T}_\# \mu = \nu$, meaning $\nu(A) = \mu(\mathcal{T}^{-1}(A))$ for all Borel sets $A$. The map $\mathcal{T}$ is chosen to minimize the total transport cost
\begin{equation}
    \inf_{\mathcal{T}_\# \mu = \nu} \int_{\mathbb{T}^d} c(x, \mathcal{T}(x)) \, d\mu(x),
\end{equation}
where $c(x, y)$ is typically the squared geodesic distance on the torus, $d_{\mathbb{T}^d}(x,y)^2$.

Due to the non-convex constraints of the Monge formulation, it is reasonable to consider the \textbf{Kantorovich relaxation}. We seek a transport plan $\gamma \in \mathcal{P}(\mathbb{T}^d \times \mathbb{T}^d)$ whose marginals are $\mu$ and $\nu$
\begin{equation} \label{eq:kantorovich_primal}
    \min_{\gamma \in \Pi(\mu, \nu)} \int_{\mathbb{T}^d \times \mathbb{T}^d} c(x, y) \, d\gamma(x, y),
\end{equation}
where $\Pi(\mu, \nu)$ is the set of couplings with marginals $\mu$ and $\nu$.

A central tool is \textbf{Kantorovich duality}. The minimum in \eqref{eq:kantorovich_primal} equals the maximum of the dual problem
\begin{equation}
    \sup_{(\varphi, \psi)} \left( \int_{\mathbb{T}^d} \varphi(x) \, d\mu(x) + \int_{\mathbb{T}^d} \psi(y) \, d\nu(y) \right),
\end{equation}
where the potentials satisfy $\varphi(x) + \psi(y) \le c(x, y)$. The optimal function $\varphi$ is called the \textbf{Kantorovich potential}. For convenience, we denote the Kantorovich potential from $\mu$ to $\nu$ by $\varphi_{\mu\to\nu}$. Similarly, the dual Kantorovich potential $\psi$ from $\nu$ to $\mu$ will be denoted by $\psi_{\nu\to\mu}$.

A fundamental result by Brenier \cite{brenier1991polar}, which establishes that the optimal transport map between sufficiently regular measures is uniquely characterized as the gradient of a convex function, was extended to the torus by Cordero-Erausquin \cite{cordero1999transport}. If $\mu$ is absolutely continuous with respect to the Lebesgue measure, the optimal transport map $\mathcal{T}$ exists and is unique. On the torus, it is given by
\begin{equation}
    \mathcal{T}(x) = x - \nabla \varphi(x),
\end{equation}
where $\varphi$ is a generalized convex function (or $c$-convex).

Finally, the minimal cost defines the $p$-Wasserstein distance as follows
\begin{equation}
    W_p(\mu, \nu) := \left( \inf_{\gamma \in \Pi(\mu, \nu)} \int_{\mathbb{T}^d \times \mathbb{T}^d} d_{\mathbb{T}^d}(x, y)^p \, d\gamma(x, y) \right)^{1/p}.
\end{equation}
This functional satisfies all the properties of a formal metric on the space of probability measures with finite $p$-th moments. In the context of Mean-Field Games, the $W_2$ metric is of particular importance, as it allows the Fokker--Planck equation to be interpreted as the gradient flow of the entropy functional. For completeness, we formally define the relevant space and metric below.

\begin{definition}[Second moment and $\mathcal{P}_2$ space]
Let $\mathcal{P}_2(\mathbb{T}^d)$ denote the space of all Borel probability measures on $\mathbb{T}^d$ with a finite second moment
\begin{equation}
    \mathcal{P}_2(\mathbb{T}^d) = \left\{ m \in \mathcal{P}(\mathbb{T}^d) : \int_{\mathbb{T}^d} d_{\mathbb{T}^d}(x, x_0)^2 \, dm(x) < +\infty \right\}.
\end{equation}
This space is metrized by the $2$-Wasserstein distance, defined for $m, m' \in \mathcal{P}_2(\mathbb{T}^d)$ by
\begin{equation}
    W_2(m, m') = \left( \inf_{\gamma \in \Pi(m, m')} \int_{\mathbb{T}^{2d}} d_{\mathbb{T}^d}(x, y)^2 \, d\gamma(x, y) \right)^{1/2},
\end{equation}
where $\Pi(m, m')$ denotes the set of all couplings of $m$ and $m'$ (i.e., probability measures on $\mathbb{T}^{2d}$ with marginals $m$ and $m'$).
\end{definition}

\section{Motivation}\label{sec: Motivation}

In many large-population systems, the desired outcome is specified
at the level of the population distribution rather than through
prescribed trajectories for individual agents. Examples include
steering a swarm toward a desired spatial configuration, allocating
mobile resources according to a distributed demand, and positioning
a defender population relative to an attacking population. In these
settings, individual agents incur movement or control costs, whereas
the effectiveness of their collective behavior depends on how the
entire population is distributed in space. This motivates models
that combine individual control objectives with macroscopic,
distribution-dependent criteria.

To illustrate this distinction, consider a representative agent
whose state evolves on $\mathbb{T}^d$ according to
\begin{equation}
\label{eq:motivation-dynamics}
    dX_t = \alpha_t\,dt + \sqrt{2}\,dW_t,
    \quad
    m_t = \operatorname{Law}(X_t),
    \quad
    \operatorname{Law}(X_0)=m_0,
\end{equation}
where $\alpha_t$ is the control and $W_t$ is a standard
$d$-dimensional Brownian motion. A population-level optimization
problem may take the form
\begin{equation}
\label{eq:macroscopic-objective}
\begin{aligned}
    \inf_{\alpha}\Bigg\{
    &\mathbb{E}\left[
        \int_0^T
        \left(
            \frac{1}{2}|\alpha_t|^2 + V(t,X_t)
        \right)\,dt
    \right]
    + \int_0^T \mathcal{F}(m_t,\nu_t)\,dt
    + \mathcal{G}(m_T)
    \Bigg\},
\end{aligned}
\end{equation}
subject to \eqref{eq:motivation-dynamics}.
Here, the quadratic term penalizes control effort,
$V$ describes a local environmental cost, and
$\mathcal{F}$ and $\mathcal{G}$ measure population-level
performance. The prescribed distribution $\nu_t$ may represent
a desired configuration or another population with which the
controlled population interacts.

The natural state variable for this centralized problem is the
probability distribution $m_t$. Consequently, dynamic programming
leads to a Bellman equation on an infinite-dimensional space of
probability measures, creating substantial analytical and
computational challenges; see \cite{PhamWei2017}.
In mean-field game theory, dependence of the equilibrium value
on arbitrary population states is described by the master equation.
For a fixed initial distribution and in the absence of common
noise, however, an equilibrium flow can be characterized directly
through a coupled Hamilton--Jacobi--Bellman and Fokker--Planck
system; see \cite{Cardaliaguet2016}.

Connecting a macroscopic objective to such a game requires
identifying the corresponding individual incentives. Simply adding
the scalar quantity $\mathcal{F}(m_t,\nu_t)$ to every agent's running
cost does not influence an individual's best response when the
population flow is held fixed. A spatially differentiated incentive
is instead provided by the first variation
\begin{equation}
\label{eq:marginal-coupling}
    f[m_t](x)
    =
    \frac{\delta \mathcal{F}}{\delta m}(m_t,\nu_t)(x),
\end{equation}
which measures the marginal effect of redistributing population
mass. More precisely, whenever the first variation exists, it
satisfies
\begin{equation}
\label{eq:first-variation}
\left.
\frac{d}{d\varepsilon}
\mathcal{F}\bigl((1-\varepsilon)m+\varepsilon\widetilde m,\nu\bigr)
\right|_{\varepsilon=0^+}
=
\int_{\mathbb{T}^d}
\frac{\delta\mathcal{F}}{\delta m}(m,\nu)(x)
\,d(\widetilde m-m)(x).
\end{equation}
This is the basic mechanism underlying potential mean-field games.
It connects population optimization with decentralized
decision-making without identifying arbitrary mean-field game
equilibria with socially optimal controls; see \cite{Graber2025}.

A particularly natural macroscopic objective is the quadratic
Wasserstein mismatch
\begin{equation}
\label{eq:wasserstein-objective}
    \mathcal{F}(m,\nu)
    =
    \frac{\kappa}{2}W_2^2(m,\nu),
    \quad \kappa>0,
\end{equation}
where
\begin{equation}
\label{eq:wasserstein-distance}
    W_2^2(m,\nu)
    =
    \inf_{\pi\in\Pi(m,\nu)}
    \int_{\mathbb{T}^d\times\mathbb{T}^d}
    d_{\mathbb{T}^d}(x,y)^2\,d\pi(x,y),
\end{equation}
and $\Pi(m,\nu)$ denotes the set of probability measures with
marginals $m$ and $\nu$.
Unlike a pointwise comparison of densities, this quantity accounts
for the spatial displacement required to rearrange one distribution
into another. It therefore distinguishes configurations according
to their geometry, without prescribing which individual agent
should occupy which target position.

Under assumptions ensuring uniqueness of the Kantorovich potential
up to additive constants, the first variation satisfies
\begin{equation}
\label{eq:wasserstein-first-variation}
    \frac{\delta}{\delta m}
    \left[\frac{1}{2}W_2^2(m,\nu)\right](x)
    =
    \phi_{m\to\nu}(x)+C,
\end{equation}
where $\phi_{m\to\nu}$ is the Kantorovich potential associated with
the transport cost
\[
    c(x,y)=d_{\mathbb{T}^d}(x,y)^2,
\]
and $C$ is independent of $x$.
After fixing an additive normalization, the resulting coupling is
\begin{equation}
\label{eq:kantorovich-coupling}
    f[m_t](x)=\kappa\phi_{m_t\to\nu_t}(x).
\end{equation}
Thus, optimal transport provides a direct link between a global
distribution-matching objective and a spatially varying individual
cost. The first-variation characterization
\eqref{eq:wasserstein-first-variation} follows from Kantorovich
duality; see \cite[Proposition~4.8]{Santambrogio2017}.

For quadratic control costs, the running objective
\eqref{eq:wasserstein-objective}, and the terminal functional
\[
    \mathcal{G}(m)
    =
    \int_{\mathbb{T}^d}\Psi(x)\,dm(x),
\]
the corresponding prototypical equilibrium system is
\begin{equation}\label{eq:prototype-kantorovich-mfg}
\begin{cases}
    -\partial_t u-\Delta u+\frac{1}{2}|D_xu|^2
    =
    V(t,x)+\kappa\phi_{m_t\to\nu_t}(x),
    \\
    \partial_t m-\Delta m
    -\operatorname{div}(mD_xu)
    =0,
    \\
    m(0,\cdot)=m_0,\;
    u(T,\cdot)=\Psi.
\end{cases}
\end{equation}
The equations hold on $(0,T)\times\mathbb{T}^d$.
The optimal feedback is
\begin{equation}
\label{eq:optimal-feedback}
    \alpha^*(t,x)=-D_xu(t,x),
\end{equation}
and the forward equation imposes consistency between individual
optimal behavior and the population distribution.

The location of the Wasserstein objective in
\eqref{eq:macroscopic-objective} determines where the Kantorovich
potential appears in the equilibrium system. A running Wasserstein
objective produces the coupling in the Hamilton--Jacobi--Bellman
equation. In contrast, a terminal objective of the form
\[
    \mathcal{G}(m)
    =
    \int_{\mathbb{T}^d}\Psi(x)\,dm(x)
    +
    \frac{\eta}{2}W_2^2(m,\nu_T),
    \quad \eta>0,
\]
produces the terminal condition
\begin{equation}
\label{eq:wasserstein-terminal-condition}
    u(T,x)
    =
    \Psi(x)+\eta\phi_{m_T\to\nu_T}(x),
\end{equation}
with a fixed normalization of the potential.

In fact we can define a general class of such MFG problems as follows
\begin{equation}\label{eq:gen-kantorovich-mfg}
\begin{cases}
    -\partial_t u-\Delta u+\frac{1}{2}|D_xu|^2
    = F(t,x,m,\nu,W_2(m,\nu),\phi_{m_t\to\nu_t})
    \\
    \partial_t m-\Delta m
    -\operatorname{div}(mD_xu)
    =0,
    \\
    m(0,\cdot)=m_0,\;
    u(T,\cdot)=\Psi,
\end{cases}
\end{equation}
where $F$ satisfies certain suitable conditions.

This class is particularly relevant to attacker--defender
models, where defensive effectiveness depends on the spatial
relationship between populations, not merely on their total sizes
or the positions of isolated agents. Distribution-dependent
interaction laws can therefore lead naturally to couplings
involving Kantorovich potentials. These couplings are nonlocal
and implicit: evaluating the cost at a single position requires
information about the full distributions and the solution of an
associated optimal transport problem.

The resulting analytical difficulty is to control how the
Kantorovich potential varies both in space and with the evolving
population densities.

Motivated by this issue, the present paper
establishes local-in-time existence of classical solutions for
the model under consideration. The argument combines regularity
and stability estimates for optimal transport potentials with
parabolic Schauder estimates, maximum-principle bounds for the
density, and a Schauder fixed-point construction.

\section{Population-wise Models}\label{sec: Pop-wise-models}
In this section, we introduce a population-wise attrition rate function that captures macroscopic interactions at the level of distributions. Based on this we discuss a large-scale $N, M$-player attacker-defender model.

\subsection{Population-wise Attrition Rate Function}\label{section 3}

To describe interactions between large populations of agents, we introduce a population-wise attrition rate function. This function captures the effect of one population on another at the distribution level.

Suppose that we are given two populations: a shooter population and a target population. We denote their spatial distributions by $\sigma_s$ and $\sigma_t$, respectively.
Let $R(\sigma_s,\sigma_t)$ be a statistical distance between the two population distributions. We define the population-wise attrition rate as
\begin{equation*}
    d^{att}(\sigma_s, \sigma_t) \equiv d^{att}(R(\sigma_s, \sigma_t)).
\end{equation*}
This formulation ensures that the attrition rate depends solely on the statistical distance between the populations, reflecting the premise that attack effectiveness is primarily dictated by the relative spatial configuration and proximity of the involved agents.

For example, one may choose an exponentially decaying attrition function of the form
\begin{equation*}
    d^{att}\left(R(\sigma_s, \sigma_t)\right) = \lambda \exp\left(-\frac{R^2(\sigma_s, \sigma_t)}{\sigma}\right),
\end{equation*}
where $\lambda>0$ represents the maximal attrition strength, while $\sigma>0$ controls how fast the attrition decreases as the population distance increases.
This attrition rate attains its maximum value when the statistical distance is zero and decreases to $0$ as the distance between populations increases.
Compared with classical agent-wise attrition models (see \cite{tsatsanifos2021}), where interactions are usually described between individual agents, the present formulation works directly with population distributions. This makes it more suitable for large-scale systems, where tracking all pairwise interactions can be computationally expensive.

The remaining modeling choice is the selection of the statistical distance $R(\cdot, \cdot)$.
In this work, we use the Wasserstein-2 distance.
A detailed comparison with other possible metrics is given in \cite{arakelyan2026mfg}.
As demonstrated in that analysis, while most distance metrics successfully capture the straightforward scenario of populations moving closer to one another, the Wasserstein-2 metric is suited to capture the complex, non-trivial effects of distribution variance.
Specifically, the effectiveness of a shooter population depends heavily on the spatial spread of populations, high concentration may fail to adequately cover a target population, moderate dispersion optimizes effectiveness, and excessive dispersion dilutes the weapon's impact.

\subsection{N, M - Player Model}

In this section, we  consider a population of $N$ attackers and $M$ defenders. 
Our focus is on scenarios where both populations are large, so that the individual influence of each agent becomes negligible.
To describe how agents interact through their surrounding population, we consider the empirical distributions
\begin{equation}\label{emp:distrib}
    \mu^N(t,x) =\frac{1}{N}\sum_{i=1}^{N}\delta_{X_t^i}(t,x),\quad m^M(t,x) =\frac{1}{M}\sum_{j=1}^{M}\delta_{S_t^j}(t,x),
\end{equation}
of attackers and defenders at time $t$.
These measures encode the mean-field effect—each agent responds not to individual others, but to the aggregate distribution of its own population and possibly to that of the opposing population.

We assume that the dynamics of the attackers and defenders are given by the following stochastic differential equations (SDE)
\begin{equation}\label{main_sde}
    \begin{split}
       \begin{cases}
        \d X_t^i = f_a(X_t^i, \mu_{X_t}^N, u_i^{att})\dt + \sigma\d W^i_t,& \quad i=1,\dots,N\\
        \d S^k_t = f_d(S^k_t, m_{S_t}^M, u_k) \dt +\sigma\d W^k_t, &\quad k=1,\dots,M,
    \end{cases}
    \end{split}
\end{equation}
where $f_a$ and $f_d$ stand for the dynamics of the defenders and attackers, respectively.
The \(X_t^i\) is the state of attacker \(i\), and \(S_t^k\) is the state of defender 
\(k\).
Here $W^j_t$ are independent Brownian motions. 

We assume that the dynamics of the attackers are known; that is, the controls $u_i^{att}$, $i=1,\dots,N$ are given.
Suppose we have several HVU's  that must be defended. Then following \eqref{emp:distrib} we define empirical distribution for those HVU's as follows
\[
 \nu^{HVU}(t,x) =\frac{1}{K}\sum_{i=1}^{K}\delta_{S_{HVU}^i}(t,x),
\]
where $K$ is a  number of HVU's.   

\begin{remark}
The naive approach is to apply the established model from \cite{walton2022defense}. However, in this model attrition dynamics  become computationally expensive as the number of agents (attackers and defenders) grows. The situation is even worse when we increase the number of HVU's as well.
Therefore, a new approach is needed to handle large-scale scenarios, which we propose in this work.
\end{remark}

Here along with an attrition rate distance defined over the full populations of attackers and defenders, we adapt the model from \cite{kang2025online}. To this aim, we let the survival probability of the attacker population from the   defender population over $[t,\; t+\Delta t]$ be 
$$
1 - d^{att}(R(m^M(t,\cdot), \mu^N(t,\cdot)))\Delta t.
$$
Then we let $Q(t)$ to be a survival probability of the population of attackers over $[0, t],$ i.e. $Q(t)$ is a total survival probability of $N$ attackers from $M$ defenders over $[0, t].$ Following the original model in \cite{kang2025online} we set
\[
\dot{Q}(t) = -Q(t)\cdot d^{att}\left(R(m^M(t,\cdot), \mu^N(t,\cdot))\right).
\]
Note that
$$
Q(t) = e^{-\int_0^t d^{att}(R(m^M(\tau,\cdot), \mu^N(\tau,\cdot)))d\tau}.
$$

We assume that the defenders are protecting a so-called high-value unit (HVU). 
Let $P(t)$ be a survival probability of HVU from the entire population of the attackers over $[0, t].$  We  have
\begin{equation*}
    { \dot{P}(t) = -P(t)\cdot d^{HVU}\left(R(\nu^{HVU}(t,\cdot), \mu^N(t,\cdot))\right)\cdot Q(t),}
\end{equation*}
where $Q(t)$ is a  survival probability of the population of attackers over [0, t] defined above section.
Hence
$$
P(t) = e^{-\int_0^t d^{HVU}(R(\nu^{HVU}(\tau,\cdot), \mu^N(\tau,\cdot)))Q(\tau)d\tau}.
$$
Recalling the explicit form of $Q(t),$  we have
$$
P(t) = e^{-\bigints_0^t d^{HVU}(R(\nu^{HVU}(\tau,\cdot), \mu^N(\tau,\cdot)))\cdot e^{-\bigintsss_0^\tau d^{att}(R(m^M(\theta,\cdot), \mu^N(\theta,\cdot)))d\theta}d\tau}.
$$
Let $\nu_{HVU}$ be a given probability distribution function, and HVU is described by $X_t^H\sim \nu_{HVU}$. The attackers and defenders dynamics are given by the system \eqref{main_sde}. Each agent  seeks to minimize
\begin{equation}\label{main_minimizer}
\min_{u_{k}} (1-P(t_{f})).
\end{equation}
Using the convexity property of exponent function we can rewrite
$$
\underset{u_k}{\min}\;(1-P(t_f))= \underset{u_k}{\max}\;(P(t_f))=\underset{u_k}{\max}\;(\log P(t_f))= \underset{u_k}{\min}\;(-\log P(t_f)).
$$
In order to involve explicitly the control $u_k$ of defenders, for some fixed $\alpha_k\in[0,1]$ we consider the following Pareto optimality version of minimizer 

\begin{equation}\label{pareto-v}
    \underset{u_k}{\min}\;\left((1-\alpha_k)(-\log P(t_f))+ \alpha_k\int_0^{t_f} u_k^2d\tau\right).
\end{equation}
Thus, recalling $P(t_f)$, the minimizer \eqref{pareto-v} reads as follows

\begin{equation}\label{eq:conv-cost}
\begin{split}
\underset{u_k}{\min}\; \bigints_0^{t_f}\left[(1-\alpha_k)\cdot d^{HVU}(R(\nu^{HVU}(\tau,\cdot), \mu^N(\tau,\cdot)))\cdot e^{-\bigintsss_0^\tau d^{att}(R(m^M(\theta,\cdot), \mu^N(\theta,\cdot)))d\theta}+\alpha_k u_k^2\right]d\tau.
\end{split}
\end{equation}
When the number of agents is large the minimization problem \eqref{pareto-v} very difficult to solve. Hence, in the next section, we consider mean-field limit to approximate the solution to the problem.

\section{Mean-Field Game Model}\label{sec: MFG model}

In this section, leveraging the population-wise attrition framework we define the macroscopic attacker-defender engagement. By formulating the optimal defense strategy as a stochastic control problem we consider corresponding MFG system.

\subsection{Mean-Field Control Problem}

Following the logic of the agent-wise models discussed earlier, we now pose the high-value unit (HVU) protection problem at the population level.
We consider populations of defenders, attackers, and  HVUs, tracked through the states of representative agents $S_t, X_t$, and $Y_t$, respectively.
To ensure well-posedness in our subsequent mathematical analysis, we define the spatial domain as a boundary-free torus, that is  $X_t, S_t, Y_t \in \mathbb{T}^d$.
Accordingly, let $m(t)$ denote the distribution of the defenders ($S_t \sim m(t)$), $\mu(t)$ the distribution of the attackers ($X_t \sim \mu(t)$), and $\nu_{\mathrm{H}}$ the  distribution of the HVUs ($Y_t = Y \sim \nu_{\mathrm{H}}(t)$). We also recall that as a statistic distance $R(\cdot,\cdot)$ for populations we take Wasserstein-2 metric.

The state evolution of the  defender populations is governed by the following  stochastic differential equation
\begin{equation}\label{McKean-Vlasov}
\mathrm{d}S_t = u\big(t, S_t\big)\,\mathrm{d}t + \sqrt{2\varepsilon}\,\mathrm{d}W_t^S,\\
\end{equation}
where $W_t^S$ is independent standard Brownian motion.
We assume the attackers drift field $u_{\mathrm{att}}$ is known and purely driven by the objective of destroying the HVU. For specific definitions of this drift field, we refer the reader to Examples \ref{ex:formation} and \ref{ex:op-tr}.
The defender's control strategy, $u$, must be optimally selected to maximize the HVU's survival probability.

In the rest of the paper, for any time-dependent distribution $\sigma(t,x)$, we will often use the shorthand notation $\sigma(t)$, omitting the spatial variable $x$ unless explicit spatial dependence needs to be emphasized.

Let $Q(t)$ represent the survival probability of the attacker population over the interval $[0,t]$.

Assuming the survival rate decays proportionally to the population-wise attrition rate $d^{\mathrm{att}}(\mu(t), m(t))$, we have
\begin{equation} \label{eq-Q-exp}
Q(t) = \exp\!\left(-\int_0^t d^{\mathrm{att}}(\mu(\tau), m(\tau))\,\mathrm{d}\tau\right).
\end{equation}
Similarly, let $P(t)$ be the survival probability of the HVU, governed by the attrition rate $d^{H}(\nu_{\mathrm{H}}, \mu(t))$ inflicted by the surviving attackers. Then
\begin{equation}\label{eq:-P-exp}
    P(t) = \exp\!\left(-\int_0^t d^{H}(\nu_{\mathrm{H}}(\tau), \mu(\tau))\, Q(\tau)\, \mathrm{d}\tau\right).
\end{equation}
By substituting \eqref{eq-Q-exp} into \eqref{eq:-P-exp}, the nested dependency of the HVU's survival on both the attackers' proximity and the defenders' continuous interception is fully characterized.

The primary objective of the defender population is to minimize the probability of the HVU's destruction, $\mathbb E\!\left[1 - P(T)\right]$. 
Given the strict monotonicity of the exponential function, this is mathematically equivalent to minimizing the negative logarithm of the survival probability
\begin{equation}
    \min_{u} \; \mathbb E\!\left[1 - P(T)\right] = \min_{u} \; \mathbb E\!\left[-\log P(T)\right].
\end{equation}
 We penalize the control effort with a weight $\alpha \in (0,1)$.
Combining this penalty with the simplified objective yields the following optimization problem: find the control $u$ that minimizes
\begin{equation}\label{minimizer_revised-v}
 \mathbb E\!\left[\int_0^{T} \left( (1-\alpha)d^{H}(\nu_{\mathrm{H}}(\tau), \mu(\tau)) Q(\tau) + \alpha u^2 \right) \mathrm{d}\tau\right],
\end{equation}
subject to the dynamics in \eqref{McKean-Vlasov}.

This formulation represents a McKean-Vlasov control problem because the running cost is intimately coupled with the macroscopic spatial distribution $m(t,\cdot)$.
Under standard dynamic programming, the resulting value function $U(t,x,m)$ leads to an infinite-dimensional master equation on the Wasserstein space $\mathcal{P}_2(\mathbb{T}^d)$.
Because directly solving the master equation is notoriously delicate and computationally prohibitive, we instead apply a representative-agent reduction in the following section.
This approach yields a decentralized mean-field game (MFG) system of coupled PDEs that is both analytically tractable and amenable to numerical approximation.


\subsection{Attacker-Defender Mean-Field Game Formulation}
\label{sec:attacker_defender_mfg}

Building on the agent-wise models discussed above, we define our MFG framework by evaluating the population distributions at any given time $t$. We represent the density-dependent attrition rate as $d^{att}(\mu(t), m(t))$, where $\mu$ and $m$ are the probability density functions of the attackers and defenders, respectively.

Let $Q(t)$ be the survival probability of the attackers over $[0, t]$. The evolution is defined by $\dot{Q}(t) = -Q(t)d^{att}(\mu(t), m(t))$, which yields:
\begin{equation}
Q(t) = e^{-\int_0^t d^{att}(\mu(\tau), m(\tau))d\tau}.
\end{equation}

Similarly, let $P(t)$ be the survival probability of the HVU facing attackers. We define the HVU attrition function as $d^{HVU}( \nu(t),\mu(t))$, where $\nu(t)$ is the distribution of the HVU. Evolving according to $\dot{P}(t) = -P(t)d^{HVU}(\nu(t),\mu(t))Q(t)$, we obtain
\begin{equation}
P(t) = \exp\left(-\bigints_0^t d^{HVU}(\nu(\tau), \mu(\tau))\cdot e^{-\bigintsss_0^\tau d^{att}(\mu(\theta), m(\theta))d\theta}d\tau\right).
\end{equation}
Our mean-field game minimization problem is framed as
\begin{equation}\label{mfg-mzr}
\mbox{Minimize}\; \int_t^{t_f}L(\tau, x, \nu_{HVU},\mu, m, u) d\tau,
\end{equation}
subject to
\begin{equation*}
        m_t -\varepsilon\Delta m-\nabla\cdot(u\cdot m)=0,\\
\end{equation*}
where the Lagrangian incorporates a Pareto optimality constant $\alpha\in(0,1)$ and
\begin{multline*}
    L(t, x, \nu_{HVU},\mu, m, u) = \\ = (1-\alpha)\cdot  d^{HVU}\left(W_2(\nu_{HVU}(t), \mu(t))\right)\cdot \frac{\delta}{\delta m}e^{-\bigintsss_0^t d^{att}(W_2(\mu(\theta), m(\theta))d\theta} + \alpha u^2.
\end{multline*}
The corresponding Hamiltonian, $H(t, x, \nu_{HVU}, \mu, m, p) = \max_q(-L - q\cdot p)$, is given by
\begin{multline}\label{mfg: Hamilton}
        H(t, x, \nu_{HVU}, \mu, m, p) = \\ = \frac{p^2}{4\alpha}-(1-\alpha)\cdot  d^{HVU}\left(W_2(\nu_{HVU}(t), \mu(t))\right)\cdot \frac{\delta}{\delta m}e^{-\bigintsss_0^t d^{att}(W_2(\mu(\theta), m(\theta)))d\theta}.
\end{multline}

Applying the chain rule to the exponential term requires the Lions derivative of the squared Wasserstein distance. The Lions derivative  of the squared Wasserstein distance $W_2^2(\mu_1, \mu_2)$ with respect to the measure $\mu_2$ is precisely twice the Kantorovich potential, $2\varphi_{\mu_2\to\mu_1},$ see \cite{carmona2018probabilistic}. In other words
 $$
 \frac{\delta}{\delta\mu_2}W_2^2(\mu_1, \mu_2)=2\varphi_{\mu_2\to\mu_1}.
 $$
We define
\[
    F(t,x, \nu_{HVU},\mu, m)  = - D(t)\cdot E(t, m,\mu)(d^{att})'(W_2(\mu(t), m(t)))\frac{\delta}{\delta m}[W_2^2(\mu, m)],
\]
where $D(t) = (1-\alpha)\cdot d^{HVU}\left(W_2(\nu_{HVU}(t), \mu(t))\right)$ and $E(t,m,\mu) = e^{-\int_0^t d^{att}(W_2(\mu(\theta), m(\theta))d\theta}$.

Our comprehensive attacker-defender model is encapsulated in the following problem:

\begin{problem}\label{prob-MFG}
Given $\nu_{HVU},$  let $\mu$ be the attackers dynamics. Find $(w,m)$ such that it solves the following coupled system of PDEs on $(0,T)\times \T^d$:
\begin{equation}\label{eq:MFG}
       \begin{cases}
        -w_t-\varepsilon\Delta w + \frac{|Dw|^2}{4\alpha}=F(t,x, \nu_{HVU},\mu, m)\\
         m_t -\varepsilon\Delta m-\frac{1}{2\alpha}\nabla\cdot(Dw\cdot m)=0\\
         w(T,x) = G(x,m(T)),\; m(0,x)=m_0,  \; x\in\T^d.
\end{cases}
\end{equation}
\end{problem}

\section{Existence of a Classical Solution}\label{sec: existence}
This section establishes the local existence of a classical solution to our Mean Field Game system \eqref{eq:MFG} via Schauder's Fixed-Point theorem. We first apply Caffarelli's regularity theory and stability estimates to show the time-Hölder and spatial smoothness of the populations distribution. Next, using a Hopf-Cole transformation, Bernstein's method, and Grönwall inequalities, we derive  uniformly bounds for the MFG system. This enables us to construct a continuous fixed-point mapping on a compact, convex space of measure trajectories over $\mathbb{T}^d$, providing a classical solution for a sufficiently small time horizon $T > 0$.

We start by writing down the necessary assumptions that make this analysis possible.

\subsection{Assumptions}\label{sec:mean_field_limits}
Here we outline the assumptions required to establish our main result and provide illustrative examples to demonstrate their applicability. These hypotheses are carefully chosen to obtain the necessary  well-posedness of the system.

We begin with the assumptions regarding the attrition functions, focusing on their smoothness.

\begin{hyp}\label{ass:attr}
 Let  $  d^{HVU}:\R\to \R$ be a non-negative Lipschitz continuous function, $d^{att}:\R\to \R$ be a bounded non-negative $C^2$ function.
\end{hyp}

 For instance, one may take $d^{HVU}(x) = \sqrt{1+x^2},$ which is Lipschitz and $d^{att}(x) = e^{-x^2}$ is bounded in $C^2$.

\begin{hyp}\label{ass:att_HVU_densities}
For every $t\in[0,T]$ let $\nu_{HVU}(t,\cdot)$ and $\mu(t,\cdot)$ are in $\Ptwo(\T^d)$ and $1/2$ Holder  continuous in time; that is, satisfying
$$
W_2\left(\mu(t,\cdot),\mu(s,\cdot)\right)\le C |t-s|^{\frac{1}{2}} \;\;\mbox{and}\;\; W_2\left(\nu_{HVU}(t,\cdot),\nu_{HVU}(s,\cdot)\right)\le C |t-s|^{\frac{1}{2}},
$$
for all $t,s \in[0,T]$ and some $C>0.$
\end{hyp}

\begin{hyp}\label{ass:term_cost}
For the terminal coupling function $G(x,\eta)$ in MFG system \eqref{eq:MFG} we assume
\[
\sup_{\eta\in\mathcal P(\T^d)}
\|G(\cdot,\eta)\|_{C^{2+1/(2d+4)}(\T^d)}
\le M_G<\infty,
\]
and
\[
W_2(\eta_n,\eta)\longrightarrow0
\quad\Longrightarrow\quad
G(\cdot,\eta_n)\longrightarrow G(\cdot,\eta)
\quad\text{in }C^2(\T^d).
\]
\end{hyp}

\begin{remark}
    Note that the solutions to  both Fokker-Planck equations from Examples \ref{ex:formation} and \ref{ex:op-tr} are satisfying Assumption \ref{ass:att_HVU_densities}.
\end{remark}

\begin{definition}\label{ass: Caffarelli-cond}
We say a time-dependent measure $\sigma(t,x)$ is  a Caffarelli measure if, for given constants $0<\lambda<\Lambda<\infty$, $K>0$, and $\gamma\in(0,1)$, it satisfies
\[
\lambda \le \sigma(t,x) \le \Lambda \quad \text{and} \quad \sup_{t\in[0,T]} \|\sigma(t,\cdot)\|_{C^\gamma(\T^d)} \le K,
\]
for all $(t,x)\in[0,T]\times \T^d.$
\end{definition}

\begin{hyp}\label{ass:att_Caffarelli}
We assume that the attacker distribution $\mu$ is a Caffarelli measure, i.e. satisfies the Definition \ref{ass: Caffarelli-cond}.
\end{hyp}

In the next two examples we discuss dynamics, which satisfy the above Assumptions \ref{ass:att_HVU_densities} and \ref{ass:att_Caffarelli}.

\begin{example}[Formation Dynamics]\label{ex:formation}
This example is motivated by well known formation dynamics in the finite number of agent setting. We start by presenting the finite case and pass to the limit to obtain corresponding mean-field equation. 

Consider a population of $N$ agents.
Let the position of the $i$-th agent at time $t$ be ${X}^i_t \in \T^d$, where $d$ is the spatial dimension. By following \cite{jabin2017mean, carrillo2010particle, topaz2004swarming}, each agent is governed by
\begin{equation}\label{gen_sde}
    d{X}^i_t = \left( J_i({X}^i_t) + R_i (\mathbf{X}_t)\right) dt + \xi d{W}^i_t,
\end{equation}
where  $\mathbf{X}_t=({X}^1_t,\dots,{X}^N_t)$, $ J_i$ steers the population towards a goal or along a trajectory, $R_i$ models agent interactions such as cohesion and collision avoidance, and finally $\xi d{W}^i_t$ represents stochastic noise, where $\xi$ controlling the amplitude and ${W}^i_t$ is a Wiener process.
The steering force guides agents towards a target by 
\begin{equation*}
    J_i(y) = -\nabla V^\text{tar}(\mathbf{y}),
\end{equation*}
with $V^\text{tar}$ as the potential function defined as follows
\begin{equation*}
    V^\text{tar}(\mathbf{y} ) = \frac{{b}}{2} \|\mathbf{y} - \mathbf{x}_\text{goal}\|^2.
\end{equation*}
The interaction force $R_i$ consist of cohesion and avoidance forces denoted by $S_i$ and $T_i$, respectively. Thus, the total interaction force will be
\begin{equation*}
   R_i (\mathbf{X}_t) = S_i(\mathbf{X}_t) + T_i(\mathbf{X}_t).
\end{equation*}
The cohesion force maintains proximity among agents and defined as follows
\begin{equation*}
        S_i(\mathbf{X}_t) = -\sum_{j \neq i} a_{ij} ({X}_t^i - {X}_t^j),
\end{equation*}
where the interaction weights are defined by a continuous function $a:\R^d\to \R$ scaled by the population size, i.e.  $a_{ij} = \frac{1}{N} a({X}_t^i - {X}_t^j)$.

The  avoidance force is preventing collisions and given as follows
    \begin{equation*}
        T_i(\mathbf{X}_t) = \sum_{j \neq i}Q({X}_t^i - {X}_t^j),
    \end{equation*}
where $Q(\cdot)$ is an avoidance kernel. The best known example of such avoidance kernel is a well-known Poisson kernel; that is,
$$
Q(\mathbf{x} ) =  C_d \frac{\mathbf{x} }{|\mathbf{x} |^\beta},\quad d = 2, 3, \dots
$$
where $\beta< d$ and $C_d > 0$   is a constant depending on the dimension and the physical parameters.

As the number of agents $N \to \infty$, the discrete system of interacting agents transitions to a mean-field description, where the individual agent dynamics are replaced by the evolution of a probability density function $\rho(t, \mathbf{x})$ representing the density of agents at position $\mathbf{x} \in \T^d$ at time $t$. 
\begin{remark}
Note that as $N \to \infty$ we have
     \begin{equation}
         S_i(\mathbf{X}_t)=-\sum_{j \neq i} a_{ij} ({X}_t^i - {X}_t^j) \to -\int_{\T^d} a(\mathbf{x} - \mathbf{y}) (\mathbf{x} - \mathbf{y}) \rho(t,\mathbf{y}) \, d\mathbf{y},
    \end{equation}
    where $a(\mathbf{x} - \mathbf{y})$ encodes the strength of interaction.
 Similarly, 
    \begin{equation}
         T_i(\mathbf{X}_t)=\sum_{j \neq i} Q({X}_t^i - {X}_t^j)\to \int_{\T^d} Q(\mathbf{x} - \mathbf{y})\rho(t,\mathbf{y}) \, d\mathbf{y}.
    \end{equation}
\end{remark}

In the light of the above remark, the system in \eqref{gen_sde} leads to the following Fokker-Planck equation
\begin{equation}\label{FP}
     \rho_t - \frac{\xi^2}{2}\Delta \rho +\div_x((W*\rho-\nabla V^\text{tar})\cdot\rho)=0,
\end{equation}
where $\rho=\rho(t,\mathbf{x})$ is the spatial probability density, and  
$$
(W*\rho)(t,\mathbf{x})=\int_{\T^d}W(\mathbf{x}-\mathbf{y})\rho(t,\mathbf{y})d\mathbf{y},
$$
with
$W(\mathbf{x}) = Q(\mathbf{x})- \mathbf{x}\cdot a(\mathbf{x}).$

Following the discussion in \cite[page $19$]{carrillo2014derivation} to provide  $Q*\rho$ is being Lipschitz continuous with respect to the spatial variable,  we assume 
\begin{equation}\label{ass:Lip_cond}
   Q(0)=0,\;\; |Q(\mathbf{x})|\le \frac{C}{|\mathbf{x}|^\beta}\;\; \mbox{and}\;\; |\nabla Q(\mathbf{x})|\le \frac{C}{|\mathbf{x}|^{1+\beta}},\;\; \forall \mathbf{x}\in\T^d\setminus \{0\}, 
\end{equation}
for $-1\le\beta< d-1.$

Assuming that $a(\cdot)$ is continuously differentiable and  $Q(\cdot)$ satisfies \eqref{ass:Lip_cond}, we obtain that  $(W*\rho-\nabla V^\text{tar})$ is  continuous in
time and Lipschitz continuous in $\mathbf{x}$.Consequently, by \cite[Lemma $1.6$]{achdou2019mean}, we deduce that $\rho$ is a $1/2$ Holder continuous in time, i.e.
$$
W_2(\rho(t,\cdot), \rho(s,\cdot))\le c_0(1+\|W*\rho-\nabla V^\text{tar}\|_\infty)|t-s|^{\frac{1}{2}}, \;\; \forall s,t\in[0,T].
$$
This property meets the requirements of Assumption \ref{ass:att_HVU_densities}.

Now, we  consider a slight modification of a Fokker-Planck equation \eqref{FP}, which solution $\rho(t,\cdot)$ qualifies as a Caffarelli measure in accordance with Assumption~\ref{ass:att_Caffarelli}. Let numbers $\lambda,\Lambda$ such that $0<\lambda<\Lambda<+\infty$ are fixed. We consider $\delta\cdot W(\mathbf{x})$ in Fokker-Planck equation \eqref{FP} instead of $W(\mathbf{x}),$ where $\delta>0$ is a constant that we are allowed to change. 
First, observe that $\rho(0,\mathbf{x})$, when periodized on the compact torus $\T^d$, is  continuous, and bounded from above and below. 
We take $\rho(0,\mathbf{x})$ such that 
\[
0 < \lambda e^{K_0 T} \le \rho(0,\mathbf{x}) \le \Lambda e^{-K_0T} < \infty \;\;\text{for all} \;\; \mathbf{x} \in \T^d,
\]
where $K_0$ will be chosen later.

The drift vector field $W*\rho-\nabla V^\text{tar}$  has uniformly bounded divergence. Indeed,
$$
\operatorname{div} ((\delta \cdot W)*\rho - \nabla V^\text{tar}) = \delta \cdot\operatorname{div}(W*\rho) - b d,
$$
taking the $L^\infty$ norm and using the triangle inequality we obtain
$$
\|\operatorname{div} ((\delta \cdot W)*\rho - \nabla V^\text{tar})\|_{L^\infty} \le \delta \cdot\|\operatorname{div}(W*\rho)\|_{L^\infty} + |b| d \le \delta \cdot d \cdot L_W + |b| d:=K_0.
$$
Since $K_0 \ge 0$ and $T > 0$, we have $e^{-K_0 T} \le 1 \le e^{K_0 T}$. At time $t=0$, this implies
$$
\lambda \le\lambda e^{K_0 T} \le \rho(0,\mathbf{x}) \le \Lambda e^{-K_0T} \le \Lambda.
$$
Since the divergence is bounded by $K_0$ we can apply the exponential Grönwall bounds established in Lemma \ref{lem:fp_bounds}. This implies that the evolving density cannot vanish or blow up in finite time, yielding two-sided pointwise estimates
$$
\lambda\le \rho(t,\mathbf{x}) \le \Lambda,
\quad \forall (t,\mathbf{x}) \in [0,T] \times \T^d.
$$
Thus, the density $\rho$ satisfies the two-sided pointwise estimates for the entire time interval
\[
0 < \lambda \le \rho(t,\mathbf{x}) \le \Lambda < \infty, \quad \forall (t,\mathbf{x}) \in [0,T] \times \T^d.
\]
This establishes the uniform bounds $\lambda$ and $\Lambda$. The only thing that remains is to provide the following inequality
$$
\lambda e^{K_0 T} < \Lambda e^{-K_0T}.
$$ 
To this aim, we rewrite it as 
$$
K_0 < \frac{1}{2T}\ln\frac{\Lambda}{\lambda}.
$$
On the other hand $K_0=\delta \cdot d \cdot L_W + |b| d <\frac{1}{2T}\ln\frac{\Lambda}{\lambda}.$ This can be provided by taking the constants $\delta$ and $b$ sufficiently small.


To verify the remaining spatial regularity requirement, we view the Fokker-Planck equation as a linear parabolic PDE with smooth coefficients and bounded divergence terms. Since the solution is uniformly bounded between $\lambda$ and $\Lambda$, then a classical parabolic Schauder regularity theory provides that for any chosen exponent $\gamma \in (0,1)$, the spatial Hölder norm of the density remains uniformly bounded independently of time. That is, there exists a constant $K > 0$, depending on $\lambda$, $\Lambda$, $\varepsilon$, and the domain geometry, such that $\sup_{t \in [0,T]} \|\rho(t,\cdot)\|_{C^\gamma(\T^2)} \le K$.

In particular, for our numerical depiction, we take $d=2$, $\xi=0.1$, $\beta = 0.8$, $a(x,y)=0.5$, and 
$$
Q(x,y) := \left[\frac{0.5x}{(x^2+y^2)^{\frac{\beta+1}{2}}}, \frac{0.5y}{(x^2+y^2)^{\frac{\beta+1}{2}}}\right].
$$
It is easy to check that $Q$ satisfies \eqref{ass:Lip_cond}. By taking $\T^2$ to be a macroscopic torus, we can perform our numerical simulations on the local subset $[-5,5]^2 \subset \T^2$ using the standard Euclidean metric. Therefore, for every $(x,y) \in [-5,5]^2$, the function $W$ takes the form
$$
W(x,y) = \frac{1}{2}\left[\frac{x}{(x^2+y^2)^{\frac{9}{10}}}-x, \frac{y}{(x^2+y^2)^{\frac{9}{10}}}-y\right].
$$

\begin{figure*}[htbp!]
\vspace{-0.5cm}

\centering

\subfloat[\centering 3D View $t = 0.00$]{%
\includegraphics[width=0.25\textwidth]{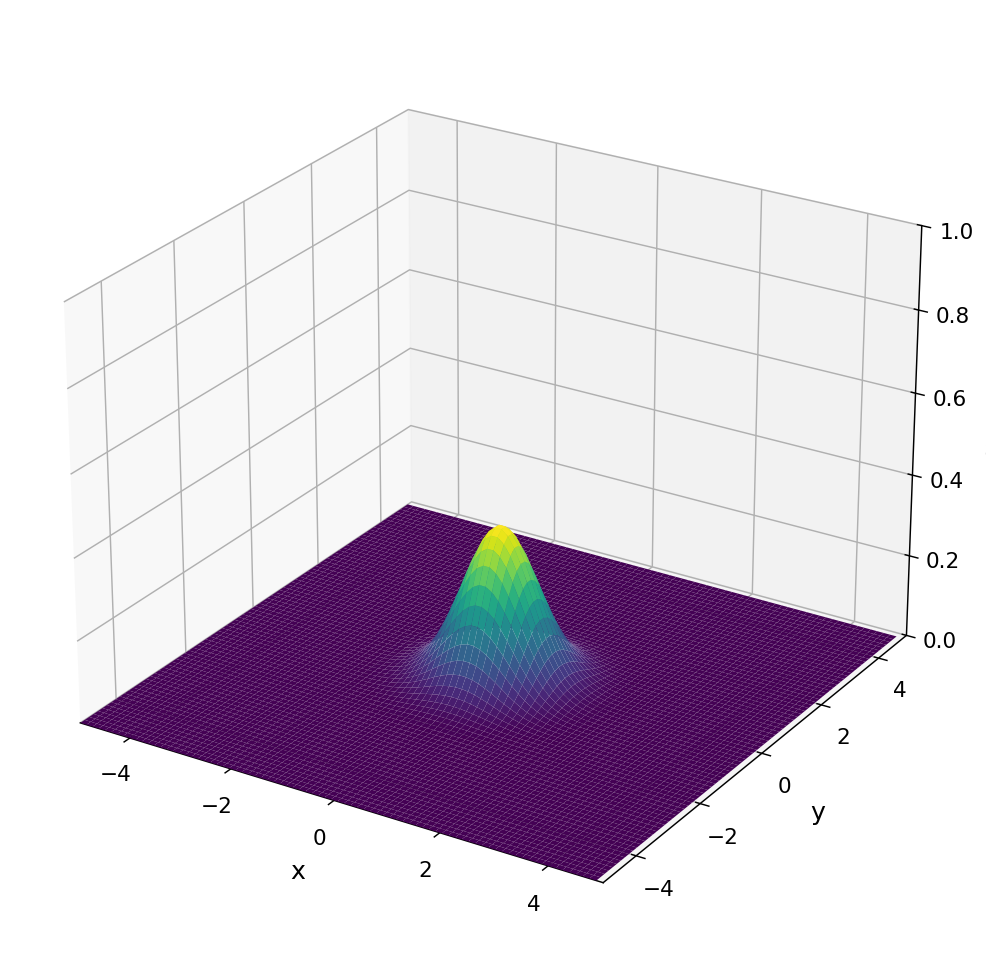}}
\hfill
\subfloat[\centering 3D View $t = 3.33$]{%
\includegraphics[width=0.25\textwidth]{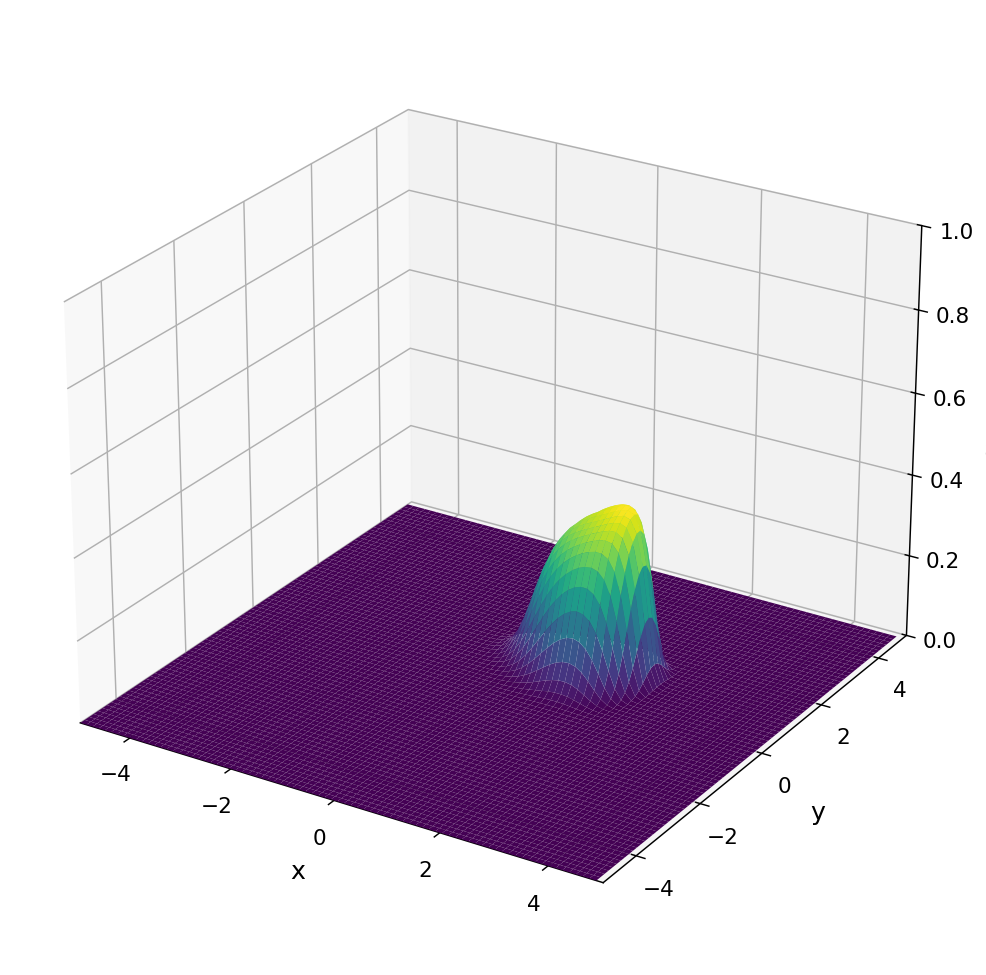}}
\hfill
\subfloat[\centering 3D View $t = 6.67$]{%
\includegraphics[width=0.25\textwidth]{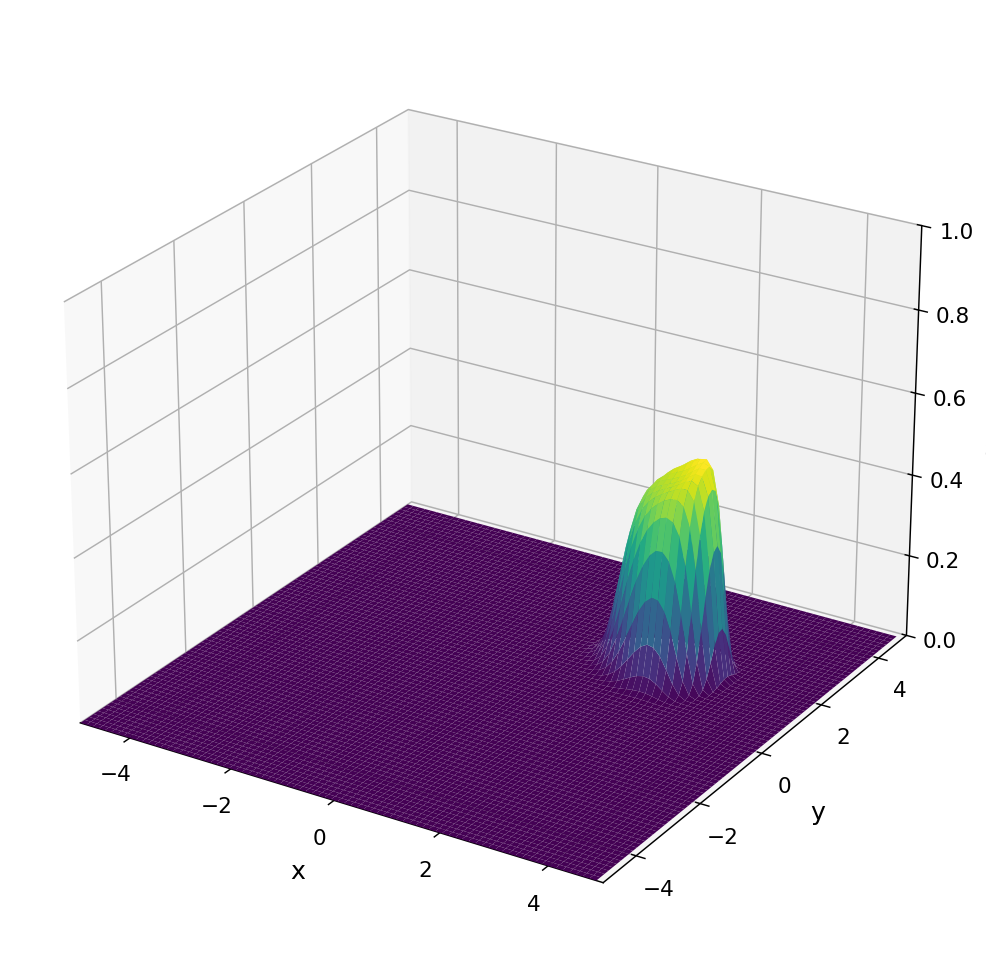}}
\hfill
\subfloat[\centering 3D View $t = 10.00$]{%
\includegraphics[width=0.25\textwidth]{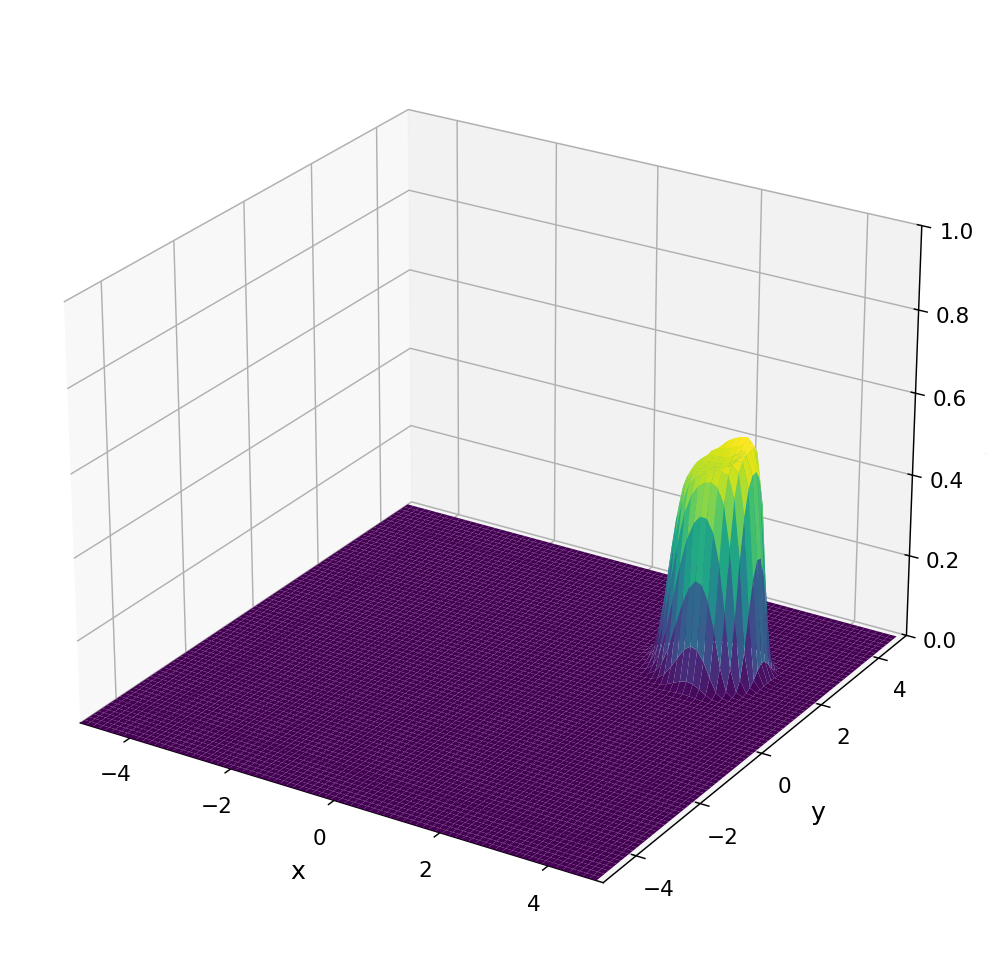}}

\caption{Simulation of formation dynamics with $T=10$, and $(x,y)\in[-5,5]^2$.}
\vspace{-0.4cm}
\label{opt_control}
\end{figure*}
\end{example}

In Figure \ref{opt_control}, we present the numerical simulation of this formation process.
As time increases, the population moves toward the target.

\begin{example}[Optimal Control Dynamics]\label{ex:op-tr}
We consider a  dynamics based on optimal control. The system is given by
\begin{equation}\label{pre_mfg_sys}
    \begin{cases}
        -w_t-\varepsilon\Delta w + \frac{|D w|^2}{2}=0,\quad &(t,x,y)\in(0,T)\times\T^2\\
        \rho_t - \varepsilon\Delta \rho -\nabla\cdot(D w\cdot \rho)=0,\quad &(t,x,y)\in(0,T)\times\T^2\\
        \rho(0,x,y)=\rho_0,\quad & w(T,x,y)= \delta\cdot ((x-x_0)^2+(y-y_0)^2),
\end{cases} 
\end{equation}
where $\rho_0=\mathcal{N}(M_0^{att},1)$ with $M_0^{att}=(0,0)$ and $\delta>0$ is a constant. As time increases, the population successfully moves toward the target.

To show that this optimal control framework is well-posed within our theoretical setting, we  verify that the evolving probability density $\rho(t,\cdot)$ satisfies both the temporal continuity requirements of Assumption \ref{ass:att_HVU_densities} and the spatial regularity   Assumption  \ref{ass:att_Caffarelli}.

Let $(w,\rho)$ be the solution of \eqref{pre_mfg_sys}.
The classical Hamilton-Jacobi theory implies that $w \in C^{1,2}([0,T] \times \T^2)$. In particular, by applying  the Hopf-Cole transformation, the spatial gradient is uniformly bounded over the entire domain, meaning there exists a finite constant $C_0 > 0$, depending only on the diffusion  $\varepsilon$ and the terminal data $\|Dw(T,\cdot)\|_{L^\infty(\T^2)}$, such that $\sup_{t \in [0,T]} \|Dw(t,\cdot)\|_{L^\infty(\T^2)} \le C_0$. Furthermore, interior regularity for parabolic equations ensures that the Laplacian of the value function is similarly bounded, yielding
$$
\|\Delta w\|_{L^\infty([0,T] \times \T^2)} \le \delta K_0,
$$ 
for some positive constant $K_0$ not depending on $\delta$.

 Following the analytical framework of \cite[Lemma $1.6$]{achdou2019mean}, the $L^\infty$ bound on the drift velocity $Dw$ controls the rate of change of the Wasserstein metric, establishing that the solution $\rho(t,\cdot)$ is uniformly $1/2$-Hölder continuous in time
$$
W_2\left(\rho(t,\cdot),\rho(s,\cdot)\right)\le c_0\left(1+\|Dw\|_{L^\infty([0,T] \times \T^2)}\right)|t-s|^{\frac{1}{2}} \le c_0(1+C_0)|t-s|^{\frac{1}{2}},
$$
for all $s,t \in [0,T]$, where $c_0 > 0$ is a constant depending only on the domain $\T^2$ and the time horizon $T$. This inequality rigorously confirms that $\rho$ fulfills the exact time-Hölder condition  in Assumption ~\ref{ass:att_HVU_densities}.

Proceeding the same steps as in the previous example and again applying Lemma \ref{lem:fp_bounds}, we can show that if we take 
\[
0 < \lambda e^{\delta K_0 T} \le \rho(0,x,y) \le \Lambda e^{-\delta K_0T} < \infty \;\;\text{for all} \;\; (x,y) \in \T^2,
\]
then
$$
0 < \lambda  \le \rho(t,x,y) \le \Lambda < \infty, \quad \forall (t,x,y) \in [0,T] \times \T^2.
$$
In this case we need to provide
$$
\delta K_0< \frac{1}{2T}\ln\frac{\Lambda}{\lambda},
$$
which can be done by assuming $\delta < \frac{1}{2TK_0}\ln\frac{\Lambda}{\lambda}.$

By identifying the attacking density as $\mu(t,\cdot) \equiv \rho(t,\cdot)$,
the above two-sided inequality provides that $\rho$ is satisfying Assumption \ref{ass:att_Caffarelli}.

 In Figure \ref{opt_ctrl}  we show the numerical implementation of the system  \eqref{pre_mfg_sys}. Again we assume $\T^2$ to be a macroscopic torus, so  we can perform our numerical simulations on the local subset $[-5,5]^2 \subset \T^2$ using the standard Euclidean metric.

\begin{figure*}[htbp!]

\centering

\subfloat[\centering 3D View $t = 0.00$]{%
\includegraphics[width=0.25\textwidth]{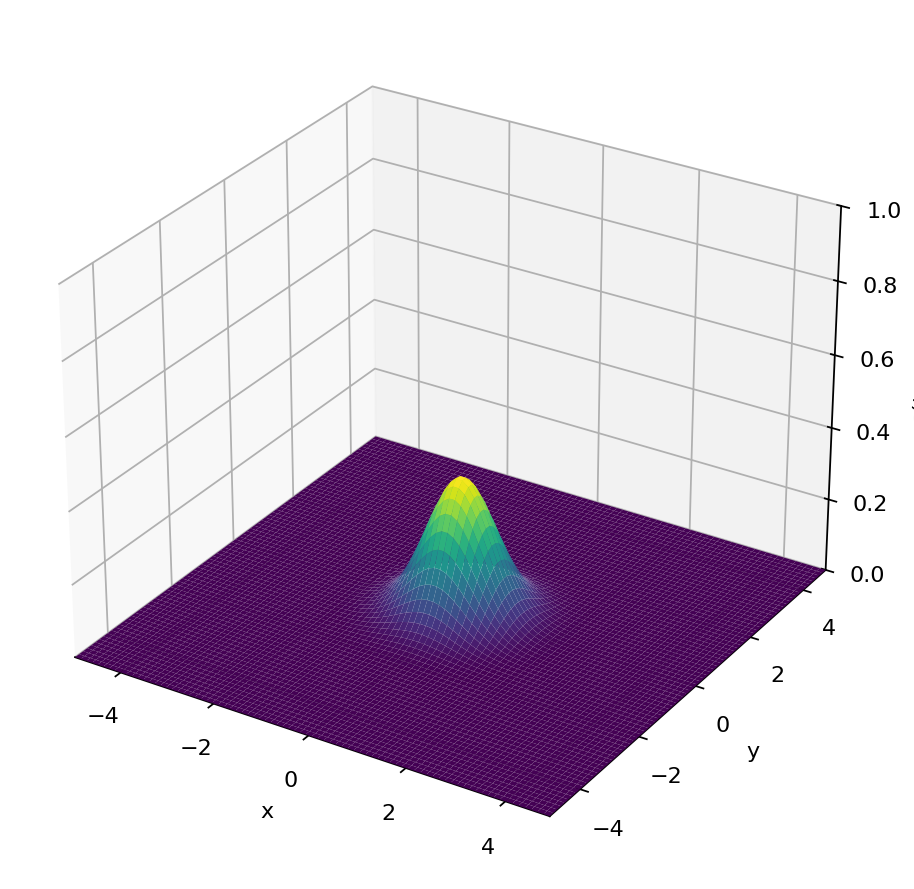}}
\hfill
\subfloat[\centering 3D View $t = 3.33$]{%
\includegraphics[width=0.25\textwidth]{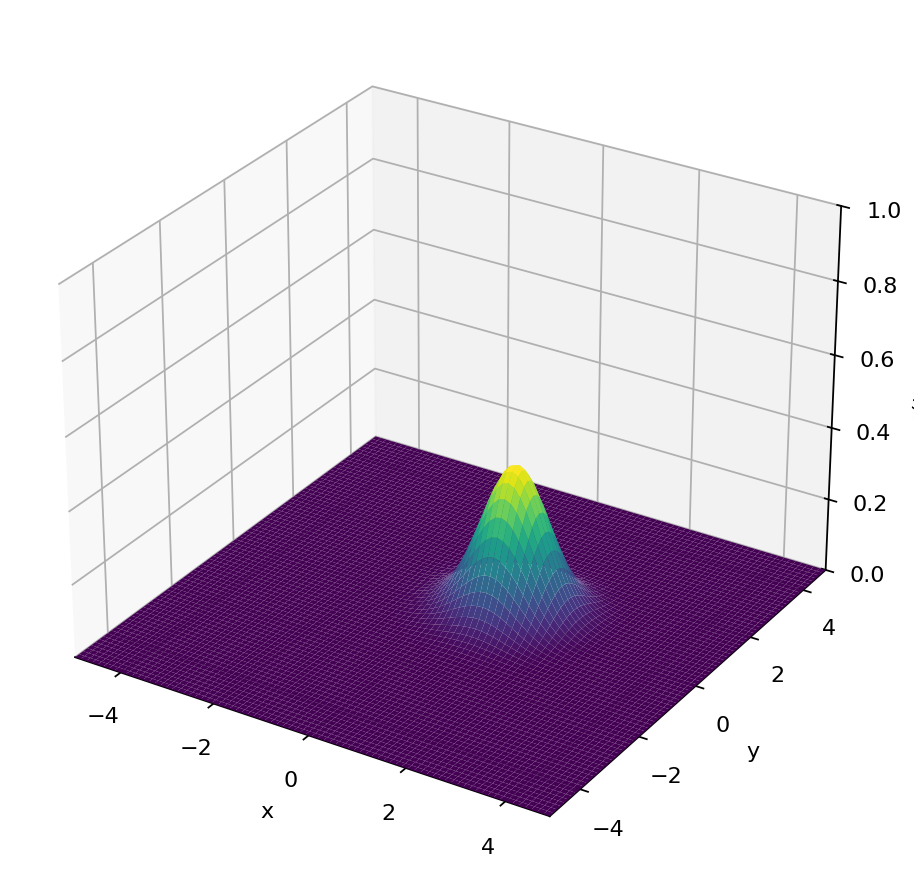}}
\hfill
\subfloat[\centering 3D View $t = 6.67$]{%
\includegraphics[width=0.25\textwidth]{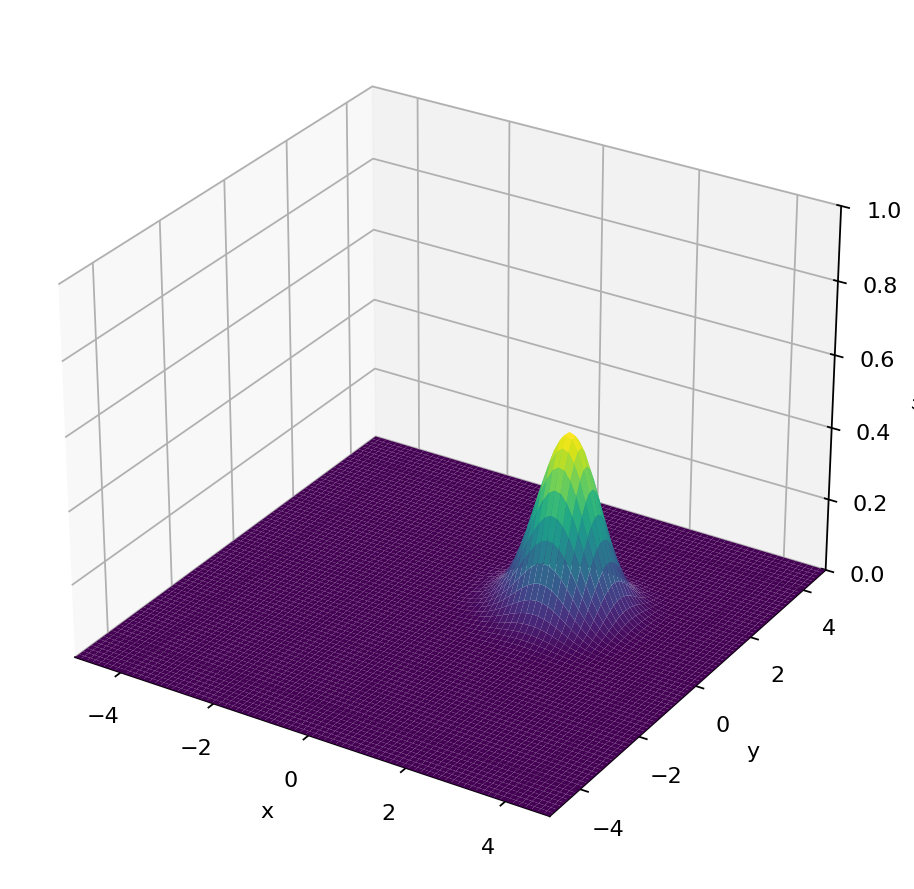}}
\hfill
\subfloat[\centering 3D View $t = 10.00$]{%
\includegraphics[width=0.25\textwidth]{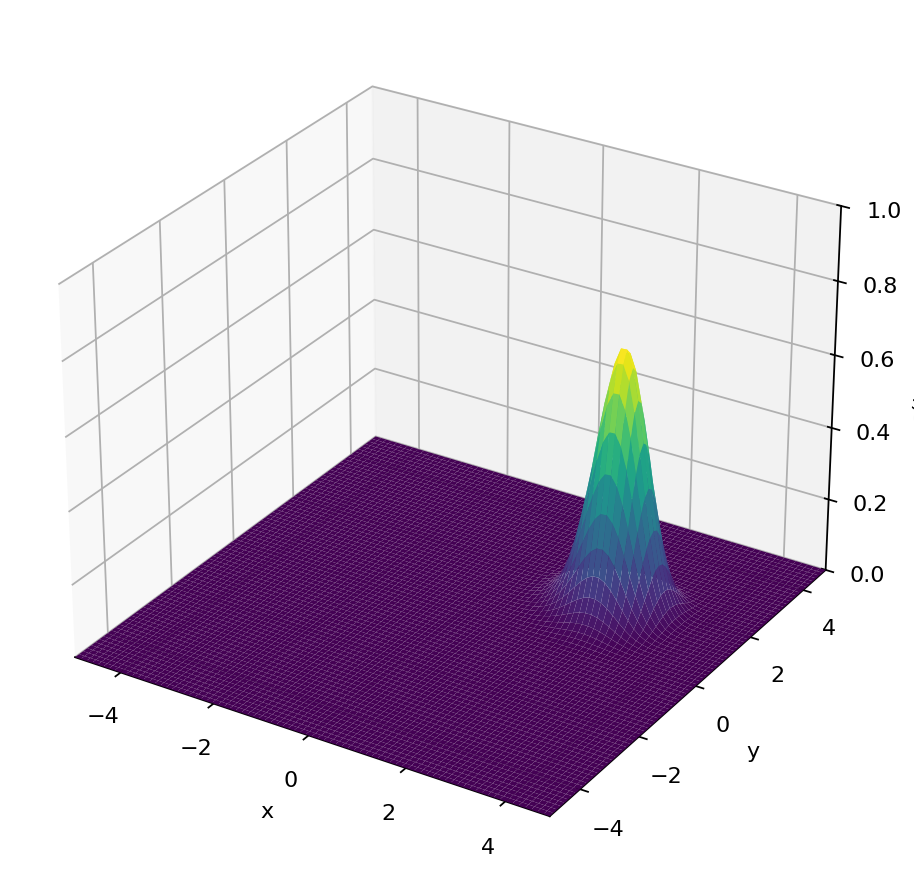}}

    \caption{Simulation  of optimal control dynamics with $\varepsilon =1/200$, $T=10$, $(x_0,y_0)=(4,4)$ and $(x,y)\in[-5,5]^2.$}
\vspace{-0.4cm}
\label{opt_ctrl}
\end{figure*}
\end{example}

In the formation model, both the cohesion and avoidance terms are taken with coefficient \(0.5\). This gives a stronger interaction between the agents and helps the density maintain a wider formation instead of becoming highly concentrated. As time increases, the density moves toward the target while keeping a broader and less peaked shape. In contrast, in the optimal-control model, the density moves toward the target while approximately preserving its initial peaked shape.

\subsection{Regularity of $F$}
In this section, we show regularity results for $F$ defined in \eqref{eq:MFG} with respect to spatial variable $x$ and time variable $t.$ To this aim, we provide general regularity and stability results for Kantorovich potentials of given distributional trajectories. For given $\sigma_1$ and $\sigma_2$ we define $\varphi_{\sigma_1(t)\to\sigma_2(t)}$ to be the Kantorovich potential of optimal transport between the measures $\sigma_1(t,\cdot)$ to $\sigma_2(t,\cdot).$ Note that 
$$
\varphi_{\sigma_1(t) \to \sigma_2(t)}\equiv\varphi_{\sigma_1\to \sigma_2}(t,x):[0,T]\times\T^d\to\R.
$$

\begin{proposition}\label{potential-estimate}
Let $\sigma_i(t,x)$ satisfy Definition \ref{ass: Caffarelli-cond} with  $0<\lambda_i,\Lambda_i<\infty,K>0,\gamma\in(0,1),i=1,2.$ Then 
$$
\sup_{t} \|D_x^2 \varphi_{\sigma_1(t)\to\sigma_2(t)}\|_{L^\infty(\T^d)} \le C_{\mathrm{OT}}(d, \gamma, \lambda_i, \Lambda_i,K),
$$
and  
$$
\|\varphi_{\sigma_1(t) \to \sigma_2(t)} - \varphi_{\sigma_1(s) \to \sigma_2(s)}\|_{L^\infty(\T^d)} \le C_{\T^d,d} \left(W_2(\sigma_1(s), \sigma_1(t))^{\frac{1}{2}} +W_2(\sigma_2(s),\sigma_2(t))^\frac{1}{2}\right)^\frac{2}{d+2}.
$$
\end{proposition}

\begin{proof}
The estimate 
\[
\sup_{t} \|D_x^2 \varphi_{\sigma_1(t)\to\sigma_2(t)}\|_{L^\infty(\T^d)} \le C_{\mathrm{OT}}(d, \gamma, \lambda_i, \Lambda_i,K),
\]
is  a direct consequence of Caffarelli regularity result \cite{caffarelli1992}, because for each $t\in[0,T]$ the distributions  $\sigma_i$ satisfy the Caffarelli measure bounds, with the same $\lambda_i$ and $\Lambda_i.$

Our aim is to prove $1/(2d+4)$ Holder regularity with respect to $t.$ To this purpose, we need to apply the stability results for Kantorovich potentials. In recent years such  results have been obtained by various authors  under more general settings\cite{kitagawa2025,delalande2023,letrouit2024,gallouet2025}. Also, the interested reader is referred to the lecture notes \cite{letrouit2025} and references therein.

For the second inequality of the proposition we assume  Kantorovich potential has $0$ average, i.e. $\int_{\T^d}\varphi=0.$ We can do this since the Kantorovich potential is defined uniquely up to a constant.  This is also a standard assumption in Optimal Transport theory  to provide uniqueness of such potentials. 

For fixed $t$ and $s$ we define a new optimal transport which maps $\sigma_1(s)$ to a target $\sigma_2(t).$ Let set corresponding intermediate Kantorovich potential by $\tilde{\varphi}_{\sigma_1(s) \to \sigma_2(t)}.$ Then we have
\begin{multline}\label{ineq: intermediate}
    \|\varphi_{\sigma_1(t) \to \sigma_2(t)} - \varphi_{\sigma_1(s) \to \sigma_2(s)}\|_{L^2(\T^d)}\le\\\le 
    \|\varphi_{\sigma_1(t) \to \sigma_2(t)} - \tilde{\varphi}_{\sigma_1(s) \to \sigma_2(t)}\|_{L^2(\T^d)} + \| \tilde{\varphi}_{\sigma_1(s) \to \sigma_2(t)} - \varphi_{\sigma_1(s) \to \sigma_2(s)}\|_{L^2(\T^d)}.
\end{multline}
Now we estimate each term in the above sum.
According to Theorem $6$ in \cite{gallouet2025}, which shows the stability estimate for transport maps, and using Poincare inequality (due to mean over torus is $0$), we get 
\begin{multline}\label{ineq: first-term}
    \| \tilde{\varphi}_{\sigma_1(s) \to \sigma_2(t)} - \varphi_{\sigma_1(s) \to \sigma_2(s)}\|_{L^2(\T^d)} \le\\\le C\|\nabla \tilde{\varphi}_{\sigma_1(s) \to \sigma_2(t)} - \nabla\varphi_{\sigma_1(s) \to \sigma_2(s)}\|_{L^2(\T^d)}\le\\ \le C_{\lambda,\Lambda,d,\T^d} W_1(\sigma_2(s),\sigma_2(t))^\frac{1}{2} \le C_{\lambda,\Lambda,d,\T^d} W_2(\sigma_2(s),\sigma_2(t))^\frac{1}{2}.
\end{multline}
For the  remaining term, we consider the dual potentials to $\varphi_{\sigma_1(t) \to \sigma_2(t)}$ and $\tilde{\varphi}_{\sigma_1(s) \to \sigma_2(t)}$ denoted by $\psi_{\sigma_2(t) \to \sigma_1(t)}$ and $\tilde{\psi}_{\sigma_2(t) \to \sigma_1(s)},$ respectively. 
For these dual potentials, the source is fixed at $\sigma_2(t)$ (which is bounded away from zero because it satisfies the Definition \ref{ass: Caffarelli-cond}), and the target changes from $\sigma_1(t)$ to $\sigma_1(s)$. We can apply the exact same target-stability estimate  to these dual potentials as in \eqref{ineq: first-term}. We get 
\begin{equation}\label{ineq:second-term}
    \| \psi_{\sigma_2(t) \to \sigma_1(t)} - \tilde{\psi}_{\sigma_2(t) \to \sigma_1(s)}\|_{L^2(\T^d)} \le C_{\lambda,\Lambda,d,\T^d} W_2(\sigma_1(s),\sigma_1(t))^\frac{1}{2}.
\end{equation}

We now transfer this estimate to the forward potentials.
Their relation to the dual potentials is given by the
$c$-transform. More precisely,
\[
\begin{aligned}
\bigl(\psi_{\sigma_2(t)\to\sigma_1(t)}\bigr)^c(x)
&=
\inf_{y\in\mathbb T^d}
\bigl\{c(x,y)-\psi_{\sigma_2(t)\to\sigma_1(t)}(y)\bigr\},
\\
\bigl(\widetilde\psi_{\sigma_2(t)\to\sigma_1(s)}\bigr)^c(x)
&=
\inf_{y\in\mathbb T^d}
\bigl\{c(x,y)-\widetilde\psi_{\sigma_2(t)\to\sigma_1(s)}(y)\bigr\},
\end{aligned}
\]
where $c(x,y)=d_{\T^d}(x,y)^2.$

The exact $c$-transforms need not have zero average.
By symmetry of the cost and uniqueness of normalized potentials, we have
\begin{equation}\label{eq:p56-recenter}
\begin{aligned}
\varphi_{\sigma_1(t)\to\sigma_2(t)}
&=
\bigl(\psi_{\sigma_2(t)\to\sigma_1(t)}\bigr)^c
-\int_{\mathbb T^d}
\bigl(\psi_{\sigma_2(t)\to\sigma_1(t)}\bigr)^c\,dx,
\\
\widetilde\varphi_{\sigma_1(s)\to\sigma_2(t)}
&=
\bigl(\widetilde\psi_{\sigma_2(t)\to\sigma_1(s)}\bigr)^c
-\int_{\mathbb T^d}
\bigl(\widetilde\psi_{\sigma_2(t)\to\sigma_1(s)}\bigr)^c\,dx.
\end{aligned}
\end{equation}

Let $\mathcal{T}$ and $ \mathcal{\widetilde T}$ be the optimal maps for the transports
$\sigma_1(t)\to\sigma_2(t)$ and $\sigma_1(s)\to\sigma_2(t)$,
respectively. Thus
\[
\mathcal{T}_\#\sigma_1(t)=\sigma_2(t),
\quad
 \mathcal{\widetilde T}_\#\sigma_1(s)=\sigma_2(t).
\]
The definition of the $c$-transform implies
\[
\begin{aligned}
&\bigl(\psi_{\sigma_2(t)\to\sigma_1(t)}\bigr)^c(x)
-\bigl(\widetilde\psi_{\sigma_2(t)\to\sigma_1(s)}\bigr)^c(x)
\ge
-\bigl(\psi_{\sigma_2(t)\to\sigma_1(t)}
-\widetilde\psi_{\sigma_2(t)\to\sigma_1(s)}\bigr)( \mathcal{T}(x)),
\\[0.4em]
&\bigl(\psi_{\sigma_2(t)\to\sigma_1(t)}\bigr)^c(x)
-\bigl(\widetilde\psi_{\sigma_2(t)\to\sigma_1(s)}\bigr)^c(x)
\le
-\bigl(\psi_{\sigma_2(t)\to\sigma_1(t)}
-\widetilde\psi_{\sigma_2(t)\to\sigma_1(s)}\bigr)( \mathcal{\widetilde T}(x)).
\end{aligned}
\]
Squaring the resulting pointwise bound, integrating, and using
the lower bounds on the source densities, we obtain
\[
\begin{aligned}
\left\|
\bigl(\psi_{\sigma_2(t)\to\sigma_1(t)}\bigr)^c
-\bigl(\widetilde\psi_{\sigma_2(t)\to\sigma_1(s)}\bigr)^c
\right\|_{L^2(\mathbb T^d)}^2&
\le
\frac1{\lambda_1}
\int_{\mathbb T^d}
\left|
\bigl(\psi_{\sigma_2(t)\to\sigma_1(t)}
-\widetilde\psi_{\sigma_2(t)\to\sigma_1(s)}\bigr)( \mathcal{T}(x))
\right|^2\,d\sigma_1(t)(x)
\\
&+
\frac1{\lambda_1}
\int_{\mathbb T^d}
\left|
\bigl(\psi_{\sigma_2(t)\to\sigma_1(t)}
-\widetilde\psi_{\sigma_2(t)\to\sigma_1(s)}\bigr)( \mathcal{\widetilde T}(x))
\right|^2\,d\sigma_1(s)(x)
\\
&=
\frac2{\lambda_1}
\int_{\mathbb T^d}
\left|
\psi_{\sigma_2(t)\to\sigma_1(t)}(y)
-\widetilde\psi_{\sigma_2(t)\to\sigma_1(s)}(y)
\right|^2\,d\sigma_2(t)(y)
\\
&\le
\frac{2\Lambda_2}{\lambda_1}
\|\psi_{\sigma_2(t)\to\sigma_1(t)}
-\widetilde\psi_{\sigma_2(t)\to\sigma_1(s)}\|_{L^2(\mathbb T^d)}^2.
\end{aligned}
\]
Thus, \eqref{eq:p56-recenter} implies
\[
\begin{aligned}
\|\varphi_{\sigma_1(t)\to\sigma_2(t)}
-\widetilde\varphi_{\sigma_1(s)\to\sigma_2(t)}
\|_{L^2(\mathbb T^d)}^2&=
\int_{\mathbb T^d}
\left|
\bigl(\psi_{\sigma_2(t)\to\sigma_1(t)}\bigr)^c(x)
-\bigl(\widetilde\psi_{\sigma_2(t)\to\sigma_1(s)}\bigr)^c(x)
\right|^2\,dx
\\
&-
\frac{2}{|\mathbb T^d|}
\left(
\int_{\mathbb T^d}
\left[
\bigl(\psi_{\sigma_2(t)\to\sigma_1(t)}\bigr)^c(x)
-\bigl(\widetilde\psi_{\sigma_2(t)\to\sigma_1(s)}\bigr)^c(x)
\right]\,dx
\right)^2
\\
&+
\frac{1}{|\mathbb T^d|}
\left(
\int_{\mathbb T^d}
\left[
\bigl(\psi_{\sigma_2(t)\to\sigma_1(t)}\bigr)^c(x)
-\bigl(\widetilde\psi_{\sigma_2(t)\to\sigma_1(s)}\bigr)^c(x)
\right]\,dx
\right)^2
\\
&=
\left\|
\bigl(\psi_{\sigma_2(t)\to\sigma_1(t)}\bigr)^c
-\bigl(\widetilde\psi_{\sigma_2(t)\to\sigma_1(s)}\bigr)^c
\right\|_{L^2(\mathbb T^d)}^2
\\
&-
\frac{1}{|\mathbb T^d|}
\left(
\int_{\mathbb T^d}
\left[
\bigl(\psi_{\sigma_2(t)\to\sigma_1(t)}\bigr)^c(x)
-\bigl(\widetilde\psi_{\sigma_2(t)\to\sigma_1(s)}\bigr)^c(x)
\right]\,dx
\right)^2
\\
&\le
\left\|
\bigl(\psi_{\sigma_2(t)\to\sigma_1(t)}\bigr)^c
-\bigl(\widetilde\psi_{\sigma_2(t)\to\sigma_1(s)}\bigr)^c
\right\|_{L^2(\mathbb T^d)}^2.
\\
&\le
\left(\frac{2\Lambda_2}{\lambda_1}\right)
\|\psi_{\sigma_2(t)\to\sigma_1(t)}
-\widetilde\psi_{\sigma_2(t)\to\sigma_1(s)}\|^2_{L^2(\mathbb T^d)}.
\end{aligned}
\]
Hence,
\[
\begin{aligned}
\|\varphi_{\sigma_1(t)\to\sigma_2(t)}
-\widetilde\varphi_{\sigma_1(s)\to\sigma_2(t)}\|_{L^2(\mathbb T^d)}&
\le
\left(\frac{2\Lambda_2}{\lambda_1}\right)^{1/2}
\|\psi_{\sigma_2(t)\to\sigma_1(t)}
-\widetilde\psi_{\sigma_2(t)\to\sigma_1(s)}\|_{L^2(\mathbb T^d)}.
\end{aligned}
\]
Combining this $L^2$ estimate with the \eqref{ineq:second-term} we obtain
$$
 \|\varphi_{\sigma_1(t) \to \sigma_2(t)} - \tilde{\varphi}_{\sigma_1(s) \to \sigma_2(t)}\|_{L^2(\T^d)}\le  C_{\lambda,\Lambda,d,\T^d} W_2(\sigma_1(s),\sigma_1(t))^\frac{1}{2}.
$$
Thus, in view of the above inequality and \eqref{ineq: first-term}, we conclude
$$
   \|\varphi_{\sigma_1(t) \to \sigma_2(t)} - \varphi_{\sigma_1(s) \to \sigma_2(s)}\|_{L^2(\T^d)}\le C_{\lambda,\Lambda,d,\T^d}(W_2(\sigma_1(s),\sigma_1(t))^\frac{1}{2} +W_2(\sigma_2(s),\sigma_2(t))^\frac{1}{2}).
$$
Our next step is  switch to $L^\infty$ norm. To get this done
we apply standard Nirenberg interpolation inequality (see \cite{nirenberg1959})
\begin{multline*}\label{Nirenberg-inqlt}
    \|\varphi_{\sigma_1(t) \to \sigma_2(t)} - \varphi_{\sigma_1(s) \to \sigma_2(s)}\|_{L^\infty(\T^d)} \le \\\le C  \max_{t\in[0,T]}\|D (\varphi_{\sigma_1(t) \to \sigma_2(t)} - \varphi_{\sigma_1(s) \to \sigma_2(s)})\|_{L^\infty(\T^d)}^{\frac{d}{d+2}} \|\varphi_{\sigma_1(t) \to \sigma_2(t)} - \varphi_{\sigma_1(s) \to \sigma_2(s)}\|_{L^2(\T^d)}^{\frac{2}{d+2}}\le \\\le C\cdot(2M)^{\frac{d}{d+2}}\left(W_2(\sigma_1(s), \sigma_1(t))^{\frac{1}{2}}+W_2(\sigma_2(s),\sigma_2(t))^{\frac{1}{2}}\right)^{\frac{2}{d+2}},
\end{multline*}
therefore
$$
\|\varphi_{\sigma_1(t) \to \sigma_2(t)} - \varphi_{\sigma_1(s) \to \sigma_2(s)}\|_{L^\infty(\T^d)} \le C_{\T^d,d} \left(W_2(\sigma_1(s), \sigma_1(t))^{\frac{1}{2}} +W_2(\sigma_2(s),\sigma_2(t))^\frac{1}{2}\right)^\frac{2}{d+2}.
$$
This completes the proof.
\end{proof}

\begin{lemma}\label{lem:1/2d+4 holder}
Let  $\sigma_i(t,x)$ are   Caffarelli measures, which satisfy the Assumptions \ref{ass:attr} and $\nu_{HVU}$ satisfies Assumption \ref{ass:att_HVU_densities}. Then $F(t,x, \nu_{HVU},\sigma_1, \sigma_2)$ is $\frac{1}{2d+4}$ Holder continuous with respect to $t,$ provided 
\begin{equation}\label{ass: 1/2 holder}
    W_2\left(\sigma_i(t),\sigma_i(s)\right)\le C |t-s|^{\frac{1}{2}},
\end{equation}
for all $ t,s \in[0,T]$ and some $C>0.$
\end{lemma}

\begin{proof}
  First, we  show $1/2$ Holder continuity for each multiplied component of the expanded $F,$ except the potential part $\frac{\delta}{\delta \sigma_2}[W_2^2(\sigma_1, \sigma_2)].$ To this aim, we recall the definition of $F$, which equals to
  \[
    F(t,x, \nu_{HVU},\sigma_1, \sigma_2)  = - D(t)\cdot E(t, \sigma_1,\sigma_2)\cdot(d^{att})'(W_2(\sigma_1(t), \sigma_2(t)))\cdot\frac{\delta}{\delta \sigma_2}[W_2^2(\sigma_1, \sigma_2)],
\]
where $D(t) = (1-\alpha)\cdot d^{HVU}\left(W_2(\nu_{HVU}(t), \sigma_1(t))\right)$ and $E(t,\sigma_1,\sigma_2) = e^{-\int_0^t d^{att}(W_2(\sigma_1(\theta), \sigma_2(\theta)))d\theta}$.

Due to $\sigma_1$ and $\nu_{HVU}$ satisfy \eqref{ass: 1/2 holder}, we conclude
\begin{multline*}
|d^{HVU}(W_2(\sigma_1(t),\nu_{HVU}(t)))-d^{HVU}(W_2(\sigma_1(s),\nu_{HVU}(s)))|\le \\ \le C|W_2(\sigma_1(t),\nu_{HVU}(t))-W_2(\sigma_1(s),\nu_{HVU}(s))|\le \\\le C (W_2(\sigma_1(t),\sigma_1(s))+W_2(\nu_{HVU}(t),\nu_{HVU}(s)))\le
 C|t-s|^{1/2}.
\end{multline*}
According to Assumption \ref{ass:attr}, we have $d^{att}$ function is bounded. Hence, for some $L>0$ we have

\begin{multline*}
\left|e^{-\bigintsss_0^t d^{att}(W_2(\sigma_1(\theta), \sigma_2(\theta)))d\theta} - e^{-\bigintsss_0^s d^{att}(W_2(\sigma_1(\theta), \sigma_2(\theta)))d\theta}\right|\le \\\le L\left|\bigintsss_s^t d^{att}(W_2(\sigma_1(\theta), \sigma_2(\theta)))d\theta \right| \le L\bigintsss_s^t |d^{att}(W_2(\sigma_1(\theta), \sigma_2(\theta)))|d\theta \le\\\le L\cdot \underset{\theta\in[0,T]}{\max}\;d^{att}(W_2(\sigma_1(\theta), \sigma_2(\theta))) |t-s|.
\end{multline*}
Thus the integral part in $F$ is Lipschitz with respect to $t,$ therefore also 
 $1/2$ Holder continuous. 
 
 In the light of \eqref{ass: 1/2 holder}, we know that $\sigma_1(t)$ is $1/2$ Hölder continuous with respect to time in $W_2$ metric.
Then we  claim that  the distance $W_2(\sigma_1(t), \sigma_2(t))$ is Hölder $1/2$ in time. Indeed, 
by the triangle inequality of the metric $W_2,$ we have 
$$
W_2(\sigma_1(t), \sigma_2(t)) - W_2(\sigma_1(s), \sigma_2(s)) \le W_2(\sigma_1(t), \sigma_1(s)) + W_2(\sigma_2(t), \sigma_2(s)).
$$
Similarly,
$$
W_2(\sigma_1(s), \sigma_2(s)) - W_2(\sigma_1(t), \sigma_2(t)) \le W_2(\sigma_1(s), \sigma_1(t)) + W_2(\sigma_2(t), \sigma_2(s)).$$
Combining the inequalities  above, we get 
$$
\left| W_2(\sigma_1(t), \sigma_2(t)) - W_2(\sigma_1(s), \sigma_2(s)) \right| \le W_2(\sigma_1(t), \sigma_1(s)) + W_2(\sigma_2(t), \sigma_2(s)).
$$
Thus, 
$$
\begin{aligned}
|W_2(\sigma_1(t), \sigma_2(t)) - W_2(\sigma_1(s), \sigma_2(s))| &\le C_{\sigma_1} |t-s|^{1/2} + C_{\sigma_2} |t-s|^{1/2} \\
&= (C_{\sigma_1} + C_{\sigma_2}) |t-s|^{1/2}.
\end{aligned}
$$
The composition with a smooth function $(d^{att})'(\cdot)$ preserves this Hölder regularity.

Now, we focus on the  potential part. We define by Kantorovich potential  
$$
\varphi_{\sigma_1(t) \to \sigma_2(t)}(x)=\frac{1}{2}\frac{\delta}{\delta \sigma_1}[W_2^2(\sigma_1(t), \sigma_2(t))].
$$
According to Proposition \ref{potential-estimate} and condition \eqref{ass: 1/2 holder}, we have 
\begin{multline}
    \|\varphi_{\sigma_1(t) \to \sigma_2(t)} - \varphi_{\sigma_1(s) \to \sigma_2(s)}\|_{L^\infty(\T^d)} \le\\\le C_{\T^d,d} \left(W_2(\sigma_1(s), \sigma_1(t))^{\frac{1}{2}} +W_2(\sigma_2(s),\sigma_2(t))^\frac{1}{2}\right)^\frac{2}{d+2}
\le C_{d,\T^d} |t-s|^{\frac{1}{2(d+2)}}.
\end{multline}
\end{proof}
\begin{lemma}\label{lm: F-deriv-bounds}
Take $\sigma_i(t,x)$ as in Proposition \ref{potential-estimate}, which satisfy the Assumptions \ref{ass:attr}.
Let $\nu_{HVU}$ satisfies Assumption \ref{ass:att_HVU_densities}. Then there exists a positive constant $C_F$ that depends only on the fixed model data and
the torus geometry, and is independent of the spatial
H\"older bound $K$ in Definition~\ref{ass: Caffarelli-cond}, such that
$$
\|F(\cdot,\cdot,\nu_{HVU},\sigma_1,\sigma_2)
\|_{L^\infty([0,T]\times\T^d)}
+
\|D_xF(\cdot,\cdot,\nu_{HVU},\sigma_1,\sigma_2)
\|_{L^\infty([0,T]\times\T^d)}\le C_F.
$$

\end{lemma}
\begin{proof}
  Set $\mathcal{D}=\operatorname{diam}(\T^d)$. For all $x,x',y\in\T^d$, we have
\[
\begin{aligned}
|c(x,y)-c(x',y)|
=
\left|d_{\T^d}(x,y)^2-d_{\T^d}(x',y)^2\right|&\le
\left(d_{\T^d}(x,y)+d_{\T^d}(x',y)\right)
d_{\T^d}(x,x')\\
&\le 2\mathcal{D}\,d_{\T^d}(x,x').
\end{aligned}
\]
Hence, by the $c$-transform representation, for $i\ne j$, we get
\[
\begin{aligned}
&\left|
\varphi_{\sigma_i(t)\to\sigma_j(t)}(x)
-\varphi_{\sigma_i(t)\to\sigma_j(t)}(x')
\right|\le
\sup_{y\in\T^d}|c(x,y)-c(x',y)|
\le2 \mathcal{D}\,d_{\T^d}(x,x'),
\end{aligned}
\]
and therefore
\[
\sup_{t\in[0,T]}
\|D_x\varphi_{\sigma_i(t)\to\sigma_j(t)}\|_{L^\infty(\T^d)}
\le 2\mathcal{D}.
\]
Using the normalization
$\int_{\T^d}\varphi_{\sigma_i(t)\to\sigma_j(t)}\,dx=0$,
we also obtain
\[
\begin{aligned}
\left|\varphi_{\sigma_i(t)\to\sigma_j(t)}(x)\right|
&=
\left|
\frac1{|\T^d|}
\int_{\T^d}
\left(
\varphi_{\sigma_i(t)\to\sigma_j(t)}(x)
-\varphi_{\sigma_i(t)\to\sigma_j(t)}(y)
\right)\,dy
\right|\\
&\le
\frac{2\mathcal{D}}{|\T^d|}
\int_{\T^d}d_{\T^d}(x,y)\,dy
\le 2\mathcal{D}^2.
\end{aligned}
\]
Therefore,
\[
\sup_{t\in[0,T]}
\left(
\|\varphi_{\sigma_i(t)\to\sigma_j(t)}\|_{L^\infty(\T^d)}
+
\|D_x\varphi_{\sigma_i(t)\to\sigma_j(t)}\|_{L^\infty(\T^d)}
\right)
\le 2\mathcal{D}^2+2\mathcal{D}.
\]

For any $\eta,\zeta\in\mathcal P(\mathbb T^d)$,
let the product measure $\eta\otimes\zeta$ describes choosing $x$ according to $\eta$ and $y$ according to $\zeta$, independently.
It is the probability measure on $\mathbb T^d\times\mathbb T^d$
defined by
\[
(\eta\otimes\zeta)(A\times B)=\eta(A)\,\zeta(B)
\]
for all Borel sets $A,B\subseteq\mathbb T^d$.
This product measure
$\eta\otimes\zeta$ belongs to $\Pi(\eta,\zeta)$, since its
marginals are $\eta$ and $\zeta$. Therefore,
\[
\begin{aligned}
W_2^2(\eta,\zeta)
=
\inf_{\pi\in\Pi(\eta,\zeta)}
\int_{\mathbb T^d\times\mathbb T^d}
d_{\mathbb T^d}(x,y)^2\,d\pi(x,y)
&\le
\int_{\mathbb T^d\times\mathbb T^d}
d_{\mathbb T^d}(x,y)^2\,d\eta(x)\,d\zeta(y)
\\
&\le
\mathcal{D}^2
\int_{\mathbb T^d\times\mathbb T^d}
d\eta(x)\,d\zeta(y)
=\mathcal{D}^2,
\end{aligned}
\]
where we used $d_{\mathbb T^d}(x,y)\le \mathcal{D}$ and
$\eta(\mathbb T^d)=\zeta(\mathbb T^d)=1$.
Thus the scalar factors in $F$ satisfy
\[
|D(t)|=
|1-\alpha|\,
\left|d^{HVU}\!\left(
W_2(\nu_{HVU}(t),\sigma_1(t))
\right)\right|
\le
|1-\alpha|\,\|d^{HVU}\|_{L^\infty([0,\mathcal{D}])},
\]
\[
0<E(t,\sigma_1,\sigma_2)=
\exp\left(
-\int_0^t
d^{att}\!\left(W_2(\sigma_1(\theta),\sigma_2(\theta))\right)
\,d\theta
\right)
\le 1,
\]
\[
\left|
(d^{att})'\!\left(W_2(\sigma_1(t),\sigma_2(t))\right)
\right|\le
\|(d^{att})'\|_{L^\infty([0,\mathcal{D}])}.
\]
Since $D$, $E$, and the attrition-rate derivative are independent
of $x$, the definition of $F$ gives
\begin{equation}\label{eq:F-first-derivatives}
\begin{aligned}
\|F(\cdot,\cdot,\nu_{HVU},\sigma_1,\sigma_2)
\|_{L^\infty([0,T]\times\T^d)}
&+
\|D_xF(\cdot,\cdot,\nu_{HVU},\sigma_1,\sigma_2)
\|_{L^\infty([0,T]\times\T^d)}
\\
&\le
2|1-\alpha|\,
\|d^{HVU}\|_{L^\infty([0,\mathcal{D}])}\,
\|(d^{att})'\|_{L^\infty([0,\mathcal{D}])}
\\
&\times
\sup_{t\in[0,T]}
\left(
\|\varphi_{\sigma_2(t)\to\sigma_1(t)}\|_{L^\infty(\T^d)}
+
\|D_x\varphi_{\sigma_2(t)\to\sigma_1(t)}\|_{L^\infty(\T^d)}
\right)
\\
&\le
4|1-\alpha|\,
\|d^{HVU}\|_{L^\infty([0,\mathcal{D}])}\,
\|(d^{att})'\|_{L^\infty([0,\mathcal{D}])}
(\mathcal{D}^2+\mathcal{D})
=:C_F.
\end{aligned}
\end{equation}
It is not hard to see that the constant $C_F$ satisfies the lemma requirements.
This completes the proof.
\end{proof}

\subsection{The Main Result}
In this section, we prove the local existence of a classical solution to our Attacker-Defender MFG system \eqref{eq:MFG}. To this aim, we first  show several auxiliary results, which are important ingredients of the proof of the main result.

\begin{lemma}\label{lem: div_drift_cont}
Take $\sigma_i(t,x)$ as in Proposition \ref{potential-estimate}, which satisfy the Assumptions \ref{ass:attr}
and assume that
\[
W_2(\sigma_i(t),\sigma_i(s))
\le C|t-s|^{1/2},
\quad i=1,2,
\]
for all $s,t\in[0,T]$.
Let $\nu_{HVU}$ satisfies Assumption \ref{ass:att_HVU_densities}. Consider the following HJB equation
\begin{equation}\label{thm:HJB-eq}
     \begin{cases}
        -w_t-\varepsilon\Delta w + \frac{|Dw|^2}{4\alpha}=F(t,x, \nu_{HVU},\sigma_1, \sigma_2), \quad (t,x)\in[0,T]\times \T^d,\\
         w(T,x) = G(x),\; x\in\T^d.
\end{cases}
\end{equation}
 Suppose that
$G\in C^{2+1/(2d+4)}(\T^d).$

Then there is a positive constant $K^*\equiv K^*(d,\varepsilon,\alpha,T,\lambda_i,\Lambda_i,\|G\|_{C^2}),$ such that
\[
\|D^2w\|_{L^\infty([0,T]\times\T^d)}\le K^*.
\]
\end{lemma}

\begin{proof}
First, we apply the Lemma \ref{lm: F-deriv-bounds} to find a constant $C_F$ such that 
$$
\|F(\cdot,\cdot,\nu_{HVU},\sigma_1,\sigma_2)
\|_{L^\infty([0,T]\times\T^d)}
+
\|D_xF(\cdot,\cdot,\nu_{HVU},\sigma_1,\sigma_2)
\|_{L^\infty([0,T]\times\T^d)}\le C_F.
$$
By applying Hopf-Cole transformation using $u= e^{-\alpha^* w},$ where $\alpha^* = \frac{1}{4\alpha\varepsilon}$ and $w$ solves \eqref{thm:HJB-eq}
we arrive at 
$$
\begin{cases}
     u_t+\varepsilon\Delta u =\alpha^*uF(t,x, \nu_{HVU},\sigma_1,\sigma_2 ),\; (t,x)\in[0,T]\times\T^d\\
     u(T,x) = e^{-\alpha^*G(x)}, \; x\in\T^d .
\end{cases}
$$
Following a classical theory of linear parabolic pdes (see \cite[Theorem $9.2.3$]{krylov1996lectures}) and Lemma \ref{lem:1/2d+4 holder} we conclude that $F$ is $1/(2d+4)$ Holder continuous with respect to $t$ and Lipschitz continuous in $x.$ This also implies  that $u\in C^{2+1/(2d+4)}([0,T]\times\T^d)$. The resulting linear parabolic equation allows us to obtain \textit{a priori} gradient estimates via the Bernstein technique. Specifically, applying the differential operator to the function $Z = |Du|^2 + u^2,$  and utilizing the maximum principle ensures that the gradient remains uniformly bounded over $[0,T] \times \mathbb{T}^d$, provided the right-hand side $F$ is bounded. The standard derivation of such bounds using the maximum principle can be found in \cite[Chapter 3]{protter1984}. Thus, $Du$ is uniformly bounded. 
This in turn implies $w\in C^{2+1/(2d+4)}([0,T]\times\T^d),$ 
 moreover, $w$ is bounded and Lipschitz continuous, with  gradient uniform bound $|Dw| \le C_0.$
 To see it, we examine the relationship between the gradient of the solution $w$ and the transformed function $u$. Recall the Hopf-Cole transformation defined by $u(t,x) = e^{-\alpha^* w(t,x)}$. Differentiating this expression with respect to the spatial variable $x$, we obtain
\begin{equation*}
    Du = -\alpha^* e^{-\alpha^* w} Dw = -\alpha^* u Dw.
\end{equation*}
Rearranging this equation allows us to express $Dw$ in terms of $u$ and $Du$:
\begin{equation*}
    Dw = -\frac{1}{\alpha^*} \frac{Du}{u}.
\end{equation*}
Taking the magnitude of both sides, we arrive at the estimate
\begin{equation} \label{eq:gradient_relation}
    |Dw| = \frac{1}{|\alpha^*|} \frac{|Du|}{|u|}.
\end{equation}
To prove that $|Dw|$ is uniformly bounded, it suffices to know that $|Du|$ is bounded from above and show that $|u|$ is strictly bounded away from zero.

Since the gradient  of $u$ is uniformly bounded, then there exists a constant $C < \infty$ such that
\begin{equation*}
    \sup_{t\in [0,T]}\|Du(t,\cdot)\|_{L^\infty(\T^d)} \le C.
\end{equation*}
Regarding the denominator in \eqref{eq:gradient_relation}, we consider the linear equation satisfied by $u$
\begin{equation*}
    u_t + \varepsilon\Delta u - \alpha^* F u = 0.
\end{equation*}
The terminal condition is given by $u(T,x) = e^{-\alpha^* G(x)}$. Since the exponential function is strictly positive and $G$ is bounded on the torus, we have $u(T,x) \ge \delta_0 > 0$ for some constant $\delta_0$. Furthermore, applying the parabolic minimum principle to the linear equation with bounded coefficients implies that the solution $u$ cannot attain a zero value inside the domain given strictly positive boundary/terminal data. Consequently, there exists a constant $\delta > 0$ such that
$$
    \inf_{(t,x) \in [0,T]\times\mathbb{T}^d} |u(t,x)| \ge \delta\;\; \text{and}\;\; \sup_{(t,x) \in [0,T]\times\mathbb{T}^d} |u(t,x)| \le \frac{1}{\delta}.
$$
Finally, substituting these bounds back into \eqref{eq:gradient_relation}, we deduce that
\begin{equation*}
    |Dw(t,x)| \le \frac{1}{|\alpha^*|} \frac{C}{\delta}.
\end{equation*}
Setting $C_0 = \frac{C}{|\alpha^*|\delta}$, we conclude that $w$ is Lipschitz continuous with $|Dw| \le C_0$.

To establish the uniform bound on $D^2w$, we continue by relating it to the second derivatives of the transformed function $u$. By differentiating the identity $w = -\frac{1}{\alpha^*} \ln u$ twice with respect to the spatial variables $x_i$ and $x_j$, we obtain
$$
\partial_{x_i x_j}^2 w= -\frac{1}{\alpha^*} \left( \frac{\partial_{x_i x_j}^2 u}{u} - \frac{\partial_{x_i} u\cdot \partial_{x_j} u}{u^2} \right).
$$
Taking the absolute value yields the estimate
$$
|D^2w| \le \frac{1}{|\alpha^*|} \left( \frac{|D^2u|}{|u|} + \frac{|Du|^2}{|u|^2} \right).
$$
We have already established that $\inf |u| \ge \delta > 0$ and $\|Du\|_{L^\infty} \le C$. Therefore, to prove that $\|D^2w\|_{L^\infty}$ is uniformly bounded, it suffices to show that the second spatial derivatives of $u$ are uniformly bounded, i.e., $\|D^2u\|_{L^\infty} \le C'$.We evaluate this by looking at the equation governing the second derivatives of $u$. Let $u_{ij}\equiv\partial_{x_i x_j}^2 u$. Differentiating the linear equation $u_t + \varepsilon\Delta u = \alpha^* F u$ twice with respect to the spatial variables gives
$$
(u_{ij})_t + \varepsilon\Delta u_{ij} - \alpha^* F u_{ij} = \alpha^* (F_{ij}u + F_i u_j + F_j u_i).
$$
This is a linear parabolic PDE for $u_{ij}$. 
We estimate $u_{ij}$ through the spatial Lipschitz constants of
$u_i:=\partial_{x_i}u$. The equation for $u_i$ is
\begin{equation}\label{eq:58-first-derivative}
(u_i)_t+\varepsilon\Delta u_i
=\alpha^*(Fu_i+F_i u).
\end{equation}
By \eqref{eq:F-first-derivatives} and the above bounds we get 
\[
\|u_i\|_{L^\infty([0,T]\times\T^d)}\le C,
\quad
\|\alpha^*(Fu_i+F_i u)\|_{L^\infty([0,T]\times\T^d)}
\le M:=\alpha^*C_F(C+\delta^{-1}).
\]
Both constants are independent of $K$.

Choose $r_0>0$ smaller than the injectivity radius of the fixed
flat torus. With the Frobenius norm for Hessian matrices, set
\[
L=
\|D^2(e^{-\alpha^*G})\|_{L^\infty([0,T]\times\T^d)}
+\frac{2C}{r_0}+\frac{Mr_0}{2\varepsilon},
\quad
\omega(r)=Lr-\frac{M}{4\varepsilon}r^2,
\quad 0\le r\le r_0.
\]
For each $i$, we claim that
\begin{equation}\label{eq:58-modulus}
|u_i(t,x)-u_i(t,y)|
\le\omega\bigl(d_{\T^d}(x,y)\bigr)
\;\;\text{if }d_{\T^d}(x,y)\le r_0.
\end{equation}
This is a consequence of the maximum principle applied to
spatial differences, as follows.

Fix $\zeta>0$ and consider
\[
Z_i(t,x,y)
=u_i(t,x)-u_i(t,y)
-\omega\bigl(d_{\T^d}(x,y)\bigr)-\zeta(T-t).
\]
At $t=T$, this function is nonpositive because
$u_i(T)=\partial_{x_i}(e^{-\alpha^*G})$ and
\[
\omega(r)\ge
\|D^2(e^{-\alpha^*G})\|_{L^\infty}r.
\]
It is also nonpositive when $x=y$, and when
$d_{\T^d}(x,y)=r_0$, because $\omega(r_0)\ge 2C$.
If $Z_i$ had a positive maximum, it would therefore occur at
$t<T$ and $0<r:=d_{\T^d}(x,y)<r_0$.
At this point $\partial_tZ_i\le 0$, using a right derivative
if $t=0$.

Use flat coordinates along the unique minimizing segment from
$y$ to $x$, and let $e$ be its unit direction.
At the spatial maximum, the variation $(x,y)\mapsto(x+he,y-he)$
gives
\[
D_{ee}u_i(t,x)-D_{ee}u_i(t,y)\le 4\omega''(r).
\]
Indeed, to justify this inequality, observe that since $Z_i(t, \cdot, \cdot)$ attains a local spatial maximum at $(x,y)$, its second derivative along any smooth curve passing through this point must be nonpositive. We consider the specific variation parameterized by $h \mapsto (x+he, y-he)$. Because we are working in flat coordinates and the distance $r = d_{\mathbb{T}^d}(x,y) < r_0$ is strictly less than the injectivity radius, the points $x+he$ and $y-he$ remain on the unique minimizing geodesic segment. By moving $x$ forward along $e$ and $y$ backward along $e$, the distance between them expands linearly, giving $d_{\mathbb{T}^d}(x+he, y-he) = r + 2h$ for sufficiently small $h$. We can define 
$$
g(h) := Z_i(t, x+he, y-he) = u_i(t, x+he) - u_i(t, y-he) - \omega(r+2h) - \zeta(T-t).
$$
Differentiating $g(h)$ with respect to $h$, the chain rule yields an inner derivative of $e$ for the $x$ term, $-e$ for the $y$ term, and $2$ for the distance term. This gives the first derivative
$$
g'(h) = D_e u_i(t, x+he) + D_e u_i(t, y-he) - 2\omega'(r+2h).
$$
Notice that the negative sign in front of $u_i(t, y-he)$ is cancelled by the inner derivative $-e$. Differentiating a second time, the inner derivative $-e$ emerges again on the second term, yielding
$$
g''(h) = D_{ee} u_i(t, x+he) - D_{ee} u_i(t, y-he) - 4\omega''(r+2h).
$$
Because $g(h)$ attains its maximum at $h=0$, the necessary optimality condition $g''(0) \le 0$ requires that
$$
D_{ee}u_i(t,x) - D_{ee}u_i(t,y) - 4\omega''(r) \le 0.
$$
Similarly, for each vector $e_k$ in an orthonormal basis perpendicular to
$e$, the variation $(x,y)\mapsto(x+he_k,y+he_k)$ gives
\[
D_{e_ke_k}u_i(t,x)-D_{e_ke_k}u_i(t,y)\le 0.
\]
Summing these inequalities yields
\[
\Delta u_i(t,x)-\Delta u_i(t,y)\le 4\omega''(r).
\]
Equation \eqref{eq:58-first-derivative} and the bound on its
right-hand side now imply
\[
\partial_tZ_i
\ge -4\varepsilon\omega''(r)-2M+\zeta
=\zeta>0,
\]
contradicting $\partial_tZ_i\le 0$. Thus $Z_i\le 0$.
Letting $\zeta\to 0$ and interchanging $x$ and $y$
proves \eqref{eq:58-modulus}.

Since $\omega(r)\le Lr$, we obtain
$\|Du_i\|_{L^\infty}\le L$ for every $i$.
Consequently,
\[
\begin{aligned}
\|D^2u\|_{L^\infty([0,T]\times\T^d)}
\le \sqrt d\,L
=\sqrt d\left(
\|D^2(e^{-\alpha^*G})\|_{L^\infty([0,T]\times\T^d)}
+\frac{2C}{r_0}
+\frac{\alpha^*C_Fr_0}{2\varepsilon}
(C+\delta^{-1})
\right)
=:C'.
\end{aligned}
\]
The terminal term is controlled by $\|G\|_{C^2}$, since
\[
D^2(e^{-\alpha^*G})
=e^{-\alpha^*G}
\left(
(\alpha^*)^2DG\otimes DG-\alpha^*D^2G
\right),
\]
where,
$$
DG \otimes DG = \begin{bmatrix}  (\partial_{x_1}G)^2 & (\partial_{x_1}G)(\partial_{x_2}G) & \dots & (\partial_{x_1}G)(\partial_{x_d}G) \\  (\partial_{x_2}G)(\partial_{x_1}G) & (\partial_{x_2}G)^2 & \dots & (\partial_{x_2}G)(\partial_{x_d}G) \\  \vdots & \vdots & \ddots & \vdots \\ (\partial_{x_d}G)(\partial_{x_1}G) & (\partial_{x_d}G)(\partial_{x_2}G) & \dots & (\partial_{x_d}G)^2  \end{bmatrix}.
$$

Finally, substituting this bound, along with our previous bounds $|u| \ge \delta$ and $|Du| \le C$, back into our algebraic estimate for $|D^2w|$, we deduce
$$
\|D^2w\|_{L^\infty([0,T]\times\T^d)} \le \frac{1}{|\alpha^*|} \left( \frac{C'}{\delta} + \frac{C^2}{\delta^2} \right) := K^*(d,\varepsilon,\alpha,T, \lambda_i, \Lambda_i,\|G\|_{C^2}).
$$
This completes the proof.
    
\end{proof}

\begin{lemma}\label{lem:fp_bounds}
Let $T > 0$ and let $b \in L^\infty([0,T]; W^{1,\infty}(\mathbb{T}^d))$ be a vector field with bounded divergence, satisfying
\begin{equation}
    \|\operatorname{div} b\|_{L^\infty([0,T]\times\mathbb{T}^d)} \le B.
\end{equation}
Let $m \in C^{1,2}([0,T]\times\mathbb{T}^d)$ be a classical positive solution to the Fokker-Planck equation
\begin{equation*}
    \partial_t m - \Delta m + \operatorname{div}(b(t,x) m) = 0 \quad \text{in } (0,T] \times \mathbb{T}^d,
\end{equation*}
with initial condition $m(0, \cdot) = m_0$ such that $0 < \min_{\mathbb{T}^d} m_0 \le \max_{\mathbb{T}^d} m_0 < \infty$.

Then, for all $(t,x) \in [0,T] \times \mathbb{T}^d$, the solution satisfies the following two-sided estimates
\begin{equation*}
    \left(\min_{y\in\mathbb{T}^d} m_0(y)\right) e^{-BT} \le m(t,x) \le \left(\max_{y\in\mathbb{T}^d} m_0(y)\right) e^{BT}.
\end{equation*}
\end{lemma}
\begin{proof}
    The proof repeats  the same lines as in \cite[Theorem 1]{arakelyan2026mfg}, hence is omitted.
\end{proof}



\begin{theorem}[Local Existence]\label{existence}
Suppose that Assumptions~\ref{ass:attr} --
\ref{ass:term_cost} hold on a fixed interval $[0,T_0]$,
where $T_0>0$, $\varepsilon>0$, and $\alpha>0$.
Let the initial defender distribution be a probability density
$m_0\in C^\gamma(\T^d)$ for some $\gamma\in(0,1)$, satisfying
\[
0<\lambda_0\le m_0(x)\le\Lambda_0<\infty.
\]
Furthermore, assume that the attacking density trajectory
$\mu(t,x)$ satisfies Assumption~\ref{ass:att_Caffarelli}
with the same spatial exponent $\gamma$.
Then there exists a time horizon $T\in(0,T_0]$ such that
Problem~\ref{prob-MFG} admits at least one  classical solution $(w,m).$ 
\end{theorem}

\begin{proof}
Without loss of generality, we assume that all considered Kantorovich potentials are normalized to have zero Lebesgue average.

Let $\mathcal C$ be the set of probability density trajectories
$\sigma\in X:=C([0,T]\times\T^d)$ such that
\[
\begin{gathered}
\sigma(0,x)=m_0(x),\quad
\lambda\le\sigma(t,x)\le\Lambda,\quad
\int_{\T^d}\sigma(t,x)\,dx=1,\\
\sup_{s\ne t}
\frac{W_2(\sigma(s),\sigma(t))}{|t-s|^{1/2}}\le C_1,
\quad
\sup_{t\in[0,T]}\|\sigma(t,\cdot)\|_{C^\gamma}\le R.
\end{gathered}
\]
The constants $C_1,R$ and $T\le T_0$ will be chosen below.
First fix $\lambda,\Lambda$ such that
\[
0<\lambda<
\min\left\{\lambda_0,\inf_{[0,T_0]\times\T^d}\mu\right\},
\quad
\Lambda>
\max\left\{\Lambda_0,\sup_{[0,T_0]\times\T^d}\mu\right\}.
\]
We  take
\[
R\ge\max\left\{\|m_0\|_{C^\gamma},
\sup_{t\in[0,T_0]}\|\mu(t,\cdot)\|_{C^\gamma}\right\}.
\]
In particular, the constant trajectory $\sigma(t)=m_0$
belongs to $\mathcal C$.

This set is convex and compact in $X$, and hence also compact
in $C^0([0,T];\mathcal P_2(\T^d))$.
Indeed, convexity of the Wasserstein time constraint is well-known. The other
constraints are convex and all the constraints are closed
under uniform convergence. The spatial $C^\gamma$ bound
gives a compact family of densities in $C^0(\T^d)$.
On this compact family, convergence in $W_2$ and uniform
convergence of densities induce the same topology.
Consequently, the uniform Wasserstein time modulus gives
uniform time equicontinuity in $C^0(\T^d)$.
Compactness in $X$ follows from the Arzel\`a--Ascoli theorem.

We construct a fixed-point mapping in the space of probability
density trajectories. Define $\Psi:\mathcal C\to X$ as follows;
we will verify that $\Psi(\mathcal C)\subseteq\mathcal C$.
For a given $\rho\in\mathcal C$, consider
\[
\begin{cases}
-w_t-\varepsilon\Delta w+\dfrac{|Dw|^2}{4\alpha}
=F(t,x,\nu_{HVU},\mu,\rho),
& (t,x)\in(0,T)\times\T^d,\\[0.4em]
w(T,x)=G(x,\rho(T)),
& x\in\T^d.
\end{cases}
\]
Here the whole trajectory $\rho$ enters the time-integral
factor of $F$.
By applying the Hopf--Cole transformation
$u=e^{-\alpha^*w}$, where
$\alpha^*=1/(4\alpha\varepsilon)$, we arrive at
\[
\begin{cases}
u_t+\varepsilon\Delta u
=\alpha^*uF(t,x,\nu_{HVU},\mu,\rho),
& (t,x)\in(0,T)\times\T^d,\\[0.4em]
u(T,x)=e^{-\alpha^*G(x,\rho(T))},
& x\in\T^d.
\end{cases}
\]
Since $\mu$ and $\rho$ satisfy Assumption~\ref{ass:att_Caffarelli}, then  
according to Lemma~\ref{lm: F-deriv-bounds}, we have
\[
\|F(\cdot,\cdot,\nu_{HVU},\mu,\rho)\|_{L^\infty([0,T]\times\T^d)}
+\|D_xF(\cdot,\cdot,\nu_{HVU},\mu,\rho)\|_{L^\infty([0,T]\times\T^d)}
\le C_F,
\]
where $C_F$ depends only on the fixed model data and the torus,
and is independent of the spatial H\"older radius of $\rho$.

Due to Proposition~\ref{potential-estimate} and
Lemma~\ref{lem:1/2d+4 holder}, $F$ is
$1/(2d+4)$-H\"older continuous in time.
Following classical linear parabolic theory
\cite{krylov1996lectures}, 
and using the spatial Lipschitz continuity of $F$, we obtain
a unique solution
\[
u\in C^{1+1/(4d+8),\,2+1/(2d+4)}([0,T]\times\T^d).
\]
The comparison principle gives the uniform positive bounds
\[
e^{-\alpha^*(M_G+C_FT_0)}
\le u(t,x)\le e^{\alpha^*(M_G+C_FT_0)}.
\]
This in turn implies
\[
w\in C^{1+1/(4d+8),\,2+1/(2d+4)}([0,T]\times\T^d).
\]
Moreover, by Lemma~\ref{lem: div_drift_cont}, we deduce that $w$ is bounded and Lipschitz continuous, with
\[
\|Dw\|_{L^\infty([0,T]\times\T^d)}\le C_0,
\quad
\|D^2w\|_{L^\infty([0,T]\times\T^d)}\le K^*.
\]
Here, the Hessian norm is the Frobenius norm, and the positive
constants $C_0,K^*$ are uniform for $T\le T_0$.
They depend only on the fixed data, the fixed density bounds,
and $M_G$, and are independent of $\rho$, $C_1$, and $R$.

Next, define $m=\Psi(\rho)$ as the solution to the
Fokker--Planck equation
\[
\begin{cases}
m_t-\varepsilon\Delta m
-\dfrac{1}{2\alpha}\operatorname{div}_x(Dw\,m)=0,
& (t,x)\in(0,T)\times\T^d,\\[0.4em]
m(0,x)=m_0(x),
& x\in\T^d.
\end{cases}
\]
Since this equation is linear and $Dw,\Delta w$ are
 H\"older continuous, it has a unique classical solution,
$m$. 
Because $m$ solves the Fokker--Planck equation then by \cite[Lemma $1.6$] {achdou2019mean}, we have
\begin{equation}\label{main-proof-W2}
    W_2(m(t),m(s))
\le c_0\left(1+\frac{1}{2\alpha}\|Dw\|_{L^\infty([0,T]\times\T^d)}\right)
|t-s|^{1/2},
\end{equation}
where $c_0$ depends only on $d,\varepsilon,T_0$.

Next, we obtain uniform H\"older estimate for $m$.
Let
\[
b=-\frac{Dw}{2\alpha},
\quad
B=\frac{\sqrt d}{2\alpha}K^*,
\quad
M_0=\Lambda_0e^{BT_0}.
\]
By Lemmas~\ref{lem: div_drift_cont} and \ref{lem:fp_bounds}, we have
\[
\|b\|_{L^\infty([0,T]\times\T^d)}\le\frac{C_0}{2\alpha},
\quad
\|\operatorname{div}b\|_{L^\infty([0,T]\times\T^d)}\le B,
\quad
\|m\|_{L^\infty([0,T]\times\T^d)}\le M_0.
\]
Writing the Fokker--Planck equation as
\[
m_t-\varepsilon\Delta m=-\operatorname{div}(bm),
\quad m(0)=m_0,
\]
the Hölder estimate for the heat equation (see \cite{ladyzhenskaya1968linear}) gives

\begin{equation}
    \label{main-eq-Ho}
\begin{split}
    \sup_{t\in[0,T]}\|m(t,\cdot)\|_{C^\gamma}
&\le C_\gamma\left(
\|m_0\|_{C^\gamma}+\|bm\|_{L^\infty([0,T]\times\T^d)}
\right)
\\&\le C_\gamma\left(
\|m_0\|_{C^\gamma}
+\frac{C_0}{2\alpha}\Lambda_0e^{BT_0}
\right)=:K_0.
\end{split}
\end{equation}
Here, $C_\gamma$ depends only on
$d,\varepsilon,\gamma,T_0$,
uniformly for $0<T\le T_0$.
Since $C_0$ and $K^*$ are independent of $R$,
the same is true for $K_0$.

Due to Lemmas~\ref{lem: div_drift_cont} and
\ref{lem:fp_bounds}, we also have the two-sided estimates
\[
\begin{aligned}
\lambda_0e^{-BT}
\le e^{-BT}\min_{y\in\T^d}m_0(y)
\le m(t,x)
\le e^{BT}\max_{y\in\T^d}m_0(y)
\le\Lambda_0e^{BT}.
\end{aligned}
\]
The constants $\lambda,\Lambda$ are given by Definition \ref{ass: Caffarelli-cond} and fixed.
Now, choose $T$ so that
\[
\lambda\le\lambda_0e^{-BT},
\quad
\Lambda_0e^{BT}\le\Lambda.
\]
These inequalities ensure that
\begin{equation}\label{main-eq-Caf}
    \lambda\le m(t,x)\le\Lambda
\quad\text{for all }(t,x)\in[0,T]\times\T^d.
\end{equation}
 Finally, taking
\begin{equation*}
    \begin{split}
        &C_1=c_0\left(1+\frac{C_0}{2\alpha}\right),
\quad
R=\max\left\{
K_0,\|m_0\|_{C^\gamma},
\sup_{t\in[0,T_0]}\|\mu(t,\cdot)\|_{C^\gamma}
\right\},\\
& 0<T\le T^*:=
\min\left\{
T_0,\,
\frac1B\log\frac{\lambda_0}{\lambda},\,
\frac1B\log\frac{\Lambda}{\Lambda_0}
\right\},
    \end{split}
\end{equation*}
by \eqref{main-proof-W2}, \eqref{main-eq-Ho}, \eqref{main-eq-Caf} we deduce that  $m\in\mathcal C$.
These choices are noncircular because $C_0,K^*,K_0$
are independent of  $R$, and $C_0,K^*$
are independent of $C_1$.
Therefore the mapping
$\Psi:\rho\mapsto m=\Psi(\rho)$ is a well-defined
self-map of $\mathcal C$.
Let us verify that $\Psi$ is continuous.
Let $\rho_n\to\rho$ in $\mathcal C$ with respect to the
uniform topology of $X$. For brevity, write
\[
F_n=F(\cdot,\cdot,\nu_{HVU},\mu,\rho_n),
\qquad
F=F(\cdot,\cdot,\nu_{HVU},\mu,\rho).
\]
All constants below are independent of $n$.

Since $\rho_n\to\rho$ uniformly, and recalling $\mathcal{D}=\diam(\T^d),$ we have
\[
\begin{aligned}
\sup_{t\in[0,T]}W_2(\rho_n(t),\rho(t))^2\le\mathcal{D}\sup_{t\in[0,T]}W_1(\rho_n(t),\rho(t))
&\le
\frac{\mathcal{D}^2}{2}
\sup_{t\in[0,T]}
\|\rho_n(t)-\rho(t)\|_{L^1(\T^d)}
\\
&\le
\frac{\mathcal{D}^2|\T^d|}{2}
\|\rho_n-\rho\|_X
\longrightarrow0.
\end{aligned}
\]
Thus, setting
\[
\delta_n:=
\sup_{t\in[0,T]}W_2(\rho_n(t),\rho(t)),
\]
we have $\delta_n\to0$.

The stability estimate proved in
Proposition~\ref{potential-estimate}, applied to the pairs
$(\rho_n(t),\mu(t))$ and $(\rho(t),\mu(t))$, yields
\[
\sup_{t\in[0,T]}
\|\varphi_{\rho_n(t)\to\mu(t)}
-\varphi_{\rho(t)\to\mu(t)}\|_{L^\infty(\T^d)}
\le C\delta_n^{1/(d+2)}
\longrightarrow0.
\]
Here all potentials are normalized to have zero Lebesgue
average.

For the time-integral factor, using $d^{att}\ge0$ and
$|e^{-a}-e^{-b}|\le|a-b|$ for $a,b\ge0$, we obtain
\[
\begin{aligned}
&|E(t,\mu,\rho_n)-E(t,\mu,\rho)|
\\
&\quad\le
\int_0^t
\left|
d^{att}\!\left(W_2(\mu(s),\rho_n(s))\right)
-d^{att}\!\left(W_2(\mu(s),\rho(s))\right)
\right|\,ds
\\
&\quad\le
C\int_0^t W_2(\rho_n(s),\rho(s))\,ds
\le CT_0\delta_n.
\end{aligned}
\]
Similarly, since $(d^{att})'$ is Lipschitz on the compact
range of the Wasserstein distances,
\[
\begin{aligned}
\sup_{t\in[0,T]}
\bigl|
(d^{att})'\!\left(W_2(\mu(t),\rho_n(t))\right)
-(d^{att})'\!\left(W_2(\mu(t),\rho(t))\right)
\bigr|
\le C\delta_n.
\end{aligned}
\]
The remaining factor $D(t)$ is independent of $\rho$.
Since all factors are uniformly bounded, it follows that
\[
\|F_n-F\|_{L^\infty([0,T]\times\T^d)}
\le C\left(\delta_n+\delta_n^{1/(d+2)}\right)
\longrightarrow0.
\]

Next, let $u_n=e^{-\alpha^*w_n}$ and
$u=e^{-\alpha^*w}$, where
$\alpha^*=1/(4\alpha\varepsilon)$.
Subtracting their transformed equations gives
\[
\begin{cases}
(u_n-u)_t+\varepsilon\Delta(u_n-u)
-\alpha^*F_n(u_n-u)
=\alpha^*(F_n-F)u,
\\[0.4em]
(u_n-u)(T,x)
=e^{-\alpha^*G(x,\rho_n(T))}
-e^{-\alpha^*G(x,\rho(T))}.
\end{cases}
\]
The continuity assumption on the terminal coupling gives
\[
G(\cdot,\rho_n(T))
\longrightarrow G(\cdot,\rho(T))
\qquad\text{in }C^2(\T^d),
\]
and hence
\[
\|(u_n-u)(T,\cdot)\|_{L^\infty(\T^d)}
\longrightarrow0.
\]
Applying the maximum principle after reversing time,
and using $\|F_n\|_{L^\infty}\le C_F$, we obtain
\[
\begin{aligned}
\|u_n-u\|_{L^\infty}
&\le e^{\alpha^*C_FT_0}
\bigl(
\|(u_n-u)(T,\cdot)\|_{L^\infty(\T^d)}
\\
&+\alpha^*T_0\|u\|_{L^\infty}
\|F_n-F\|_{L^\infty}
\bigr)
\longrightarrow0,
\end{aligned}
\]
where the unspecified $L^\infty$ norms are taken over
$[0,T]\times\T^d$.

The previously established positive lower bound gives
$u_n,u\ge c>0$, uniformly in $n$. Therefore,
\[
\begin{aligned}
\|w_n-w\|_{L^\infty}
&=
\frac1{\alpha^*}
\|\log u_n-\log u\|_{L^\infty}
\\
&\le
\frac1{\alpha^*c}\|u_n-u\|_{L^\infty}
\longrightarrow0.
\end{aligned}
\]
Moreover, Lemma~\ref{lem: div_drift_cont} gives
\[
\|D^2w_n\|_{L^\infty}
+\|D^2w\|_{L^\infty}\le2K^*.
\]
Taylor's formula for $(w_n-w)$ along spatial directions consequently yields,
for every sufficiently small $h>0$,
\[
\|Dw_n-Dw\|_{L^\infty}
\le
\frac{2}{h}\|w_n-w\|_{L^\infty}+K^*h.
\]
First letting $n\to\infty$, and then $h\downarrow0$, proves
\[
\|Dw_n-Dw\|_{L^\infty([0,T]\times\T^d)}
\longrightarrow0.
\]

We now consider the associated Fokker--Planck solutions
$m_n=\Psi(\rho_n)$ and $m=\Psi(\rho)$.
Set
\[
b_n=-\frac{Dw_n}{2\alpha},
\quad
b=-\frac{Dw}{2\alpha},
\quad
z_n=m_n-m.
\]
Subtracting their equations gives
\[
\begin{cases}
(z_n)_t-\varepsilon\Delta z_n
+\operatorname{div}(bz_n)
=-\operatorname{div}\bigl((b_n-b)m_n\bigr),
\\[0.4em]
z_n(0,\cdot)=0.
\end{cases}
\]
Multiplying by $z_n$, integrating over $\T^d$, and integrating
by parts, we obtain
\[
\frac12\frac{d}{dt}\|z_n\|_{L^2(\T^d)}^2
+\varepsilon\|Dz_n\|_{L^2(\T^d)}^2
=
-\frac12\int_{\T^d}(\operatorname{div}b)z_n^2\,dx
+
\int_{\T^d}(b_n-b)m_n\cdot Dz_n\,dx.
\]
Using $\|\operatorname{div}b\|_{L^\infty}\le B$ and
Young's inequality yields
\[
\begin{aligned}
\frac12\frac{d}{dt}\|z_n\|_{L^2(\T^d)}^2
+\varepsilon\|Dz_n\|_{L^2(\T^d)}^2
&\le
\frac B2\|z_n\|_{L^2(\T^d)}^2
+\frac{\varepsilon}{2}\|Dz_n\|_{L^2(\T^d)}^2
\\
&+
\frac{1}{2\varepsilon}
\|b_n-b\|_{L^\infty(\T^d)}^2
\|m_n\|_{L^2(\T^d)}^2.
\end{aligned}
\]
Since $m_n$ is a probability density bounded by $\Lambda$,
\[
\|m_n(t)\|_{L^2(\T^d)}^2
=
\int_{\T^d}m_n(t,x)^2\,dx
\le
\Lambda\int_{\T^d}m_n(t,x)\,dx
=\Lambda.
\]
Consequently,
\[
\frac{d}{dt}\|z_n(t,\cdot)\|_{L^2(\T^d)}^2
\le
B\|z_n(t,\cdot)\|_{L^2(\T^d)}^2
+\frac{\Lambda}{\varepsilon}
\|b_n(t)-b(t)\|_{L^\infty(\T^d)}^2.
\]
Since $z_n(0,\cdot)=0$, Gr\"onwall's inequality gives
\[
\begin{aligned}
\|z_n(t,\cdot)\|_{L^2(\T^d)}^2
&\le
\frac{\Lambda}{\varepsilon}
\int_0^t e^{B(t-s)}
\|b_n(s)-b(s)\|_{L^\infty(\T^d)}^2\,ds
\\
&\le
\frac{\Lambda T_0}{\varepsilon}e^{BT_0}
\|b_n-b\|_{L^\infty([0,T]\times\T^d)}^2.
\end{aligned}
\]
Thus,
\[
\sup_{t\in[0,T]}
\|m_n(t)-m(t)\|_{L^2(\T^d)}^2
\le
\frac{\Lambda T_0}{\varepsilon}e^{BT_0}
\left\|\frac{Dw_n-Dw}{2\alpha}\right\|_
{L^\infty([0,T]\times\T^d)}^2
\to 0.
\]

Finally, the uniform spatial H\"older bounds imply
\[
[z_n(t,\cdot)]_{C^\gamma(\T^d)}
\le
[m_n(t,\cdot)]_{C^\gamma(\T^d)}
+[m(t,\cdot)]_{C^\gamma(\T^d)}
\le2K_0,
\]
where
\[
[f(t,\cdot)]_{C^\gamma(\T^d)}
:=
\sup_{\substack{x,y\in\T^d\\x\neq y}}
\frac{|f(t,x)-f(t,y)|}
     {d_{\T^d}(x,y)^\gamma}.
\]
For $r>0$ smaller than the injectivity radius of $\T^d$,
and every $y\in B_r(x)$, the spatial H\"older bound gives
\[
\begin{aligned}
|z_n(t,x)|
&\le |z_n(t,y)|+|z_n(t,x)-z_n(t,y)|
\\
&\le |z_n(t,y)|+2K_0d_{\T^d}(x,y)^\gamma
\\
&\le |z_n(t,y)|+2K_0r^\gamma.
\end{aligned}
\]
Integrating this inequality with respect to $y$ over $B_r(x)$
and applying the Cauchy--Schwarz inequality, we obtain
\[
\begin{aligned}
|z_n(t,x)|
&\le
\frac1{|B_r|}
\int_{B_r(x)}|z_n(t,y)|\,dy
+2K_0r^\gamma
\\
&\le
|B_r|^{-1/2}
\|z_n(t,\cdot)\|_{L^2(B_r(x))}
+2K_0r^\gamma
\\
&\le
Cr^{-d/2}\|z_n(t,\cdot)\|_{L^2(\T^d)}
+2K_0r^\gamma,
\end{aligned}
\]
where we used $|B_r|=\omega_d r^d$, with $\omega_d$
denoting the volume of the Euclidean unit ball.

Taking the supremum over $(t,x)$ yields
\[
\|z_n\|_{L^\infty([0,T]\times\T^d)}
\le
Cr^{-d/2}
\sup_{t\in[0,T]}\|z_n(t,\cdot)\|_{L^2(\T^d)}
+2K_0r^\gamma.
\]
First letting $n\to\infty$, and then $r\downarrow0$, we obtain
\[
\|\Psi(\rho_n)-\Psi(\rho)\|_X
=
\|m_n-m\|_{L^\infty([0,T]\times\T^d)}
\longrightarrow0.
\]
Therefore $\Psi$ is continuous in $X$.

By the Schauder fixed-point theorem, $\Psi$ has a fixed
point in $\mathcal C$. At this fixed point $\rho=m$,
so the terminal condition is $w(T,x)=G(x,m(T))$.
Therefore this fixed point, together with its associated
HJB solution $w$, is a solution of
Problem~\ref{prob-MFG}.
\end{proof}

\begin{remark}
    It is noteworthy that the same existence result can be proven for general MFG problem \eqref{eq:gen-kantorovich-mfg} by assuming that the right-hand side function $F$ is Holder continuous in time and Lipschitz continuous  with respect to $x$ and the last argument.
\end{remark}

\section{Simulations}\label{sec: simulations}
In this section, we present numerical experiments for the proposed mean-field attacker-defender model. The goal is to illustrate how the model captures the interaction between attacker and defender populations, and how the survival probabilities depend on the attrition parameters, the initial population spread, and the control cost parameter. We first discuss a one-dimensional test case with an explicit solution, which is used to validate the PINN approximation. Then, we present two-dimensional simulations showing different attacker-defender scenarios.

We recall that, in the model, attrition is described at the
population level through the Wasserstein distance between the involved
distributions. More precisely, if
$\mu(t)$ denotes the attacker distribution, $m(t)$ denotes the defender
distribution, and $\nu_H$ denotes the HVU distribution, we use
\[
d^{att}(\mu(t),m(t))
=
\lambda_A
\exp\left(
-\frac{W_2^2(\mu(t),m(t))}{\sigma_A}
\right),
\]
and
\[
d^{H}(\nu_H,\mu(t))
=
\lambda_H
\exp\left(
-\frac{W_2^2(\nu_H(t),\mu(t))}{\sigma_H}
\right).
\]
Here $\lambda_A$ and $\lambda_H$ determine the strength of the attrition
effects, while $\sigma_A$ and $\sigma_H$ control the sensitivity of the
rates to the Wasserstein distance.

The survival probabilities of the attackers and the HVU are respectively given by
$$ Q(t) = \exp\!\left(-\int_0^t d^{\mathrm{att}}(\mu(\tau), m(\tau))\,d\tau\right)\;\;\text{and}\;\; 
    P(t) = \exp\!\left(-\int_0^t d^{H}(\nu_{\mathrm{H}}, \mu(\tau))\,
Q(\tau)\, d\tau\right).$$

For the defenders dynamics we need to compute Wasserstein distance effectively. To this aim, we use the so-called Sinkhorn algorithm that approximates the entropically regularized optimal transport problem. For more details see \cite{cuturi2013sinkhorn} and references therein. This algorithm approximates the $\varepsilon$ close value of the distance, rather than the exact one. This is acceptable from computational perspective.
    
\subsection{Explicit Solution in 1D}
For the sake of exposition, in  this section we discuss $1D$ mean field game attacker-defender system, where we can write its explicit solution.

We consider a coupled system
\begin{equation}\label{exp.1}  
\begin{cases}
 -u_t - \varepsilon \Delta u + \dfrac{u_x^2}{2} + V(t,x) =
 \displaystyle
\frac{\delta}{\delta m}
\Bigg(
W_2^2\big(m(t,x),w(t,x)\big)
\Bigg)W_2\big(w(t,x),\nu(t,x)\big), \\[8pt]
 m_t - \varepsilon \Delta m - (u_x m)_x = 0.
\end{cases}
\end{equation}
We define \( w \) and \( m \) as Gaussian probability densities as follows
\begin{equation}\label{exp.2}
    \begin{cases}
w(t,x) = \dfrac{1}{\sqrt{2\pi\sigma^2_w}} 
\exp\!\left( -\dfrac{(x - \mu_w(t))^2}{2\sigma_w^2} \right),\\  
m(t,x) = \dfrac{1}{\sqrt{2\pi\sigma^2_m}} 
\exp\!\left( -\dfrac{(x - \mu_m(t))^2}{2\sigma_m^2} \right),\\
\nu(t,x) = \dfrac{1}{\sqrt{2\pi\sigma^2_{\nu}}} 
\exp\!\left( -\dfrac{(x - \mu_\nu)^2}{2\sigma_{\nu}^2} \right)
\end{cases}
\end{equation}
where
\[
\mu_m(t) = x_0 \left(1 - \dfrac{t}{T}\right) + x_T \dfrac{t}{T},
\quad 
\mu_m'(t) = \frac{x_T-x_0}{T}.
\]
This choice corresponds to a linear motion of the mean from \( x_0 \) to \( x_T \)
with constant variance.

Dividing the Fokker-Planck equation by \( m \) and substituting \( y := x - \mu_m(t) \), we obtain
\begin{equation}\label{exp.3}
\frac{y(x_T-x_0)}{T\sigma_m^2}
- \varepsilon\frac{y^2}{\sigma_m^4}
+ \frac{\varepsilon}{\sigma_m^2}
- u_{xx}
+ u_x\frac{y}{\sigma_m^2} = 0.
\end{equation}
We seek a solution of the form
\begin{equation}\label{exp.4}
u_x = py + q.
\end{equation}
Substituting \eqref{exp.4} into \eqref{exp.3}, we get
\[
p = \frac{\varepsilon}{\sigma_m^2},
\;\;
q = -\frac{x_T-x_0}{T}.
\]
Therefore,
\[
u_x(t,x)
= \frac{\varepsilon}{\sigma_m^2}(x-\mu_m(t))
-\frac{x_T-x_0}{T},
\]
and integrating with respect to \( x \), we obtain
\[
u(t,x)
= \frac{\varepsilon}{2\sigma_m^2}x^2
-\frac{\varepsilon\mu_m(t)}{\sigma_m^2}x
-\frac{x_T-x_0}{T}x.
\]
We also get
\[
u_t = -\frac{\varepsilon x(x_T-x_0)}{\sigma_m^2T},
\;\;
u_{xx} = \frac{\varepsilon}{\sigma_m^2}.
\]
The Wasserstein-2 distance between $w$ and $\nu$ is given by
\begin{equation}
W_2\big(\nu(t,x), w(t,x)\big)
=
\sqrt{
\big(\mu_w(t) - \mu_{\nu}\big)^2
+
\big(\sigma_w - \sigma_{\nu}\big)^2
}.
\end{equation}
The Kantorovich potential between the
Gaussian measures \( m \) and \( w \) satisfies
\[
\Phi(t,x) =
\frac{\delta}{\delta m}
\Bigg(
\dfrac{1}{2}\,
W_2^2\big(m(t,x),w(t,x)\big)
\Bigg),
\]
and has the explicit form
\[
\Phi(t,x)
= \frac{1}{2}\Big(1-\frac{\sigma_w}{\sigma_m}\Big)x^2
+ \Big(\frac{\sigma_w}{\sigma_m}\mu_m(t)-\mu_w(t)\Big)x.
\]
Finally, by choosing $V(t,x)$ to be we 
\[
V(t,x)
= \Phi(t,x)W_2(w,\nu) + u_t + \varepsilon\,u_{xx} - \frac{u_x^2}{2},
\]
we provide that the pair $(w,m)$ defined above in \eqref{exp.2} will solve the system \eqref{exp.1}.
This provides an explicit solution of the system in the one-dimensional case.
\begin{figure}[h]
\centering
\makebox[\textwidth][c]{%
    \hspace{-0.01cm}
    \includegraphics[scale=0.5]{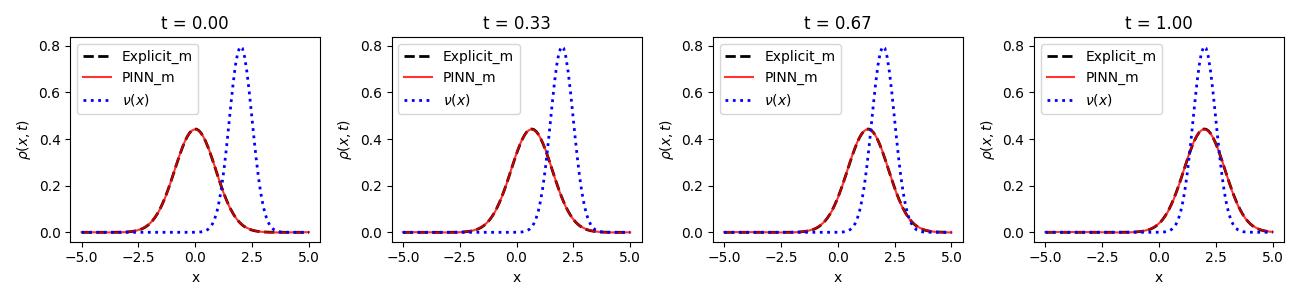}
}
\caption{Comparison of Explicit and PINN Solutions 
}
\label{fig:Explicit}
\end{figure}




Figure~\ref{fig:Explicit} shows the comparison between the explicit solution and the PINN approximation for the density $m(t,x)$. The black dashed and red curves almost coincide at all displayed time levels, which indicates that the PINN accurately captures the explicit solution.
The mean absolute error between the PINN approximation and the explicit solution remains below \(7\times10^{-4}\) at all evaluated time levels, confirming the accuracy of the numerical approximation.
In this example, the density keeps a fixed variance while its mean moves linearly from $x_0$ to $x_T$. The blue dotted curve represents the function $u(t,x)$.


\begin{figure*}[!htbp]

\centering
\vspace{-0.2cm} 

\subfloat{%
\includegraphics[width=0.33\textwidth]{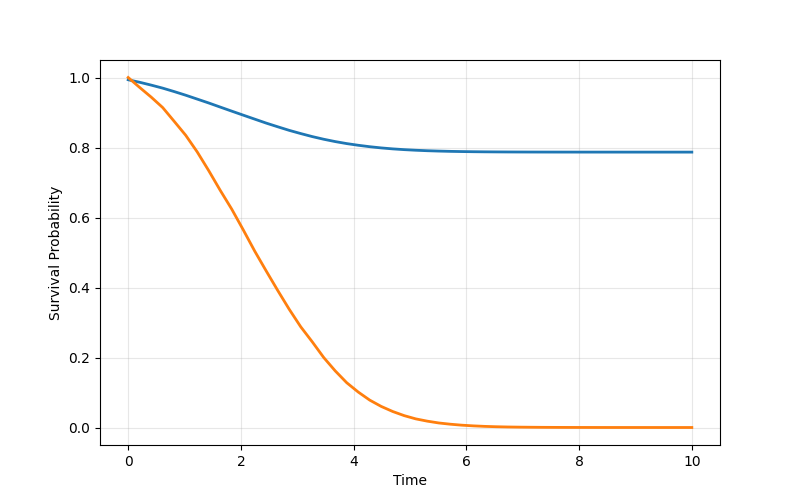}}
\hfill
{%
\includegraphics[width=0.33\textwidth]{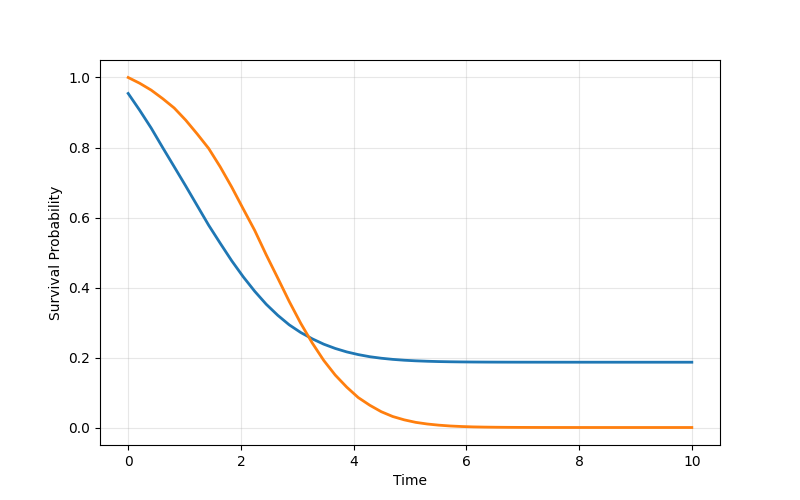}}
\hfill
{%
\includegraphics[width=0.33\textwidth]{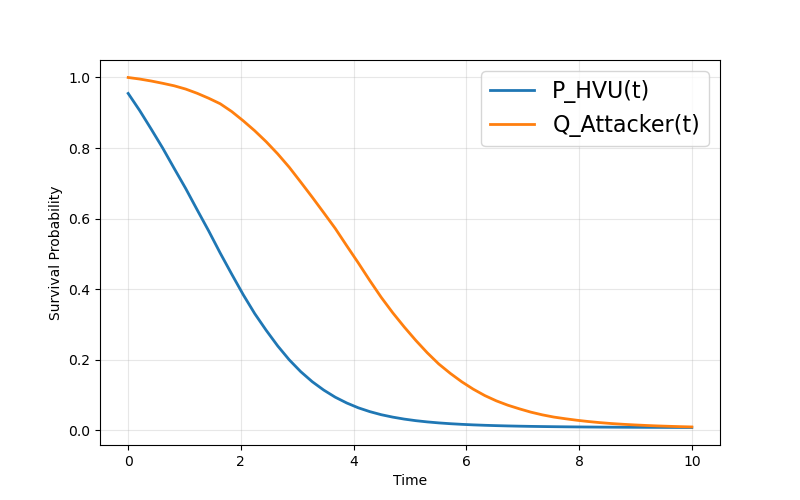}}

\caption{The blue curve represents the survival probability of the HVU, while the orange curve corresponds to the attackers.The initial position of the defenders is $(-3,-3)$  with variance $\sigma^2 = 0.85$.
The attackers start from $(-4,4)$ with variance $\sigma^2 = 0.85$ and move toward $(1,1)$, where the HVU is located with variance $\sigma^2 = 0.1$. For the first panel, the attrition rate parameters are $\sigma_A = \sigma_H = 5,\ \lambda_A = 14,\ \lambda_H = 1$. 
For the second panel, $\sigma_A = \sigma_H = 5,\ \lambda_A =\lambda_H = 7$. 
For the third panel, $\sigma_A = \sigma_H = 5,\ \lambda_A = 2,\ \lambda_H = 7$.}
\label{surv_won_lose}
\vspace{-0.5cm}
\end{figure*}


\subsection{Numerical Method and Simulation Setup}

A key challenge in numerically solving the proposed MFG system lies in
the computation of the Wasserstein-$2$ distance, which is
computationally demanding and must be evaluated at each time step. To
mitigate this difficulty, we approximate the Wasserstein-$2$ distance
using the Sinkhorn algorithm, see \cite{cuturi2013sinkhorn}. The resulting problem is then solved using
physics-informed neural networks (PINNs), see \cite{pinn}.

We consider a two-dimensional setting, with spatial domain
$\Omega = [-5,5]^2$ and time horizon $T=10$.  The spatial domain is discretized using a uniform
$60 \times 60$ grid, while the time interval is sampled uniformly with
$50$ points. These collocation points are used to train the PINN model. In our model, we set the
diffusion coefficient to $\varepsilon = 0.001$ and the Pareto weight to
$\alpha = 0.1$.

 As mentioned before, to numerically solve the mean-field game system,
we first approximate the MFG system \eqref{eq:MFG} by replacing the
Wasserstein-$2$ distance with its Sinkhorn approximation. We then solve
the resulting approximated problem using PINNs. In this approach, the loss function is built from the residuals of the coupled Hamilton–Jacobi and Fokker–Planck equations. It also includes the terminal and initial conditions of the problem.
We use two separate neural networks: one for the value function 
$w(t,x,y)$ and one for the density $m(t,x,y)$. Both networks have the same ResNet architecture. It consists of an input layer with three neurons, followed by a fully connected hidden layer with 128 neurons and Tanh activation function. This is followed by four residual blocks, each composed of two fully connected layers with 128 neurons, Tanh activation functions, and a residual connection. The network ends with a final linear output layer with one neuron. For the density network, we apply Softplus activation function. 
We do not use labeled data for training. Instead, the networks learn by minimizing the equation errors. Trainings were conducted for 15,000 epochs with the Adam optimizer, utilizing NVIDIA GeForce RTX 4080 and 5080 GPUs.

\vspace{-0.2cm}
\subsection{Simulation Scenarios}
We consider four groups of simulations. First, we vary the attrition-rate
parameters in order to study how different weapon strengths affect the
survival probabilities of the HVU and the attackers. Second, we vary the
initial spread of the defender population to show how population
dispersion influences the effectiveness of defense. Third, we study the
effect of the parameter $\alpha$, which controls the cost of motion in
the optimization problem. Finally, we consider scenarios with several
attacker groups and several HVUs, in order to illustrate how the defender
population redistributes itself in more complex configurations.

Fig.~\ref{surv_won_lose} illustrates the impact of the attrition rate
parameters on the survival probabilities of the HVU and the attackers.
In all three scenarios (first, second, and third panels), the defenders,
attackers, and HVU share the same initial positions and variances, and
the attacker dynamics are identical. The only difference lies in the
attrition rate parameters. 
In the first panel, the defenders are stronger than the attackers; in
the second, both sides have equal strength; and in the third, the
attackers are stronger than the defenders. This results in high,
moderate, and low survival probabilities for the HVU population, respectively, and
conversely for the attackers.

\begin{figure*}[!htbp]
\centering
\subfloat{%
\includegraphics[width=0.33\textwidth]{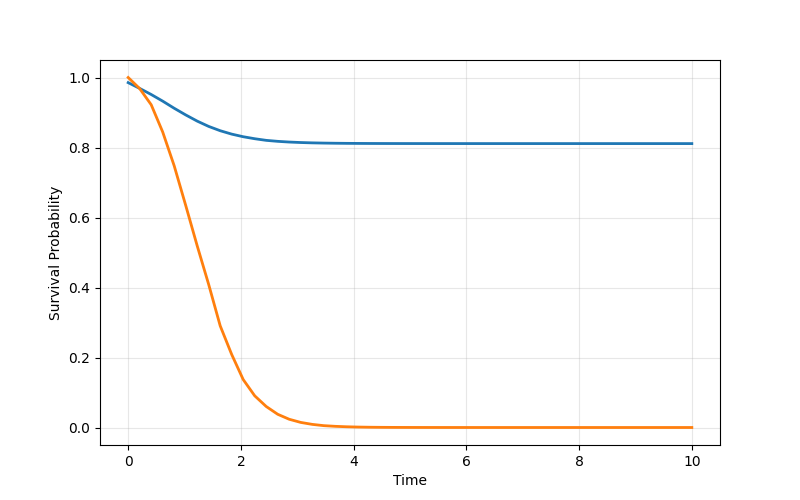}}
\hfill
\subfloat{%
\includegraphics[width=0.33\textwidth]{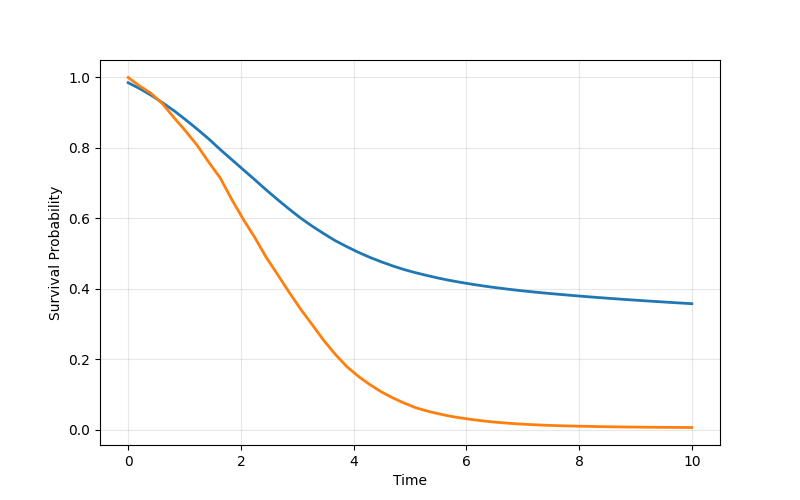}}
\hfill
\subfloat{%
\includegraphics[width=0.33\textwidth]{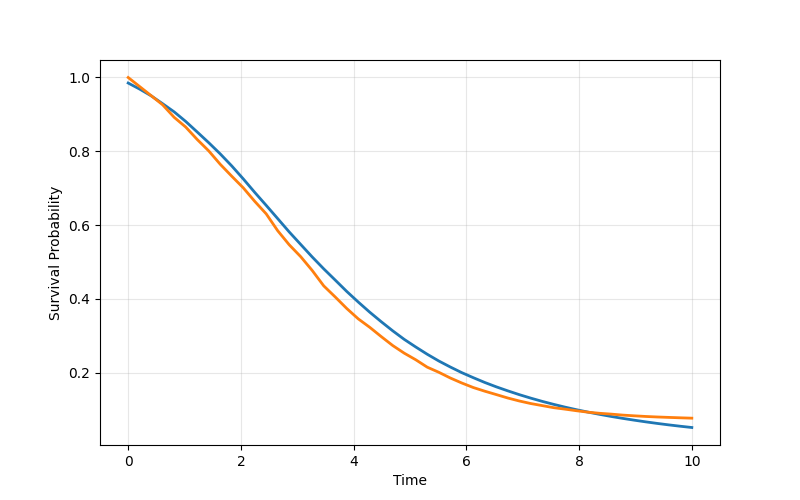}}
\hfill
\caption{Comparis  on of survival probabilities for different $\alpha$ values. 
The first, second, and third panels correspond to $\alpha = 0.1$, 
$\alpha = 0.5$, and $\alpha = 0.7$, respectively. 
The attrition parameters and initial conditions are identical across all scenarios. 
The attrition rate parameters are 
$\sigma_A = \sigma_H = 2$, $\lambda_A = \lambda_H = 5$.}
\label{surv_alpha}
\vspace{-0.4cm}
\end{figure*}

\begin{figure*}[!htbp]
\vspace{-0cm}
\centering

\subfloat[\centering Top-Down View $t = 0.00$]{%
\includegraphics[width=0.25\textwidth]{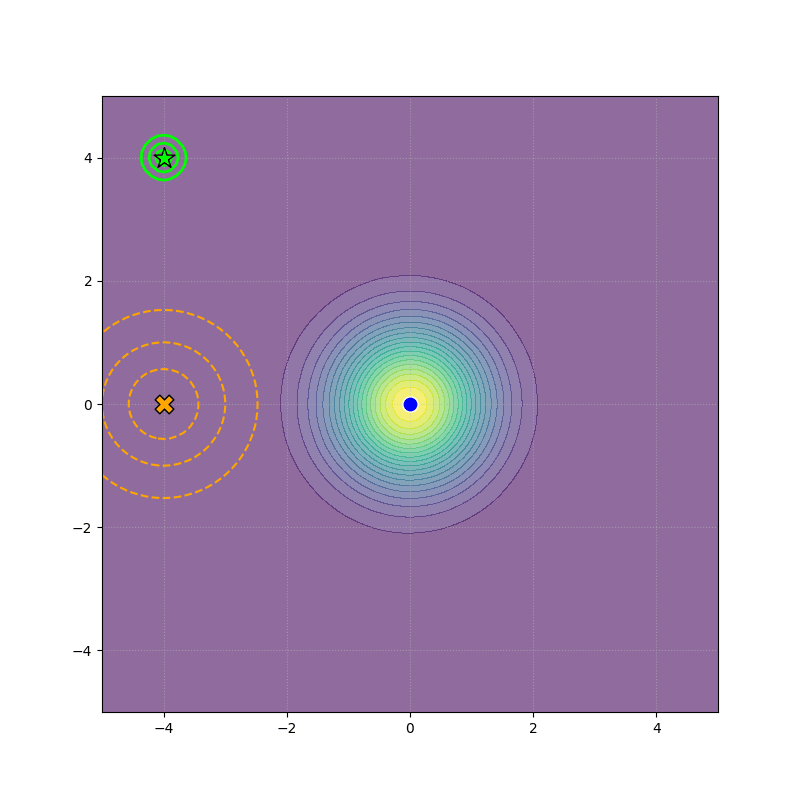}}
\hfill
\subfloat[\centering Top-Down View $t = 3.27$]{%
\includegraphics[width=0.25\textwidth]{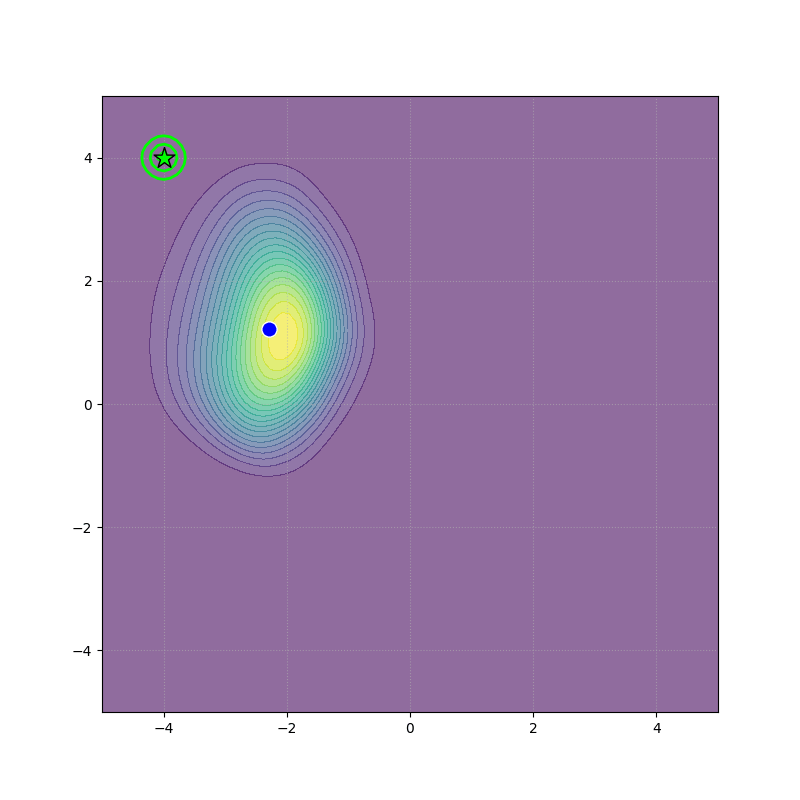}}
\hfill
\subfloat[\centering Top-Down View $t = 6.53$]{%
\includegraphics[width=0.25\textwidth]{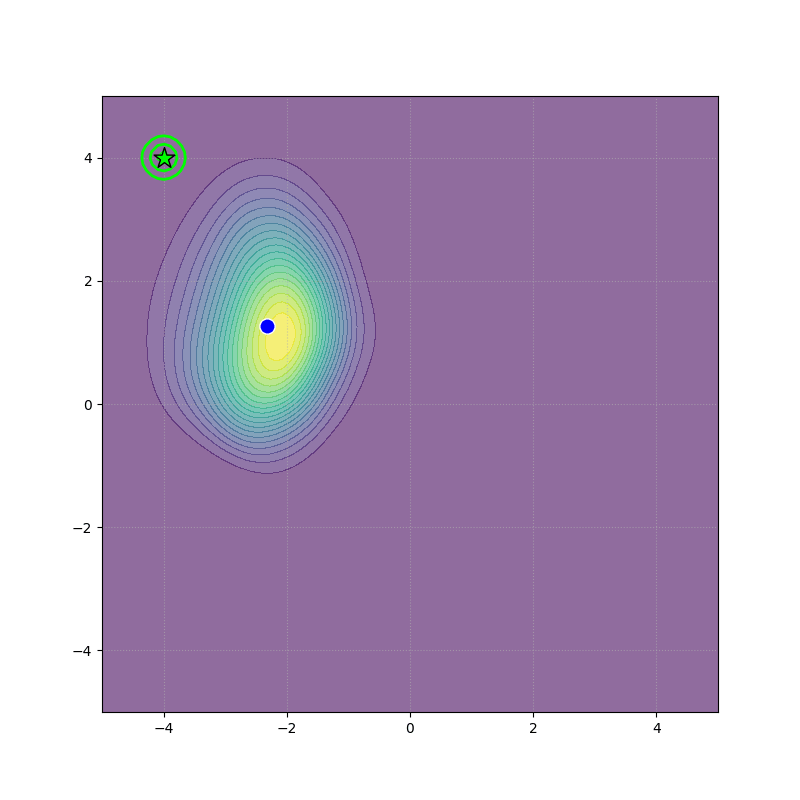}}
\hfill
\subfloat[\centering Top-Down View $t = 10.00$]{%
\includegraphics[width=0.25\textwidth]{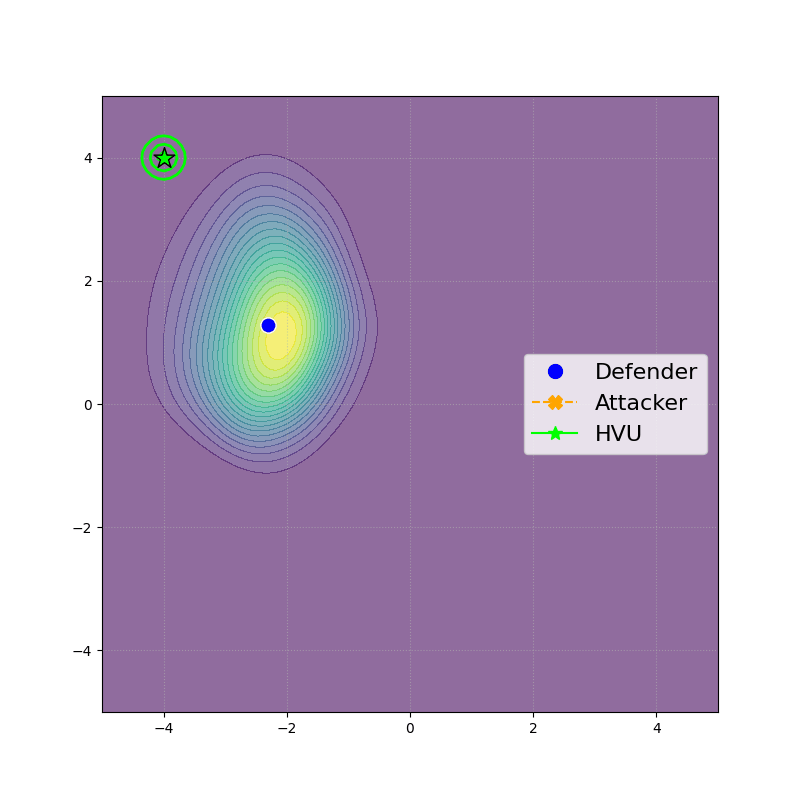}}

\caption{Top-down view of the evolution of the attacker and defender over time
for the first scenario ($\alpha = 0.1$) in Fig.~\ref{surv_alpha}.}
\label{alfa=0.1}
\end{figure*}

\begin{figure*}[!htbp]
\vspace{-0.4cm}
\centering
\subfloat[\centering Top-Down View $t = 0.00$]{%
\includegraphics[width=0.25\textwidth]{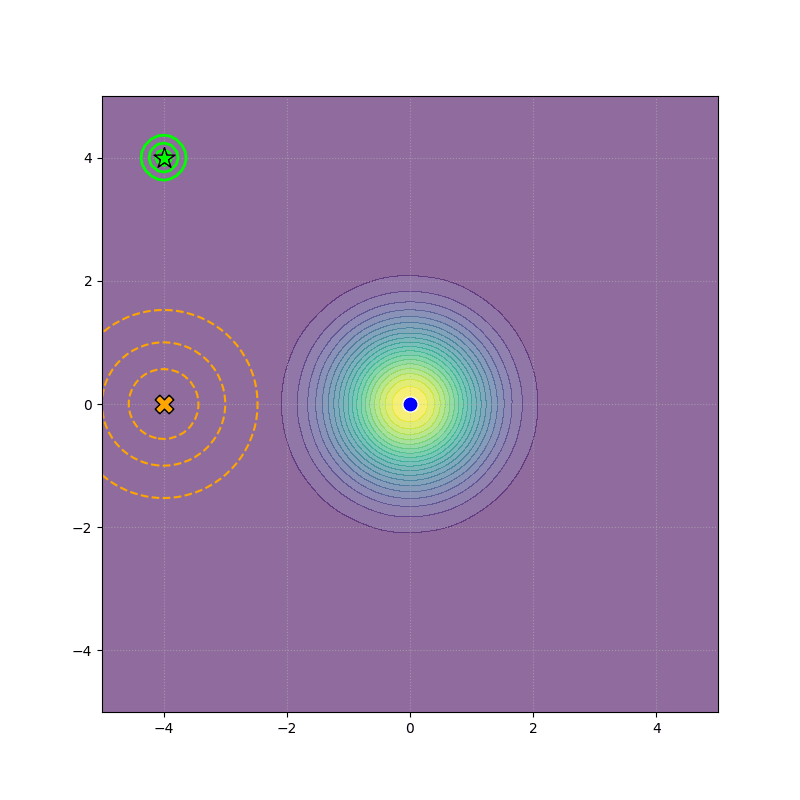}}
\hfill
\subfloat[\centering Top-Down View $t = 3.27$]{%
\includegraphics[width=0.25\textwidth]{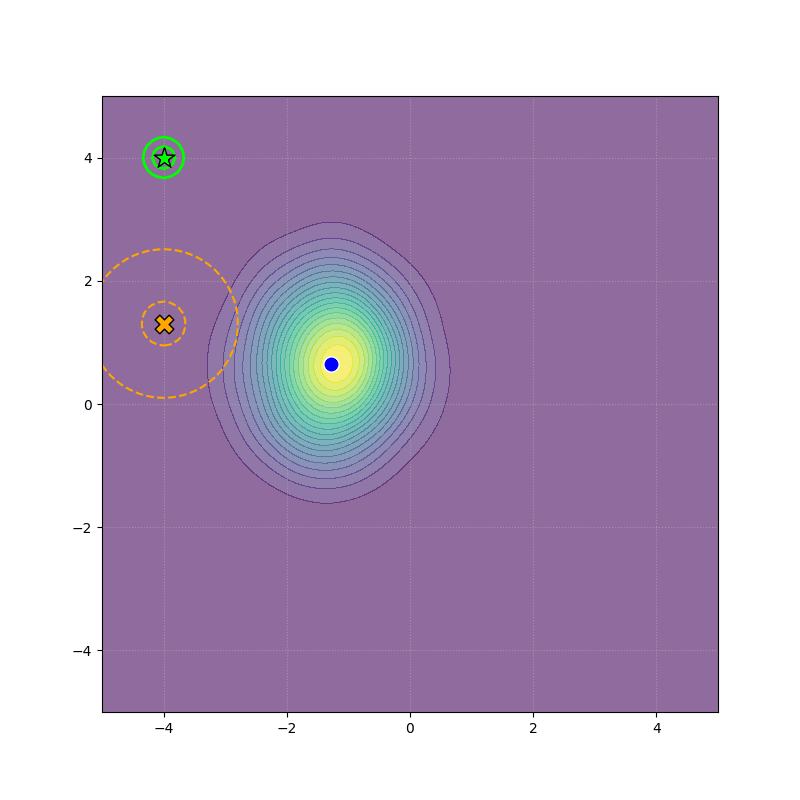}}
\hfill
\subfloat[\centering Top-Down View $t = 6.53$]{%
\includegraphics[width=0.25\textwidth]{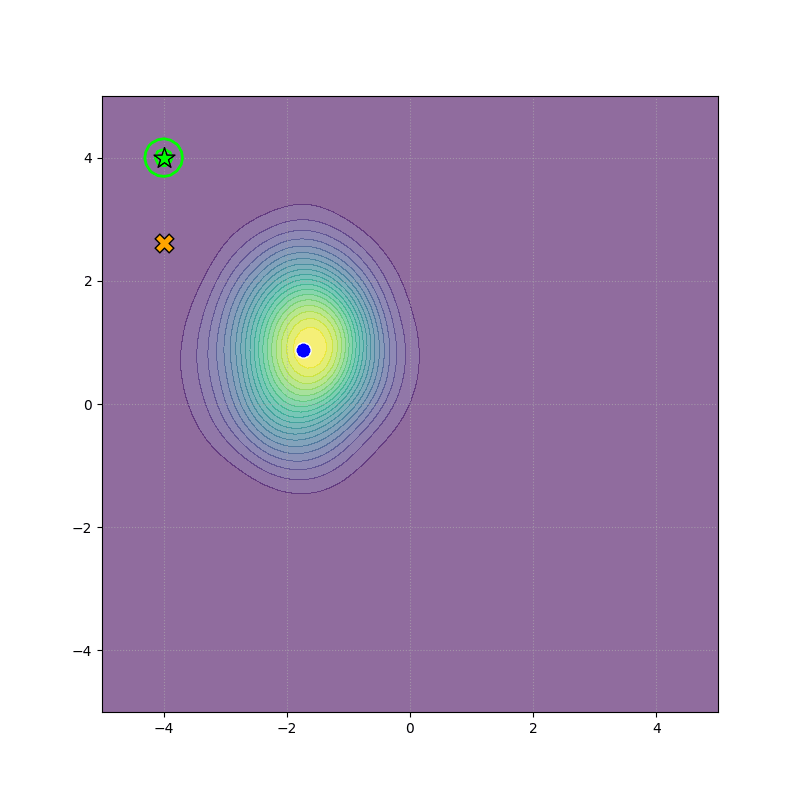}}
\hfill
\subfloat[\centering Top-Down View $t = 10.00$]{%
\includegraphics[width=0.25\textwidth]{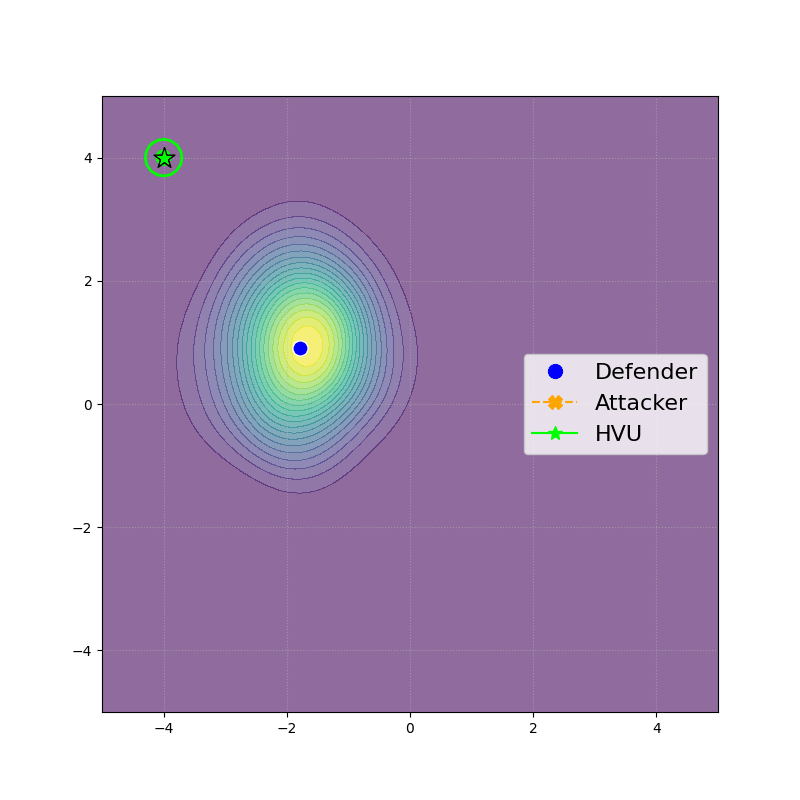}}

\caption{Top-down view of the evolution of the attacker and defender over time for the second scenario ($\alpha = 0.5$) in Fig.~\ref{surv_alpha}.}
\label{alfa=0.5}
\end{figure*}

\begin{figure*}[!htbp]
\vspace{-0.4cm}
\centering
\subfloat[\centering Top-Down View $t = 0.00$]{%
\includegraphics[width=0.25\textwidth]{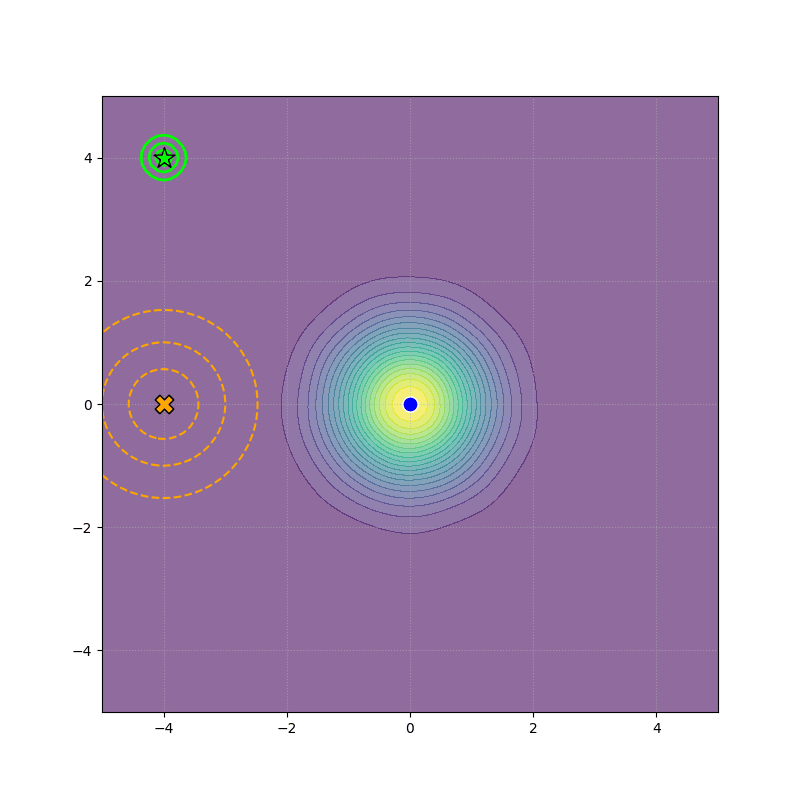}}
\hfill
\subfloat[\centering Top-Down View $t = 3.27$]{%
\includegraphics[width=0.25\textwidth]{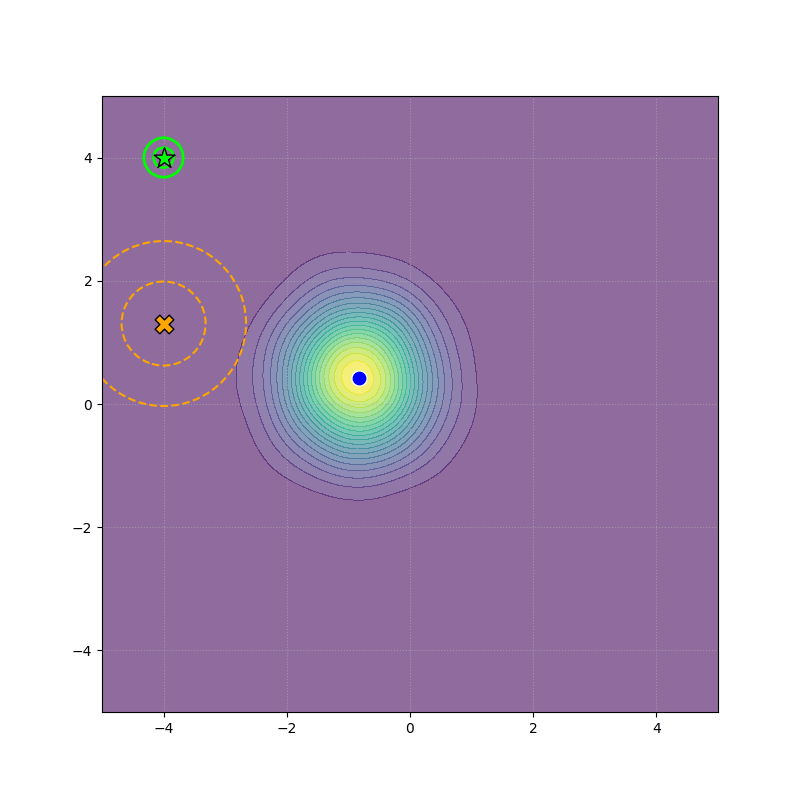}}
\hfill
\subfloat[\centering Top-Down View $t = 6.53$]{%
\includegraphics[width=0.25\textwidth]{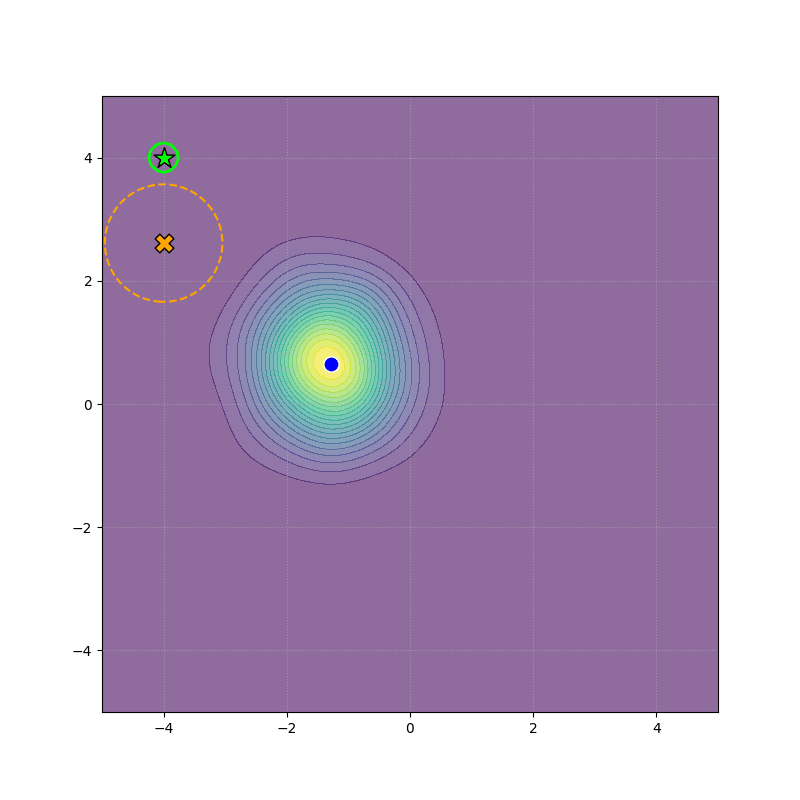}}
\hfill
\subfloat[\centering Top-Down View $t = 10.00$]{%
\includegraphics[width=0.25\textwidth]{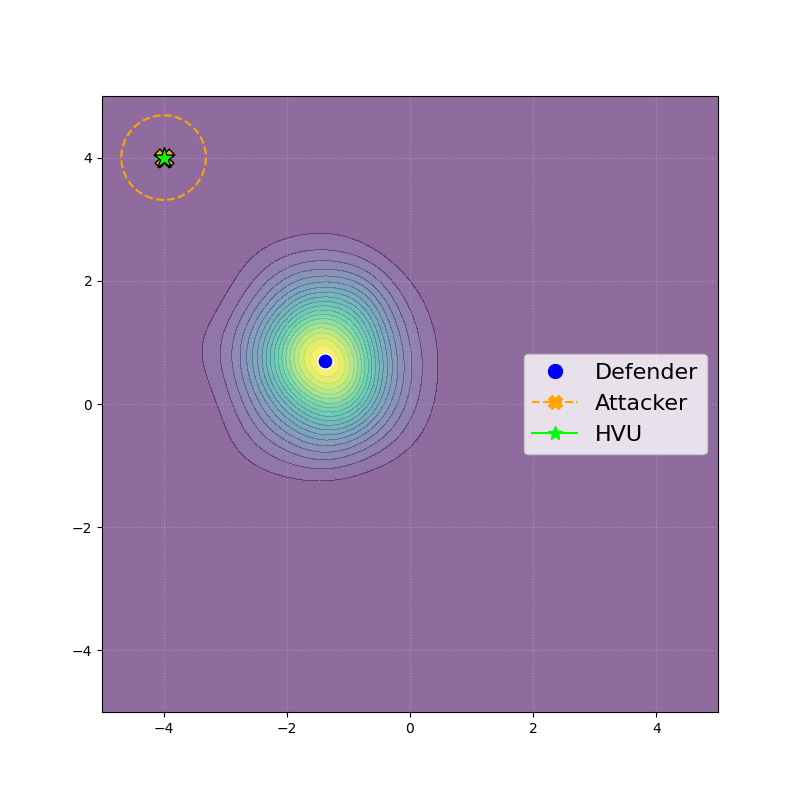}}

\caption{Top-down view of the evolution of the attacker and defender over time for the third scenario ($\alpha = 0.7$) in Fig.~\ref{surv_alpha}.}
\label{alfa=0.7}
\end{figure*}

From Fig.~\ref{surv_alpha}, we observe the differences in survival probabilities for various values of $\alpha$. When $\alpha$ is small ($\alpha = 0.1$, see Fig.~\ref{alfa=0.1}), the defenders move faster and can better protect the HVU, resulting in a higher survival probability. 
For larger values of $\alpha$, the defenders react more slowly, and the survival probability of the HVU decreases. The case $\alpha = 0.5$ (Fig.~\ref{alfa=0.5}) represents a moderate scenario, while $\alpha = 0.7$ (Fig.~\ref{alfa=0.7}) leads to worse performance compared to $\alpha = 0.5$.

\begin{figure*}[!htbp]
\vspace{-0.4cm}
\centering

\subfloat[\centering Top-Down View $t = 0.00$]{%
\includegraphics[width=0.25\textwidth]{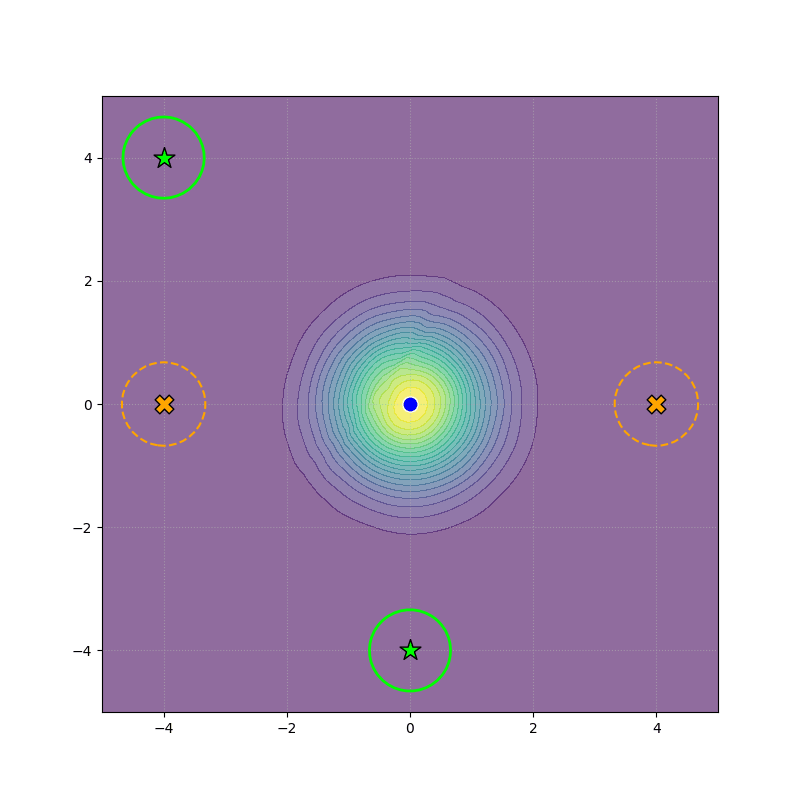}}
\hfill
\subfloat[\centering Top-Down View $t = 1.6$]{%
\includegraphics[width=0.25\textwidth]{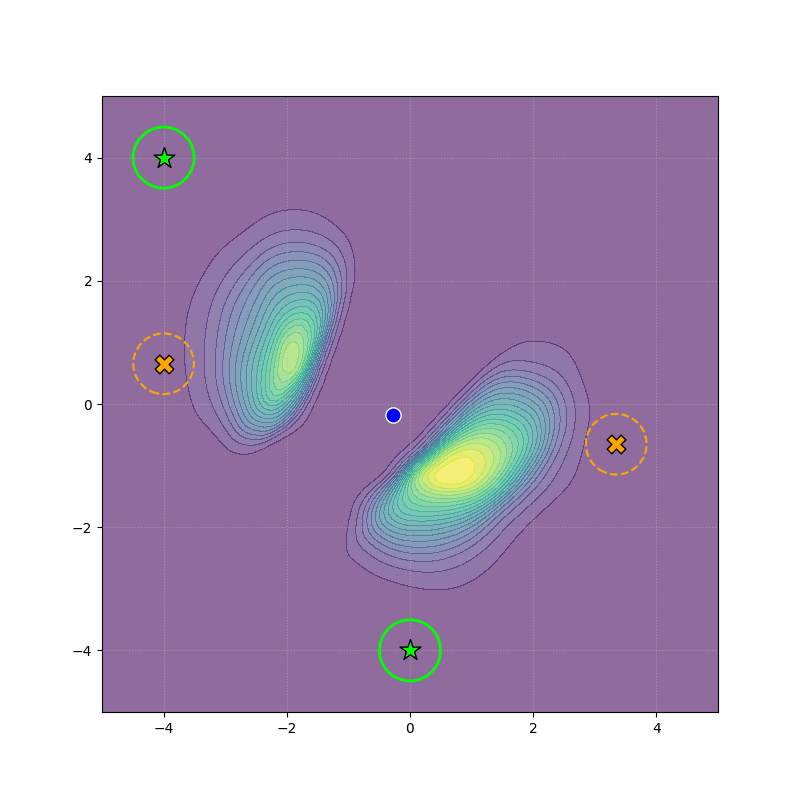}}
\hfill
\subfloat[\centering Top-Down View $t = 3.2$]{%
\includegraphics[width=0.25\textwidth]{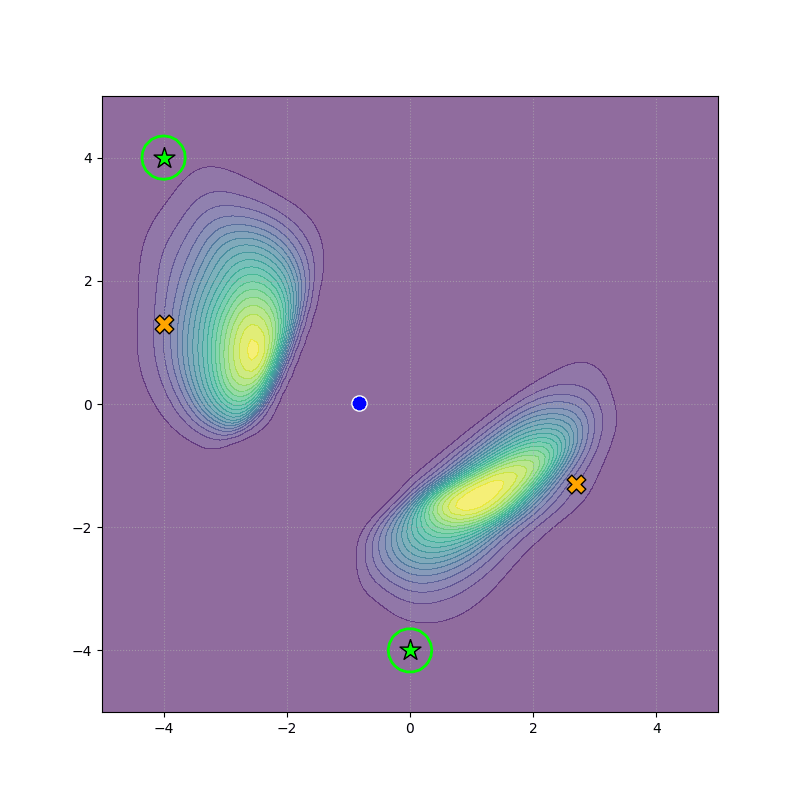}}
\hfill
\subfloat[\centering Top-Down View $t = 10.00$]{%
\includegraphics[width=0.25\textwidth]{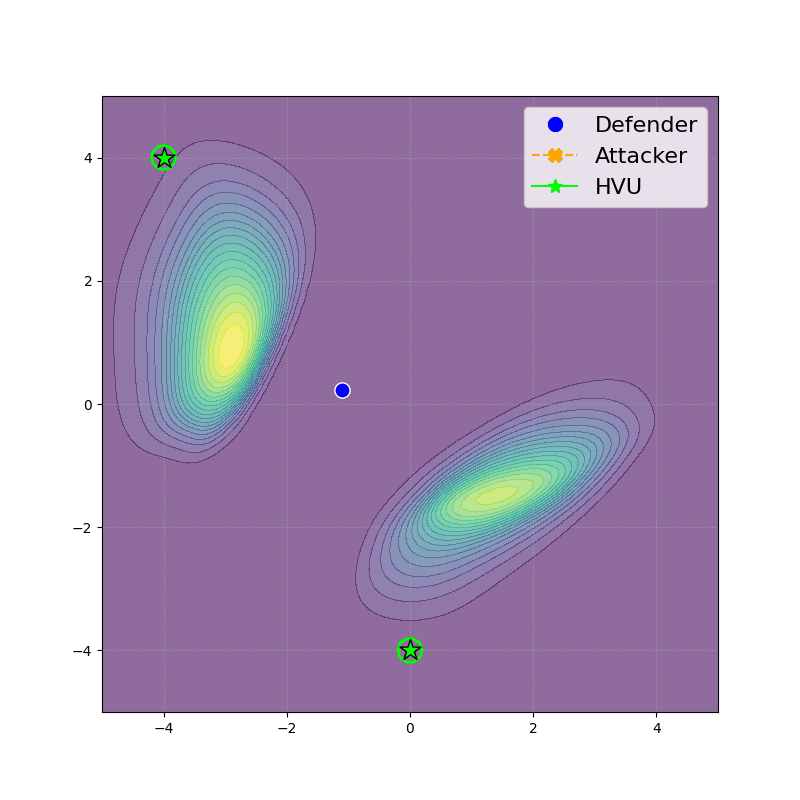}}

\caption{Top-down view of the evolution of the attacker and defender over time.There are two attacker groups and two HVUs. To protect them better, the defenders split into two parts, and each part moves to defend one HVU.
The defenders start from $(0,0)$ with variance $\sigma^2 = 0.85$.
The first group of attackers starts from $(-4,0)$ and moves toward the first HVU at $(-4,4)$. The second group starts from $(4,0)$ and moves toward the second HVU at $(0,-4)$.
The attackers and The HVUs have variance $\sigma^2 = 0.5$.Since there are two HVUs, the total mass is split equally, so each one has mass $\frac{1}{2}$. The attrition rate parameters are
$\sigma_A=\sigma_H=5$, $\lambda_A=1$, $\lambda_H=5$ and $\alpha=0.1$.}
\label{split}
\end{figure*}

\begin{figure*}[!htbp]
\vspace{-0.4cm}
\centering

\subfloat[\centering Top-Down View $t = 0.00$]{%
\includegraphics[width=0.25\textwidth]{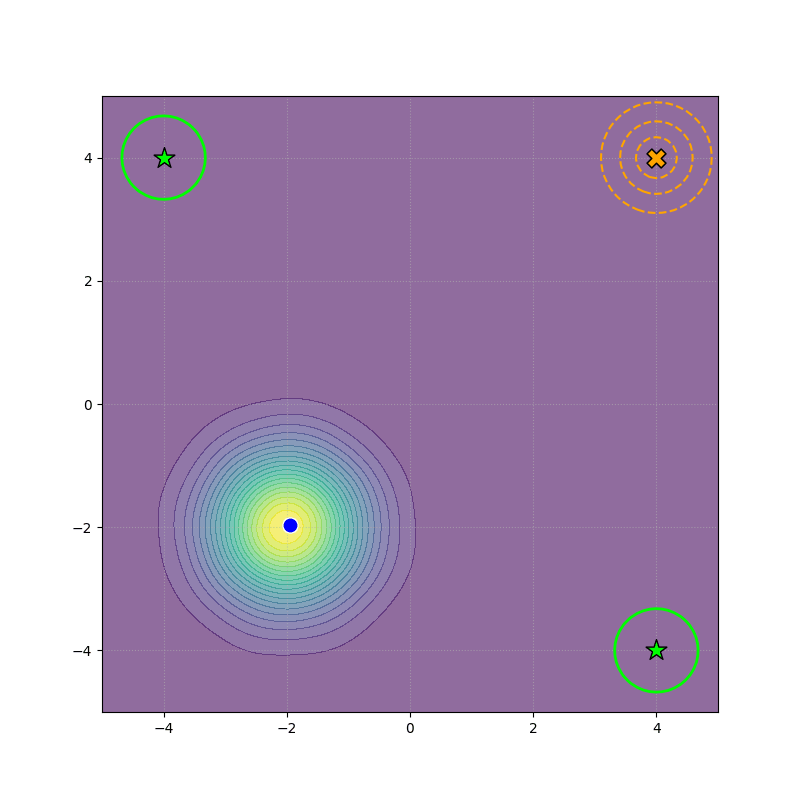}}
\hfill
\subfloat[\centering Top-Down View $t = 1.6$]{%
\includegraphics[width=0.25\textwidth]{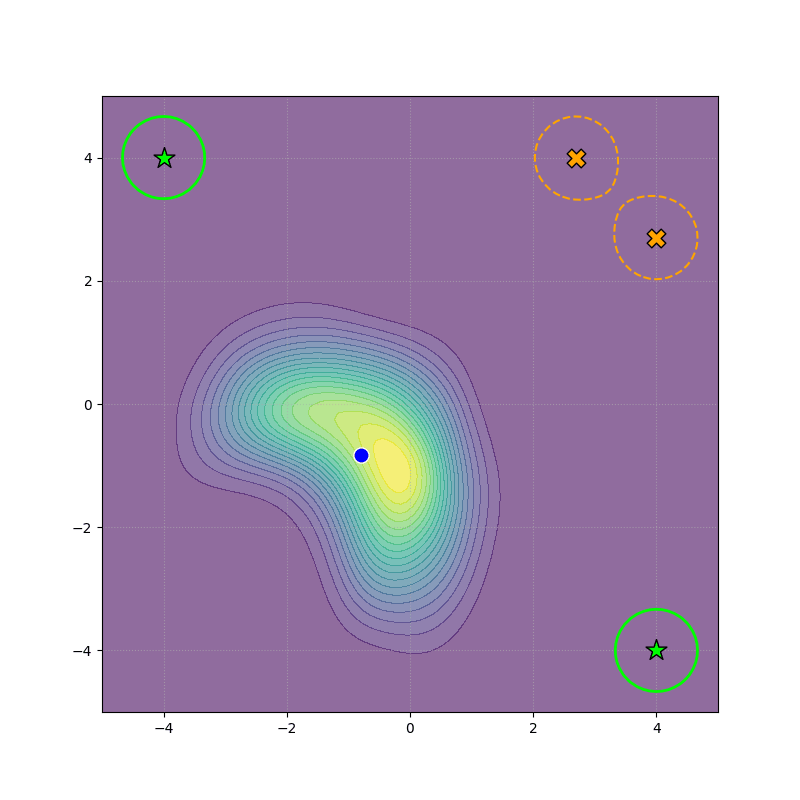}}
\hfill
\subfloat[\centering Top-Down View $t = 3.2$]{%
\includegraphics[width=0.25\textwidth]{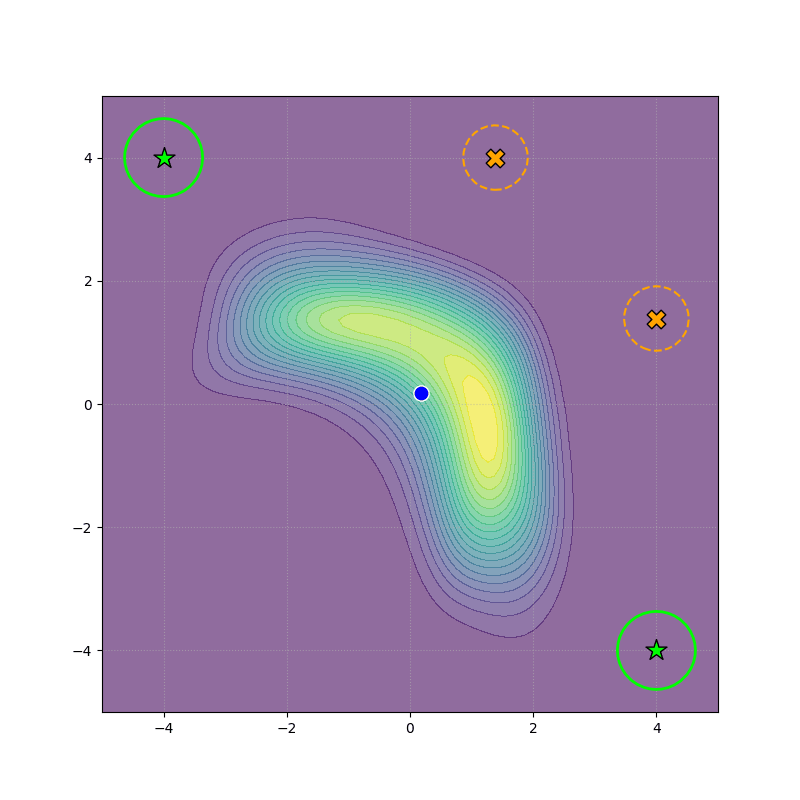}}
\hfill
\subfloat[\centering Top-Down View $t = 10.00$]{%
\includegraphics[width=0.25\textwidth]{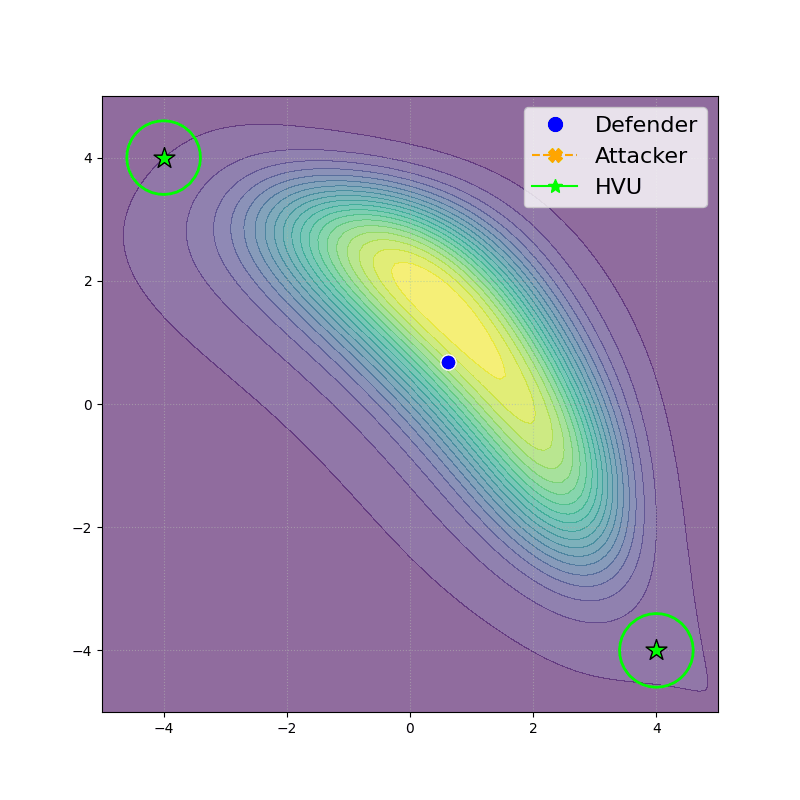}}

\caption{Top-down view of the evolution of the attacker and defender over time.There are two attacker groups and two HVUs.The defenders gradually spread out over time to intercept both attacker groups.
The defenders start from $(-2,-2)$ with variance $\sigma^2 = 0.85$.
Both attacker groups start from $(4,4)$; one moves toward the first HVU at $(-4,4)$, while the other moves toward the second HVU at $(4,-4)$.The attackers and The HVUs have variance $\sigma^2 = 0.5$.Since there are two HVUs, the total mass is split equally, so each one has mass $\frac{1}{2}$.The attrition rate parameters are
$\sigma_A=\sigma_H=5$,$\lambda_A=3$, $\lambda_H=5$ and $\alpha=0.1$}
\label{two_att}
\end{figure*}

\begin{figure*}[!htbp]
\vspace{-0.4cm}
\centering

\subfloat[\centering Top-Down View $t = 0.00$]{%
\includegraphics[width=0.25\textwidth]{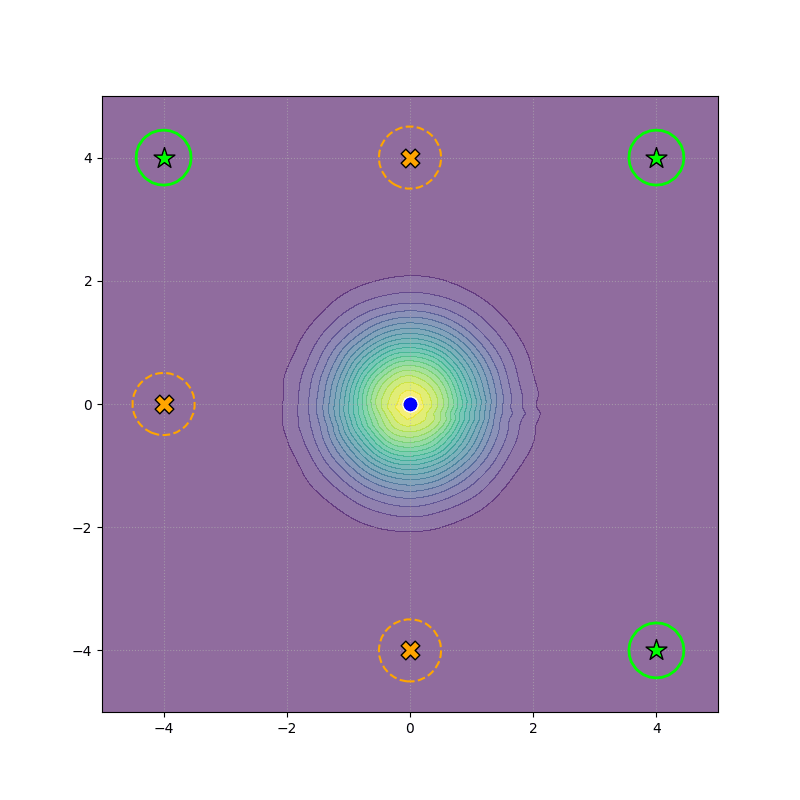}}
\hfill
\subfloat[\centering Top-Down View $t = 1.6$]{%
\includegraphics[width=0.25\textwidth]{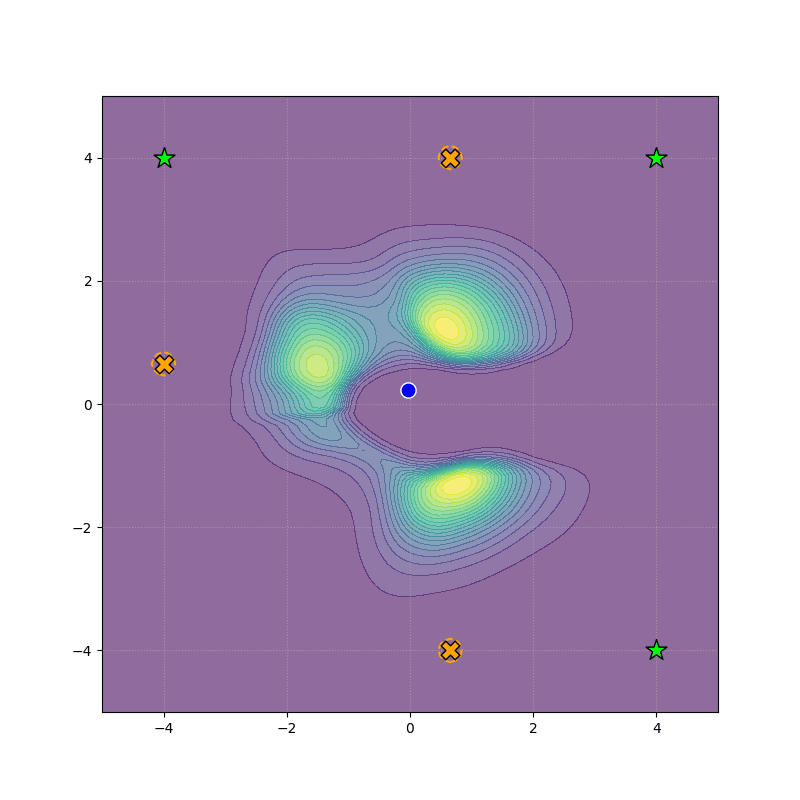}}
\hfill
\subfloat[\centering Top-Down View $t = 3.2$]{%
\includegraphics[width=0.25\textwidth]{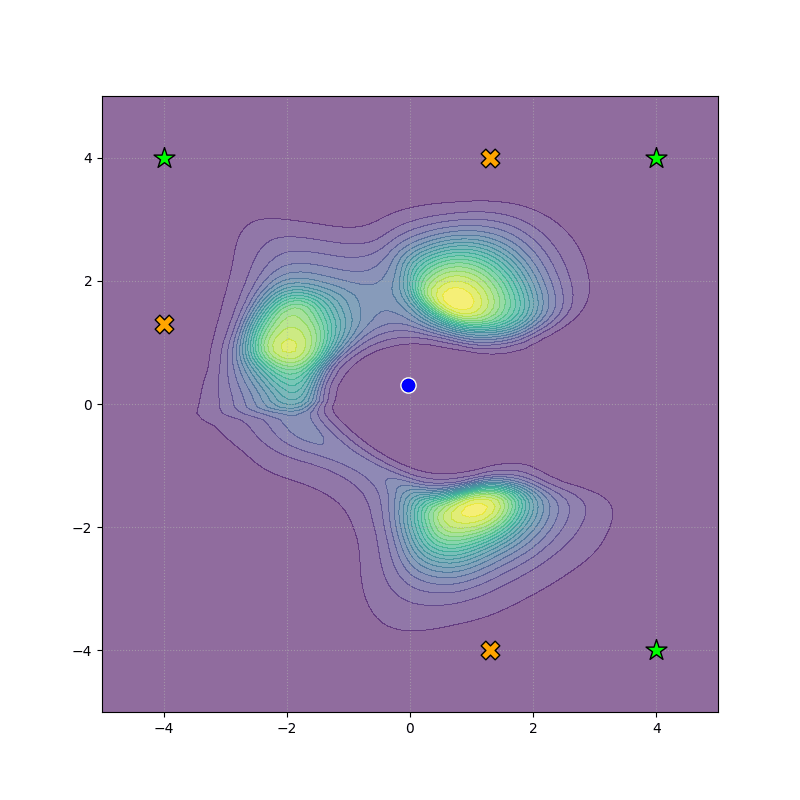}}
\hfill
\subfloat[\centering Top-Down View $t = 10.00$]{%
\includegraphics[width=0.25\textwidth]{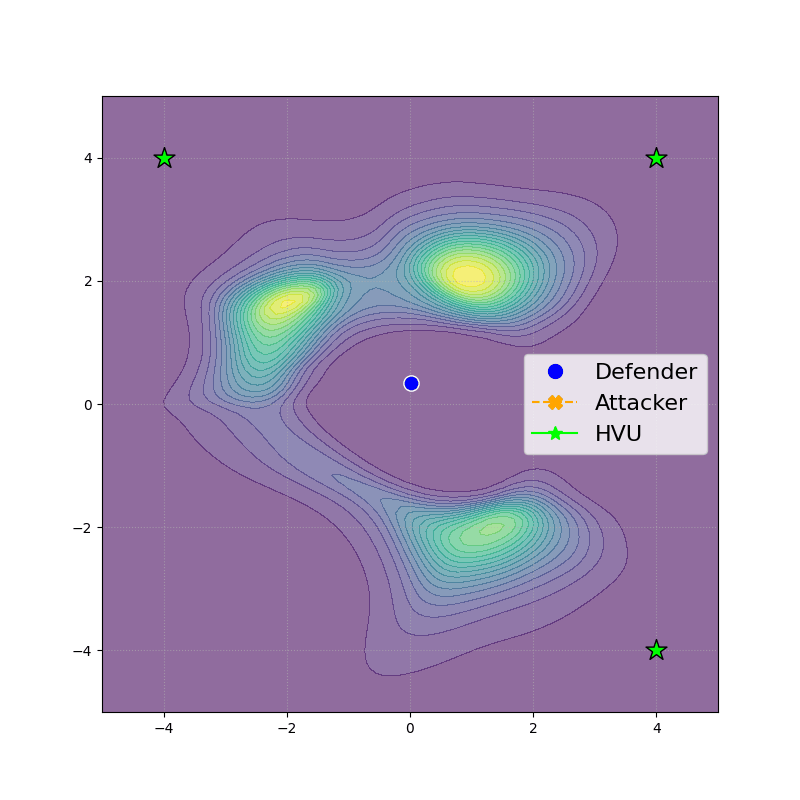}}

\caption{Top-down view of the evolution of the attacker and defender over time.There are three attacker groups and three HVUs. To protect them better, the defenders split into three parts, and each part moves to defend one HVU.
The defenders start from $(0,0)$ with variance $\sigma^2 = 0.85$.
The first group of attackers starts from $(-4,0)$ and moves toward the first HVU at $(-4,4)$. The second group starts from $(0,4)$ and moves toward the second HVU at $(4,4)$.
The third group starts from $(0,-4)$ and moves toward the third HVU at $(4,-4)$.
The attackers and The HVUs have variance $\sigma^2 = 0.5$. The attrition rate parameters are
$\sigma_A=\sigma_H=5$, $\lambda_A=1$, $\lambda_H=5$ and $\alpha=0.1$.}
\label{three_att}
\end{figure*}

\begin{figure*}[!htbp]
\centering
\subfloat{%
\includegraphics[width=0.33\textwidth]{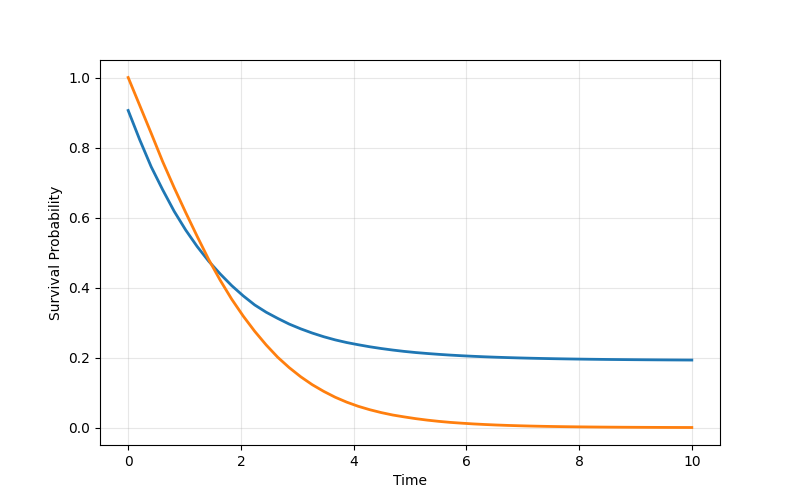}}
\hfill
\subfloat{%
\includegraphics[width=0.33\textwidth]{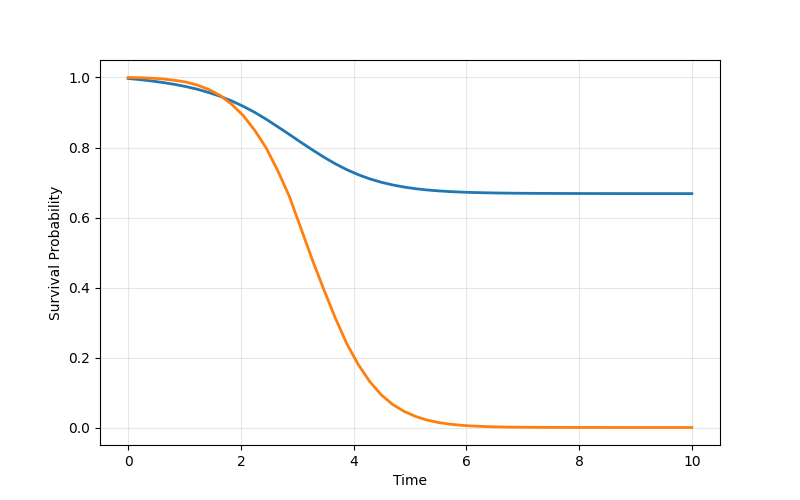}}
\hfill
\subfloat{%
\includegraphics[width=0.33\textwidth]{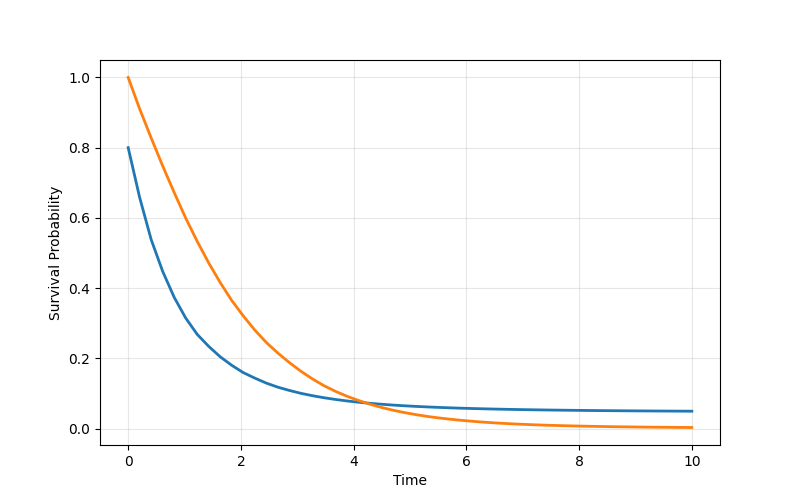}}
\hfill
\caption{Survival probabilities of the attackers and the HVUs. The first panel shows the survival probability for the scenario in Fig.~\ref{split}. The second panel shows the survival probability for the scenario in Fig.~\ref{two_att}. The third panel shows the survival probability for the scenario in Fig.~\ref{three_att}} 
\label{surv_split}
\vspace{-0.6cm}
\end{figure*}

In Figs.~\ref{split}, \ref{two_att}, and \ref{three_att}, we present several numerical simulations of the attacker--defender model with different numbers of attacker groups and HVUs. In all scenarios, the defenders dynamically redistribute themselves in order to intercept the attackers and protect the HVUs. The figures illustrate how the defender density evolves over time depending on the number of targets. Fig~\ref{surv_split} presents the corresponding survival probabilities of the attackers and the HVUs for the scenarios shown in Figs.~\ref{split}, \ref{two_att}, and \ref{three_att}.


\bibliographystyle{acm}
\bibliography{bibl}

@article{achdou2019mean,
  title={Mean field games},
  author={Achdou, Yves and Cardaliaguet, Pierre and Delarue, Fran{\c{c}}ois and Porretta, Alessio and Santambrogio, Filippo},
  journal={Lecture Notes in Math},
  volume={2281},
  year={2019},
  publisher={Springer}
}

@article{krylov1996lectures,
  title={Lectures on Elliptic and Parabolic Equations in H older spaces},
  author={Krylov, NV},
  journal={American Mathematical Society},
  year={1996}
}

@article{jabin2017mean,
  title={Mean field limit for stochastic particle systems},
  author={Jabin, Pierre-Emmanuel and Wang, Zhenfu},
  journal={Active Particles, Volume 1: Advances in Theory, Models, and Applications},
  pages={379--402},
  year={2017},
  publisher={Springer}
}

@article{carrillo2014derivation,
  title={The derivation of swarming models: mean-field limit and Wasserstein distances},
  author={Carrillo, Jos{\'e} Antonio and Choi, Young-Pil and Hauray, Maxime},
  journal={Collective Dynamics from Bacteria to Crowds: An Excursion Through Modeling, Analysis and Simulation},
  pages={1--46},
  year={2014},
  publisher={Springer}
}

@article{topaz2004swarming,
  title={Swarming patterns in a two-dimensional kinematic model for biological groups},
  author={Topaz, Chad M and Bertozzi, Andrea L},
  journal={SIAM Journal on Applied Mathematics},
  volume={65},
  number={1},
  pages={152--174},
  year={2004},
  publisher={SIAM}
}

@article{carrillo2010particle,
  title={Particle, kinetic, and hydrodynamic models of swarming},
  author={Carrillo, Jos{\'e} A and Fornasier, Massimo and Toscani, Giuseppe and Vecil, Francesco},
  journal={Mathematical modeling of collective behavior in socio-economic and life sciences},
  pages={297--336},
  year={2010},
  publisher={Springer}
}

@article{walton2022defense,
  title={Defense against adversarial swarms with parameter uncertainty},
  author={Walton, Claire and Kaminer, Isaac and Gong, Qi and Clark, Abram H and Tsatsanifos, Theodoros},
  journal={Sensors},
  volume={22},
  number={13},
  pages={4773},
  year={2022},
  publisher={MDPI}
}

@inproceedings{tsatsanifos2021,
  title={Modeling large-scale adversarial swarm engagements using optimal control},
  author={Tsatsanifos, Theodoros and Clark, Abram H and Walton, Claire and Kaminer, Isaac and Gong, Qi},
  booktitle={2021 60th IEEE Conference on Decision and Control (CDC)},
  pages={1244--1249},
  year={2021},
  organization={IEEE}
}

@article{cuturi2013sinkhorn,
  title={Sinkhorn distances: Lightspeed computation of optimal transport},
  author={Cuturi, Marco},
  journal={Advances in neural information processing systems},
  volume={26},
  year={2013}
}

@article{lasry2007mean,
  title={Mean field games},
  author={Lasry, Jean-Michel and Lions, Pierre-Louis},
  journal={Japanese journal of mathematics},
  volume={2},
  number={1},
  pages={229--260},
  year={2007},
  publisher={Springer}
}

@article{walton2018optimal,
  title={Optimal motion planning in rapid-fire combat situations with attacker uncertainty},
  author={Walton, Claire and Lambrianides, Panos and Kaminer, Isaac and Royset, Johannes and Gong, Qi},
  journal={Naval Research Logistics (NRL)},
  volume={65},
  number={2},
  pages={101--119},
  year={2018},
  publisher={Wiley Online Library}
}

@book{villani2021topics,
  title={Topics in optimal transportation},
  author={Villani, C{\'e}dric},
  volume={58},
  year={2021},
  publisher={American Mathematical Soc.}
}

@incollection{villani2009wasserstein,
  title={The wasserstein distances},
  author={Villani, C{\'e}dric},
  booktitle={Optimal transport: old and new},
  pages={93--111},
  year={2009},
  publisher={Springer}
}

@book{santambrogio2015optimal,
  title={Optimal transport for applied mathematicians},
  author={Santambrogio, Filippo},
  year={2015},
  publisher={Springer}
}

@inproceedings{kang2025online,
  title={Online Defensive Motion Planning Against Adversarial Swarm Attacks using Bernstein Polynomials-Based Model Predictive Control},
  author={Kang, Hyungsoo and Aoun, Christoph and Kaminer, Isaac and Hovakimyan, Naira},
  booktitle={2025 International Conference on Unmanned Aircraft Systems (ICUAS)},
  pages={222--227},
  year={2025},
  organization={IEEE}
}

@article{caffarelli1992,
  title={The regularity of mappings with a convex potential},
  author={Caffarelli, Luis A},
  journal={Journal of the American Mathematical Society},
  volume={5},
  number={1},
  pages={99--104},
  year={1992},
  publisher={JSTOR}
}

@article{nirenberg1959,
  title={On elliptic partial differential equations},
  author={Nirenberg, Louis},
  journal={Annali della Scuola Normale Superiore di Pisa-Scienze Fisiche e Matematiche},
  volume={13},
  number={2},
  pages={115--162},
  year={1959}
}

@article{kitagawa2025,
  title={Stability of optimal transport maps on Riemannian manifolds},
  author={Kitagawa, Jun and Letrouit, Cyril and M{\'e}rigot, Quentin},
  journal={arXiv preprint arXiv:2504.05412},
  year={2025}
}

@article{delalande2023,
  title={Quantitative stability of optimal transport maps under variations of the target measure},
  author={Delalande, Alex and Merigot, Quentin},
  journal={Duke Mathematical Journal},
  volume={172},
  number={17},
  pages={3321--3357},
  year={2023},
  publisher={Duke University Press}
}

@article{gallouet2025,
  title={Strong c-concavity and stability in optimal transport},
  author={Gallou{\"e}t, Anatole and M{\'e}rigot, Quentin and Thibert, Boris},
  journal={Journal de Math{\'e}matiques Pures et Appliqu{\'e}es},
  pages={103773},
  year={2025},
  publisher={Elsevier}
}

@article{letrouit2025,
  title={Quantitative stability of optimal transport},
  author={Letrouit, Cyril},
  journal={Notes du cours Peccot},
  volume={2025},
  year={2025}
}

@article{letrouit2024,
  title={Gluing methods for quantitative stability of optimal transport maps},
  author={Letrouit, Cyril and M{\'e}rigot, Quentin},
  journal={arXiv preprint arXiv:2411.04908},
  year={2024}
}

@article{brenier1991polar,
  title={Polar factorization and monotone rearrangement of vector-valued functions},
  author={Brenier, Yann},
  journal={Communications on pure and applied mathematics},
  volume={44},
  number={4},
  pages={375--417},
  year={1991},
  publisher={Wiley Online Library}
}

@article{cordero1999transport,
  title={Sur le transport de mesures p{\'e}riodiques},
  author={Cordero-Erausquin, Dario},
  journal={Comptes Rendus de l'Acad{\'e}mie des Sciences-Series I-Mathematics},
  volume={329},
  number={3},
  pages={199--202},
  year={1999},
  publisher={Elsevier}
}

@book{protter1984,
  title={Maximum Principles in Differential Equations},
  author={Protter, Murray H. and Weinberger, Hans F.},
  year={1984},
  publisher={Springer-Verlag},
  address={New York}
}

@article{pinn,
  title={Physics-informed neural networks: A deep learning framework for solving forward and inverse problems involving nonlinear partial differential equations},
  author={Raissi, Maziar and Perdikaris, Paris and Karniadakis, George E},
  journal={Journal of Computational physics},
  volume={378},
  pages={686--707},
  year={2019},
  publisher={Elsevier}
}

@book{carmona2018probabilistic,
  title={Probabilistic theory of mean field games with applications I-II},
  author={Carmona, Ren{\'e} and Delarue, Fran{\c{c}}ois and others},
  volume={3},
  year={2018},
  publisher={Springer}
}

@article{arakelyan2026mfg,
  title={A Mean-Field Game Model For Large-Scale Attrition in Attacker-Defender Systems},
  author={Arakelyan, Avetik and Bakaryan, Tigran and Alaverdyan, Davit and Hovakimyan, Naira and Kaminer, Isaac},
  journal={arXiv preprint arXiv:2604.02101},
  year={2026}
}

@article{PhamWei2017,
  author  = {Pham, H. and Wei, X.},
  title   = {Dynamic programming for optimal control of stochastic {McKean}--{Vlasov} dynamics},
  journal = {SIAM Journal on Control and Optimization},
  volume  = {55},
  number  = {2},
  pages   = {1069--1101},
  year    = {2017},
  doi     = {10.1137/16M1071390}
}

@article{Cardaliaguet2016,
  author  = {Cardaliaguet, P.},
  title   = {Mean field games: the master equation and the mean field limit},
  journal = {S\'eminaire Laurent Schwartz --- \'Equations aux d\'eriv\'ees partielles et applications},
  volume  = {2015--2016},
  pages   = {1--10},
  year    = {2016},
  note    = {Expos\'e no.~17},
  doi     = {10.5802/slsedp.99}
}

@article{Graber2025,
  author  = {Graber, P. J.},
  title   = {Remarks on potential mean field games},
  journal = {Research in the Mathematical Sciences},
  volume  = {12},
  pages   = {Article 13},
  year    = {2025},
  doi     = {10.1007/s40687-024-00494-3}
}

@article{Santambrogio2017,
  author  = {Santambrogio, F.},
  title   = {Euclidean, metric, and {Wasserstein} gradient flows: an overview},
  journal = {Bulletin of Mathematical Sciences},
  volume  = {7},
  pages   = {87--154},
  year    = {2017},
  doi     = {10.1007/s13373-017-0101-1}
}

@book{ladyzhenskaya1968linear,
  title={Linear and quasi-linear equations of parabolic type},
  author={Ladyzhenskaya, Olga Aleksandrovna and Solonnikov, Vsevolod Alekseevich and Uraltseva, Nina Nikolaevna},
  volume={23},
  year={1968},
  publisher={American Mathematical Soc.}
}
\end{document}